\documentclass[opre,nonblindrev]{informs3}

\usepackage{endnotes}
\let\footnote=\endnote
\let\enotesize=\normalsize
\def\notesname{Endnotes}%
\def\makeenmark{\hbox to1.275em{\theenmark.\enskip\hss}}
\def\enoteformat{\rightskip0pt\leftskip0pt\parindent=1.275em
  \leavevmode\llap{\makeenmark}}

\usepackage{subcaption}
\usepackage{algpseudocode}
\usepackage{algorithm}
\usepackage{bm}
\usepackage{epstopdf}
\graphicspath{{figures/}}

\newcommand{\prth}[1]{\left(#1\right)}

\newcommand{\norm}[1]{\left \Vert #1 \right\Vert}

\newcommand{\brac}[1]{\left[ #1 \right]}
\newcommand{\R}{\mathbb{R}}

\usepackage{verbatim}
\usepackage{color}
\usepackage{textcomp} 
\usepackage{rotating}
\usepackage{graphicx}
\usepackage{caption}
\usepackage{subcaption}
\usepackage{paralist} 
\usepackage{booktabs}

\usepackage{epsfig}

\usepackage{bbm}
\usepackage{indentfirst}
\usepackage{float}

\usepackage{natbib}
 \bibpunct[, ]{(}{)}{,}{a}{}{,}%
 \def\bibfont{\small}%
\TheoremsNumberedThrough     
\ECRepeatTheorems

\EquationsNumberedThrough    

\begin{document}


 \RUNAUTHOR{Goeva et al.}


\TITLE{Optimization-based Calibration of Simulation Input Models}

\ARTICLEAUTHORS{%
\AUTHOR{Aleksandrina Goeva}
\AFF{Broad Institute, Cambridge, MA 02142, USA. \EMAIL{agoeva@broadinstitute.org}} 
\AUTHOR{Henry Lam}
\AFF{Department of Industrial Engineering and Operations Research, Columbia University, New York, NY 10027, USA. \EMAIL{henry.lam@columbia.edu}}
\AUTHOR{Huajie Qian}
\AFF{Department of Mathematics, University of Michigan, Ann Arbor, MI 48109, USA. \EMAIL{hqian@umich.edu}}
\AUTHOR{Bo Zhang}
\AFF{IBM Research AI, Yorktown Heights, NY 10598, USA. \EMAIL{zhangbo@us.ibm.com}}
} 

\ABSTRACT{%
Studies on simulation input uncertainty often built on the availability of input data. In this paper, we investigate an inverse problem where, given only the availability of output data, we nonparametrically calibrate the input models and other related performance measures of interest. We propose an optimization-based framework to compute statistically valid bounds on input quantities. The framework utilizes constraints that connect the statistical information of the real-world outputs with the input-output relation via a simulable map. We analyze the statistical guarantees of this approach from the view of data-driven robust optimization, and show how the guarantees relate to the function complexity of the constraints arising in our framework. We investigate an iterative procedure based on a stochastic quadratic penalty method to approximately solve the resulting optimization. We conduct numerical experiments to demonstrate our performance in bounding the input models and related quantities.
}%



\maketitle

%

\section{Introduction}
Stochastic simulation takes in input models and generates random outputs for subsequent performance analyses. The accuracy of these input model assumptions is critical to the analyses' credibility. In the conventional premise in studying stochastic simulation, these input models are conferred either through physical implication or expert opinions, or observable via input data. In this paper, we answer a converse question: Given only \emph{output} data from a stochastic system, can one infer about the input model?

The main motivation for asking this question is that, in many situations, a simulation modeler plainly may not have the luxury of direct data or knowledge about the input. The only way to gain such knowledge could be data from other sources that are at the output level. For instance, one of the authors has experienced such complication when building a simulation model for a contract fulfillment center, where service agents work on a variety of processing tasks and, despite the abundant transaction data stored in the center's IT system, there is no record on the start, completion, or service times spent by each agent on each particular task. Similarly, in clinic operations, patients often receive service in multiple phases such as initial checkup, medical tests and doctor's consultation. Patients' check-in and check-out times could be accurately noted, but the ``service" times provided by the medical staff could very well be unrecorded. Clearly, these service time distributions are needed to build a simulation model, if an analyst wants to use the model for sensitivity analysis or system optimization purposes.

The problem of inferring an input model from output data is sometimes known as \emph{model calibration}. In the simulation literature, this is often treated as a refinement process that occurs together with iterative comparisons between simulation reports and real-world output data (a task known as \emph{model validation}; \cite{sargent2005verification,kleijnen1995verification}). If simulation reports differ significantly from output data, the simulation model is re-calibrated (which can involve both the input distributions and system specifications), re-compared, and the process is iterated. Suggested approaches to compare simulation with real-world data include conducting statistical tests such as two-sample mean-difference tests (\cite{balci1982some}) and the Schruben-Turing test (\cite{schruben1980establishing}). Beyond that, inferring input from output seems to be an important problem that has not been widely discussed in the stochastic simulation literature (\cite{nelson2016some}). 

The setting we consider can be briefly described as follows. We assume an input model is missing and make no particular assumptions on the form of its probability distribution. We assume, however, that a certain output random variable from a well-specified system is observable with some data. Our task is to nonparametrically infer the input distribution, or other quantities related to this input distribution (e.g., a second output measure driven by the same input distribution). One distinction between our setting and model calibration in other literature (e.g., computer experiments) is the intrinsic probabilistic structure of the system. Namely, the input and the output in stochastic simulation are represented as probability distributions, or in other words, the relation that links the observed and the to-be-calibrated objects is a (simulable) map between the spaces of distributions. Our calibration method will be designed to take such a relation into account.

Specifically, we study an optimization-based framework for model calibration, where the optimization, on a high level, entails an objective function associated with the ``input" and constraints associated with the ``output". The decision variable in this optimization is the unknown input distribution. The constraints comprise a confidence region on the  the output distribution that is compiled from the observed output statistics. By expressing the region in terms of the input distributions via the simulable map, the optimization objective, which is set to be some target input quantity, will then give rise to statistically valid confidence bounds on this target. Advantageously, this approach leads to valid bounds even if the input model is \emph{non-identifiable}, i.e., there exist more than one input model that give rise to the same observable output pattern, which may occur since the simulable map is typically highly complicated. The tightness of the bounds in turn depends on the degree of non-identifiability (which also leads to a notion of \emph{identifiability gap} that we will discuss). The idea of utilizing a confidence region as the constraint is inspired by distributionally robust optimization (DRO). However, in the conventional DRO literature, the constraints (often called collectively as the uncertainty set or the ambiguity set) are constructed based on direct observation of data. On the other hand, our constraints here serve as a tool to integrate the input-output relation, in addition to the output-level statistical noise, to effectively calibrate the input model. This leads to several new methodological challenges and solution approaches.

Under this general framework, we propose a concrete optimization formulation that balances statistical validity and the required computational efforts. Specifically, we use a nonparametric statistic, namely the Kolmogorov-Smirnov (KS) statistic, to construct the output-level confidence region. This formulation has the strengths of being statistically consistent (implied by the KS statistic) and expressible as expectation-type constraints that can be effectively solved by our subsequent algorithms. It also has an interesting additional benefit in terms of controlling the dimension of the optimization. Because of computational capacity, the decision variable, which is the unknown input distribution and potentially infinite-dimensional, needs to be suitably discretized by randomly generating a finite number of support points. A consistent statistic typically induces a large number of constraints, and one may need to use a large number of support points to retain the discretization error. However, as will be seen, it turns out that the KS constraints allow us to use a moderate support size without compromising the asymptotic statistical guarantees, thanks to their low complexity as measured by the so-called bracketing number in the empirical process theory. This thus leads us to an optimization problem with both a controllable number of decision variables and statistical validity.

Next, due to the sophisticated input-output map, the optimization programs generally involve non-convex stochastic (i.e., simulation-based) constraints. We propose and analyze a stochastic quadratic penalty method, by adding a growing penalty on the squared constraint violation. This method borrows from the quadratic penalty method used in deterministic nonlinear programming. However, while the deterministic version suggests solving a nonlinear program at each particular value of the penalty coefficient and letting the coefficient grows, the stochastic method we analyze involves a stochastic approximation (SA) that runs updates of the solution, slack variables and the penalty coefficient simultaneously. This is motivated from the typical challenge of finding good stopping times for SA, which are needed for each SA run at each penalty coefficient value if one were to mimic the deterministic procedure. Simultaneous updates of all the quantities, however, only need one SA run. We analyze the convergence guarantee of this algorithm and provide guidance on the step sizes of all the constituent updates. Our SA update uses a mirror descent stochastic approximation (MDSA) (\cite{nemirovski2009robust}), in particular the entropic descent (\cite{beck2003mirror}). 



 

The remainder of the paper is organized as follows. Section \ref{sec:lit} reviews the related literature. Section \ref{sec:formulation} introduces the problem setting and presents our general optimization-based framework. Section \ref{sec:KS} refines our framework with the KS-based formulations and demonstrates the statistical guarantees. Section \ref{sec:procedure} presents and analyzes our optimization algorithm. Section \ref{sec:numerics} reports numerical results. Section \ref{sec:discussion} concludes. The Appendix contains all the proofs. 

\section{Related literature}\label{sec:lit}
We organize the literature review in two aspects, one related to the model calibration problem, and one related to our optimization approach.

\subsection{Literature Related to Our Problem Setting}
Input modeling and uncertainty quantification in the stochastic simulation focus mostly on the input level. \cite{barton2012tutorial} and \cite{song2014input}, e.g., review some major methods in quantifying the statistical errors from finite input data. These methods include the delta or two-point method (\cite{cheng1998two,cheng2004calculation}), Bayesian methodology and model averaging (\cite{chick2001input,zouaoui2004accounting}) and resampling methods (\cite{barton2001resampling,barton2013quantifying}). Our problem is more related to model calibration. In the simulation literature, this is often considered together with model validation (\cite{sargent2005verification,kleijnen1995verification}). Conventional approaches compare simulation data with real-world historical output data according to statistical or Turing tests (\cite{balci1982some,schruben1980establishing}), conduct re-calibration, and repeat the process until the data are successfully validated (\cite{banks2000dm,kelton2000simulation}).

The model calibration problem is also known as the \emph{inverse problem} (\cite{tarantola2005inverse}) in the literature of other fields. It generally refers to the identification of parameters or functions that can only be inferred from transformed outputs. In the context where the parameters are probability distributions, \cite{kraan2005probabilistic} demonstrates theoretically the characterization of a distribution that leads to the smallest relative entropy with a reference measure, and proposes an entropy maximization to calibrate the distribution from output data. Our work relates to \cite{kraan2005probabilistic} as we also utilize a probabilistic input-output map, but we focus on maps that are evaluable only by simulation, and aim to compute confidence bounds on the true distribution instead of attempting to recover the maximum entropy distribution. 

The inverse problem also appeared in many other contexts. In signal processing, the linear inverse problem (e.g., \cite{csiszar1991least,donoho1992maximum}) reconstructs signals from measurements of linear transformations. Common approaches consist of least-square minimization and the use of penalty such as the entropy. In computer experiments (\cite{santner2013design}), surrogate models built on complex physical laws require the calibration of physical parameters. Such models have wide scientific applications such as weather prediction, oceanography, nuclear physics, and acoustics (e.g., \cite{wunsch1996ocean,shirangi2014history}). Bayesian and Gaussian process methodologies are commonly used (e.g., \cite{kennedy2001bayesian,currin1991bayesian}). We point out that Bayesian methods could be a potential alternative to the approach considered in this paper, but because of the nature of discrete-event systems, one might need to resort to sophisticated techniques such as approximate Bayesian computation (\cite{marjoram2003markov}). Other related literature include experimental design to optimize inference for input parameters (e.g., \cite{chick2002simulation}) and calibrating financial option prices (e.g., \cite{avellaneda2001weighted,glasserman2005large}).

Also related to our work is the body of research on inference problems in the context of queueing systems. The first stream, similar to our paper, aims at inferring the constituent probability distributions of a queueing model based on its output data, e.g., queue length or waiting time data, collected either continuously or at discrete time points. This stream of papers focuses on systems whose structures allow closed-form analyses or are amenable to analytic approximations via, for instance, the diffusion limit. The majority of them assume that the inferred distribution(s) comes from a parametric family and use maximum likelihood estimators (\cite{basawa1996maximum,pickands1997estimation,basawa2008parameter,fearnhead2004filtering,wang2006maximum,ross2007estimation,heckmuller2009reconstructing,whitt2012fitting}). Others work on nonparametric inference by exploiting specific queueing system structures (\cite{bingham1999non,hall2004nonparametric,moulines2007,feng2014estimating}). A related stream of literature studies point process approximation (see Section 4.7 of \cite{cooper1972introduction}, \cite{whitt1981approximating,whitt1982approximating}, and the references therein), based on a parametric approach and is motivated from traffic pattern modeling in communication networks. Finally, there are also a number of studies inspired by the ``queue inference engine" by \cite{larson1990queue}. But, instead of inferring the input models, many of these studies use transaction data to estimate the performance of a queueing system directly and hence do not take on the form of an inverse problem (see \cite{mandelbaum1998estimating} for a good survey of the earlier literature and \cite{frey2010queue} and its references for more recent progress). Several papers estimate both the queueing operational performance and the constituent input models (e.g., \cite{daley1998moment,kim2008new,park2011analysis}), and can be considered to belong to both this stream and the aforementioned first stream of literature. 


\subsection{Literature Related to Our Methodology}
Our formulation uses ideas from robust optimization (e.g., \cite{bertsimas2011theory,ben2009robust}), which studies optimization under uncertain parameters and suggests to obtain decisions that optimize the worst-case scenarios, subject to a set of constraints on the belief/uncertainty that is often known as the ambiguity set or the uncertainty set. Of particular relevance to us is the setting of distributionally robust optimization (DRO), where the uncertainty is on the probability distribution in a stochastic optimization problem. This approach has been applied in many disciplines such as stochastic control  (e.g., \cite{pjd00,doi:10.1287/moor.1120.0540,iyengar2005robust}), economics (\cite{hansen2008robustness}), finance (\cite{gx12a}), queueing control (\cite{jls10}) and dynamic pricing (\cite{ls07}). Its connection to machine learning and statistics has also been recently investigated (\cite{blanchet2016robust,shafieezadeh2015distributionally}). In the DRO literature, common choices of the uncertainty set are based on moments (\cite{delage2010distributionally,goh2010distributionally,wiesemann2014distributionally,bertsimas2005optimal,smith95,bertsimas2007semidefinite}), distances from nominal distributions (\cite{ben2013robust,bayraksan2015data,blanchet2016quantifying,esfahani2015data,gao2016distributionally}), and shape conditions (\cite{popescu2005semidefinite,lam2017tail,li2016ambiguous,hanasusanto2017ambiguous}). The literature of data-driven DRO further addresses the question of calibrating these sets, using for instance confidence regions or hypothesis testing (\cite{bertsimas2014robust}), empirical likelihood (\cite{lam2017empirical,duchi2016statistics,lam2016recovering}) and the related Wasserstein profile function (\cite{blanchet2016robust}), and Bayesian perspectives (\cite{gupta2015near}). 

For DRO in the simulation context, \cite{hu2012robust} studies the computation of robust bounds under Gaussian model assumptions, \cite{gx12b,glasserman2016bounding} study distance-based constraints to address model risks, \cite{lam2013robust,lam2016serial} study asymptotic approximations of related formulations, and \cite{ghosh2015robust} studies formulations and solution techniques for DRO in quantifying simulation input uncertainty. \cite{fan2013robust,ryzhov2012ranking} study the use of robust optimization in simulation-based decision-making. Our framework in particular follows the concept in using confidence region such that the uncertainty set covers the true distribution with high probability. However, it also involves the simulation map between input and output that serves as the key in our model calibration goal.


 
Our optimization procedure builds on the quadratic penalty method (\cite{bertsekas1999nonlinear}), which is a deterministic nonlinear programming technique that reformulates the constraints as squared penalty and sequentially tunes the penalty coefficient to approach optimality. Different from the deterministic technique, our procedure in solving the stochastic quadratic penalty formulation sequentially update the penalty parameter simultaneously together with the solution and slack variables. This involves a specialized version of MDSA proposed by \cite{nemirovski2009robust}. \cite{nemirovski2009robust} analyzed convergence guarantees on convex programs with stochastic objectives. \cite{lanzhou2017,NIPS2017_6741} investigated convex stochastic constraints, and \cite{ghadimi2013stochastic,ghadimi2015accelerated,dang2015stochastic,ghadimi2016mini} studied related schemes for nonconvex and nonsmooth objectives. \cite{wang2008stochastic} introduced a quadratic penalty method for stochastic objectives with deterministic constraints. The particular scheme of MDSA we consider uses entropic penalty, and is known as the entropic descent algorithm (\cite{beck2003mirror}).


\section{Proposed Framework}\label{sec:formulation}
Consider a generic input variate $X$ with an input probability distribution $P_X$. We let $\mathbf X=(X_1,\ldots,X_{T})$, where $X_t\in\mathcal X$, be an i.i.d. sequence of input variates each distributed under $P_X$ over a time horizon $T$. We denote the function $h(\cdot)\in\mathbb R$ as the system logic from the input sequence $\mathbf X$ to the output $h(\mathbf X)$. We assume that $h$ is completely specified and is computable, even though it may not be writable in closed-form, i.e. we can evaluate the output given $\mathbf X$. For example, $\mathbf X$ can denote the sequence of interarrival or service times for the customers in a queue, and $h(\mathbf X)$ is an average queue length seen by the $T$ customers. Note that we can work in a more general framework where $h$ depends on both $\mathbf X$ and other independent input sequences, denoted collectively as $\mathbf W$, that possess known or observable distributions. In other words, we can have $h(\mathbf X,\mathbf W)$ as the output. Our developments can readily handle this case, but for expositional convenience we will assume the absence of these auxiliary input sequences most of the time, and will indicate the modifications of our developments in handling them at various suitable places.

Consider the situation that only $h(\mathbf X)$ can be observed via data. Let $D=\{y_1,\ldots,y_n\}$ be $n$ observations of $h(\mathbf X)$. Our task is to calibrate some quantities related to $P_X$, which we call $\psi(P_X)$. Two types of target quantities we will consider are:
\begin{enumerate}
\item Restricting $X$ to real value, we consider the distribution function of $P_X$, denoted $F_X(x)$, where $x$ can take a range of values. Note that, obviously, $F_X(x)=E_{P_X}[I(X\leq x)]$ where $E_{P_X}[\cdot]$ denotes the expectation with respect to $P_X$ and $I(\cdot)$ denotes the indicator function.\label{case1}
\item We consider a performance measure $E_{P_X}[g(\mathbf X)]$ where $E_{P_X}[\cdot]$ here denotes the expectation with respect to the product measure induced by the i.i.d. sequence $\mathbf X=(X_1,\ldots,X_{S})$ over a time horizon $S$. The function $g(\mathbf X)$ can denote another output of interest different from $h(\mathbf X)$ that is unobservable, and requires information about $\mathbf X$. This case includes the first target quantity above (by choosing $g(\mathbf X)=I(X_1\leq x)$ when $X$ is real-valued), as well as other statistics of $X$ such as power moments (by choosing $g(\mathbf X)=X_1^k$ for some $k$).\label{case2}
\end{enumerate}

To describe our framework, we denote $P_Y=P_{h(\mathbf X)}$ as the probability distribution of the output $Y=h(\mathbf X)$. Since $P_Y$ is completely identified by $P_X$, we can view $P_Y$ as a transformation of $P_X$, i.e., $P_Y=\gamma(P_X)$ for some map $\gamma$ between probability distributions. We denote $\mathcal P_X$ and $\mathcal P_Y$ as the spaces of all possible input and output distributions respectively. 

On an abstract level, we use the optimization formulations
\begin{equation}
\begin{array}{ll}
\max&\psi(P_X)\\
\text{subject to}&P_Y\in\mathcal U
\end{array}\label{RO max}
\end{equation}
and
\begin{equation}
\begin{array}{ll}
\min&\psi(P_X)\\
\text{subject to}&P_Y\in\mathcal U
\end{array}\label{RO min}
\end{equation}
where the decision variable is the unknown $P_X\in\mathcal P_X$, and $\mathcal U\subset\mathcal P_Y$ is an ``uncertainty set" that covers a set of possibilities for $P_Y$. The objective function $\psi(P_X)$ refers to either $F_X(x)$ in case \ref{case1} or $E_{P_X}[g(\mathbf X)]$ in case \ref{case2} above. 

An important element in formulations \eqref{RO max} and \eqref{RO min} is that the constraints represented by $\mathcal U$ are cast on the output level. Since we have available output data, $\mathcal U$ can be constructed using these observations in a statistically valid manner (e.g., by using the confidence region on the output statistic). By expressing $P_Y=\gamma(P_X)$, the region $\mathcal U$ can be viewed as a region on $P_X$, given by $\{P_X\in\mathcal P_X:\gamma(P_X)\in\mathcal U\}$. The following result summarizes the confidence guarantee for the optimal values of \eqref{RO max} and \eqref{RO min} in bounding $\psi(P_X)$ when $\mathcal U$ is chosen suitably:
\begin{proposition}
Let $P_X^0\in\mathcal P_X$ and $P_Y^0\in\mathcal P_Y$ be the true input and output distributions. Suppose $\mathcal U$ is a $(1-\alpha)$-level confidence region for $P_Y^0$, i.e., 
\begin{equation}
\mathbb P_D(P_Y^0\in\mathcal U)\geq1-\alpha\label{confidence setup}
\end{equation}
where $\mathbb P_D(\cdot)$ denotes the probability with respect to the data $D$. Let $\overline Z$ and $\underline Z$ be the optimal values of \eqref{RO max} and \eqref{RO min} respectively. Then we have
$$\mathbb P_D(\underline Z\leq\psi(P_X^0)\leq\overline Z)\geq1-\alpha$$
Similar statements hold if the confidence is approximate, i.e., if
$$\liminf_{n\to\infty}\mathbb P_D(P_Y^0\in\mathcal U)\geq1-\alpha$$
then
$$\liminf_{n\to\infty}\mathbb P_D(\underline Z\leq\psi(P_X^0)\leq\overline Z)\geq1-\alpha$$
\label{basic guarantee}
\end{proposition}

It is worth pointing out that the same guarantee holds, without any statistical adjustment, if one solves \eqref{RO max} and \eqref{RO min} simultaneously for different $\psi(\cdot)$, say $\psi_l(\cdot),l=1,\ldots,L$, i.e., supposing that \eqref{confidence setup} holds, then the confidence statement
$$\mathbb P_D(\underline Z_l\leq\psi_l(P_X^0)\leq\overline Z_l,\ l=1,\ldots,L)\geq1-\alpha$$
holds, so does a similar statement for the limiting counterpart. We provide this extended version of Proposition \ref{basic guarantee} in the appendix (Proposition \ref{extended guarantee}). This allows us to obtain bounds for multiple quantities about the input model at the same time. Note that, in conventional statistical methods, simultaneous estimation like this sort often requires Bonferroni correction or more advanced techniques, but these are not needed in our approach. 

We mention an important feature of our framework related to the issue of \emph{non-identifiability} (e.g., \cite{tarantola2005inverse}). When there are more than one input model $P_X$ that leads to the same output distribution, it is statistically impossible to recover exactly the true $P_X$, and methods that attempt to do so may result in ill-posed problems. Our framework, however, gets around this issue by focusing on computing bounds instead of full model recovery. Even though $P_X$ can be non-identifiable, our optimization always produces valid bounds for it. One special case of interest is when we use $\mathcal U=\{P_Y^0\}$, i.e., the true output distribution is exactly known. In this case, \eqref{RO max} and \eqref{RO min} will provide the best bounds for $\psi(P_X)$ given the output. If $\underline Z<\overline Z$, then $P_X$ cannot be exactly identified, implying an issue of non-identifiabilty, but our outputs would still be valid. In fact, the difference $\overline Z-\underline Z$ can be viewed as an \emph{identifiability gap} with respect to $\psi$.

\section{Kolmogorov-Smirnov-based Constraints}\label{sec:KS}
We will now choose a specific $\mathcal U$ that is statistically consistent on the output level, i.e., $\mathcal U$ shrinks to $\{P_Y^0\}$ as $n\to\infty$ (in a suitable sense). In particular, we use $\mathcal U$ implied by the Kolmogorov-Smirnov (KS) statistic, and discuss how this choice enjoys benefits balancing statistical consistency and computation.

\subsection{Statistical Confidence Guarantee}\label{sec:KS confidence}
It is known that the empirical distribution for continuous i.i.d. data $D$, denoted $\hat F_Y(y)$, satisfies $\sqrt n\|\hat F_Y-F_Y^0\|_\infty\Rightarrow\sup_{u\in[0,1]}BB(u)$ where $F_Y^0$ is the true distribution function of $Y$, $\|\cdot\|_\infty$ denotes the sup norm over $\mathbb R$, $BB(\cdot)$ is a standard Brownian bridge, and $\Rightarrow$ denotes weak convergence. This implies that the KS-statistic $\sqrt n\|\hat F_Y-F_Y^0\|_\infty$ satisfies
$$\lim_{n\to\infty}P\left(\|\hat F_Y-F_Y^0\|_\infty\leq\frac{q_{1-\alpha}}{\sqrt n}\right)=1-\alpha$$
where $q_{1-\alpha}$ is the $(1-\alpha)$-quantile of $\sup_{u\in[0,1]}BB(u)$. Therefore, setting 
\begin{equation}
\mathcal U=\left\{P_Y\in\mathcal P_Y:\|F_Y-\hat F_Y\|_\infty\leq\frac{q_{1-\alpha}}{\sqrt n}\right\}\label{uncertainty set}
\end{equation}
ensures that \eqref{confidence setup} holds and subsequently the conclusion in Proposition \ref{basic guarantee}. As $n$ increases, the size of \eqref{uncertainty set} shrinks to zero.

The following result states precisely the implication of this construction, and moreover, describes how this leads to a more tractable optimization formulation:
\begin{theorem}
Let $\overline Z$ and $\underline Z$ be the optimal values of the optimization programs
\begin{equation}
\begin{array}{ll}
\max&\psi(P_X)\\
\text{subject to}&\|F_Y-\hat F_Y\|_\infty\leq\frac{q_{1-\alpha}}{\sqrt n}\\
&P_X\in\mathcal P_X
\end{array}\label{KS RO max direct}
\end{equation}
and
\begin{equation}
\begin{array}{ll}
\min&\psi(P_X)\\
\text{subject to}&\|F_Y-\hat F_Y\|_\infty\leq\frac{q_{1-\alpha}}{\sqrt n}\\
&P_X\in\mathcal P_X
\end{array}\label{KS RO min direct}
\end{equation}
where $q_{1-\alpha}$ is the $(1-\alpha)$-quantile of $\sup_{u\in[0,1]}BB(u)$, and $\hat F_Y$ is the empirical distribution of i.i.d. output data. Supposing the true output distribution is continuous, we have
\begin{equation}
\liminf_{n\to\infty}\mathbb P_D(\underline Z\leq\psi(P_X^0)\leq\overline Z)\geq1-\alpha\label{KS confidence guarantee}
\end{equation}
where $P_X^0$ is the true distribution of the input variate $X$. Moreover, \eqref{KS RO max direct} and \eqref{KS RO min direct} are equivalent to

\begin{equation}
\begin{array}{ll}
\max&\psi(P_X)\\
\text{subject to}&\hat F_Y(y_j+)-\frac{q_{1-\alpha}}{\sqrt n}\leq E_{P_X}[I(h(\mathbf X)\leq y_j)]\leq\hat F_Y(y_j-)+\frac{q_{1-\alpha}}{\sqrt n},j=1,\ldots,n\\
&P_X\in\mathcal P_X
\end{array}\label{KS RO max}
\end{equation}
and
\begin{equation}
\begin{array}{ll}
\min&\psi(P_X)\\
\text{subject to}&\hat F_Y(y_j+)-\frac{q_{1-\alpha}}{\sqrt n}\leq E_{P_X}[I(h(\mathbf X)\leq y_j)]\leq\hat F_Y(y_j-)+\frac{q_{1-\alpha}}{\sqrt n},j=1,\ldots,n\\
&P_X\in\mathcal P_X
\end{array}\label{KS RO min}
\end{equation}
respectively, where $\hat F_Y(y_j+)$ and $\hat F_Y(y_j-)$ refer to the right- and left-limits of the empirical distributions $\hat F_Y$ at $y_j$, and $E_{P_X}[\cdot]$ denotes the expectation taken with respect to the $T$-fold product measure of $P_X$. \label{KS guarantee}
\end{theorem}

A merit of using the depicted KS-based uncertainty set, seen by Theorem \ref{KS guarantee}, is that it can be reformulated into linear constraints in terms of the expectations $E_{P_X}[\cdot]$ of certain ``moments" of $h(\mathbf X)$. These constraints constitute precisely $n$ interval-type conditions, and the moment functions are the indicator functions of $h(\mathbf X)$ falling under the thresholds $y_j$'s. The derivation leading to the reformulation result in Theorem \ref{KS guarantee} has been used conventionally in computing the KS-statistic. Similar reformulations have also appeared in recent work in approximating stochastic optimization via robust optimization (\cite{bertsimas2014robust}).

The asymptotic of the KS-statistic is more complicated if the output distribution is discrete (this happens if the outputs we look at are for instance the queue length). In such cases, the  critical values are generally smaller than those for the continuous distribution (\cite{lehmann2006testing}). Consequently, using $q_{1-\alpha}/\sqrt n$ to calibrate the size of the uncertainty set as in \eqref{uncertainty set} is still valid, but could be conservative, i.e., we still have $P\left(\|\hat F_Y-F_Y^0\|_\infty\leq\frac{q_{1-\alpha}}{\sqrt n}\right)$ asymptotically at least $1-\alpha$, but possibly strictly higher. As a remedy, one can use bootstrapping to calibrate the size of a tighter set. Moreover, the constraint of the form $\|\hat F_Y-F_Y^0\|_\infty\leq q$  will now be written as
\begin{equation}
\hat F_Y(w_j)-q\leq E_{P_X}[I(h(\mathbf X)\leq w_j)]\leq\hat F_Y(w_j)+q,j=1,\ldots,K\label{discrete KS}
\end{equation}
where $w_j,j=1,\ldots$ are the ordered support points of $Y$, with $K=\min\{j:\hat F_Y(w_j)=1\}$. These are the points where jumps occur (and the constraints put on the first $K$ of them automatically ensure the rest). If the support size is small, an alternative is to impose constraints on each probability mass, i.e.,
\begin{equation}
\hat P(Y=w_j)-q\leq E_{P_X}[I(h(\mathbf X)=w_j)]\leq\hat P(Y=w_j)+q,j=1,\ldots,K\label{discrete moment}
\end{equation}
where $\hat P(Y=w_j)$ is the observed proportions of $Y$ being $w_j$, and $q$ can be calibrated by a standard binomial quantile and the Bonferroni correction.

The KS-statistic has several advantages over other types of uncertainty sets in our considered settings. Alternatives like $\chi^2$ goodness-of-fit tests could be used, but the resulting formulations would not come as handy when expressed in terms of $P_X$ or $h(\mathbf X)$, which would affect the efficiency of the gradient estimator that we will discuss in Section \ref{sec:gradient}. Another advantage of using KS-statistic relates to the statistical property of a discretization that is needed to feed into an implementable optimization procedure, which we shall discuss next. 





\subsection{Randomizing the Decision Space}
\label{sec:discretization}

Note that optimization programs \eqref{KS RO max} and \eqref{KS RO min} involve decision variable $P_X$ that is potentially infinite-dimensional, e.g., when $X$ is a continuous variable. This can cause algorithmic and storage issues. One could appropriately discretize the decision variable by randomly sampling a finite set of support points on $\mathcal X$. Once these support points are realized, the optimization is imposed on the probability weights on these points, or in other words on a discrete input distribution.

Our next result shows that as the support points are generated from a suitably chosen distribution, and the number of these points grows at an appropriate rate relative to the output data size, the discretized KS-implied optimization will retain the confidence guarantee as the original formulation:






\begin{theorem}
Suppose we sample $\{z_i\}_{i=1,\ldots,m}$ in the space $\mathcal X$ from a distribution $Q$. Suppose that $P_X^0$, the true distribution of $X$, is absolutely continuous with respect to $Q$ and $\|dP_X^0/dQ\|_\infty\leq C$ for some $C>0$, where $dP_X^0/dQ$ is the likelihood ratio calculated from the Radon-Nikodym derivative of $P_X^0$ with respect to $Q$, and $\|\cdot\|_\infty$ denotes the essential supremum. Using the notations as in Theorem \ref{KS guarantee}, we solve 
\begin{equation}
\begin{array}{ll}
\max&\psi(P_X)\\
\text{subject to}&\hat F_Y(y_j+)-\frac{q_{1-\alpha}}{\sqrt n}\leq E_{P_X}[I(h(\mathbf X)\leq y_j)]\leq\hat F_Y(y_j-)+\frac{q_{1-\alpha}}{\sqrt n},j=1,\ldots,n\\
&P_X\in\hat{\mathcal P}_X
\end{array}\label{KS RO max discretized}
\end{equation}
and
\begin{equation}
\begin{array}{ll}
\min&\psi(P_X)\\
\text{subject to}&\hat F_Y(y_j+)-\frac{q_{1-\alpha}}{\sqrt n}\leq E_{P_X}[I(h(\mathbf X)\leq y_j)]\leq\hat F_Y(y_j-)+\frac{q_{1-\alpha}}{\sqrt n},j=1,\ldots,n\\
&P_X\in\hat{\mathcal P}_X
\end{array}\label{KS RO min discretized}
\end{equation}
where $\hat{\mathcal P}_X$ denotes the set of distributions with support points $\{z_i\}_{i=1,\ldots,m}$. Let $\hat{\overline Z}$ and $\hat{\underline Z}$ be the optimal values of \eqref{KS RO max discretized} and \eqref{KS RO min discretized}. 

Denote $\mathbb P$ as the probability taken with respect to both the output data and the support generation for $X$. Suppose that $\psi(P_X)$ takes the form $E_{P_X}[g(\mathbf X)]$ (which subsumes both types of target measures discussed in Section \ref{sec:formulation}) where $E_{P_X^0}[g(X_{i_1},\ldots,X_{i_T})^2]<\infty$ for any $1\leq i_1,\ldots,i_T\leq T$. Also suppose that the true output distribution is continuous and that $\mathbb P(\text{for any\ }P_X\in\hat{\mathcal P}_X, \text{supp}(\gamma(P_X))\cap\{y_j\}_{j=1,\ldots,n}\neq\emptyset)=0$ where $\text{supp}(\gamma(P_X))$ denotes the support of the distribution $\gamma(P_X)$. Then, we have
$$\liminf_{n\to\infty,m/n\to\infty}\mathbb P\left(\hat{\underline Z}+O_p\left(\frac{1}{\sqrt m}\right)\leq\psi(P_X^0)\leq\hat{\overline Z}+O_p\left(\frac{1}{\sqrt m}\right)\right)\geq1-\alpha$$
\label{main guarantee}
\end{theorem}

The error term $O_p(1/\sqrt m)$ represents a random variable of stochastic order $1/\sqrt m$, i.e., $a_m=O_p(1/\sqrt m)$ if for any $\epsilon>0$, there exists $M,N>0$ such that $P(|\sqrt m a_m|\leq N)>1-\epsilon$ for $m>M$. 

Theorem \ref{main guarantee} guarantees that by solving the finite-dimensional optimization problems \eqref{KS RO max discretized} and \eqref{KS RO min discretized}, we obtain confidence bounds for the true quantity of interest $\psi(P_X^0)$, up to an error of order $O_p(1/\sqrt m)$. Note that the conclusion holds with the numbers of constraints in \eqref{KS RO max discretized} and \eqref{KS RO min discretized} growing in the data size $n$. One significance of the result is that, despite this growth, as long as one generates the supports of $X$ from a distribution with a heavier tail than the true distribution, and with a size $m$ of order higher than $n$, the confidence guarantee is approximately retained. A key element in explaining this behavior lies in the low complexity of the function class $I(h(\cdot)\leq y)$ (parametrized by $y$) appearing in the constraints and interplayed with the likelihood ratio $dP_X^0/dQ$, as measured by the bracketing number. This number captures the richness of the involved function class with the counts of neighborhoods, each formed by an upper and a lower bounding function that is known as a bracket, in covering the whole class (see the discussion in Appendix \ref{sec:complexity}). A slowly growing (e.g., polynomial in our case) bracketing number turns out to allow the statistic on the output performance measure to be approximated uniformly well with a discretized input distribution, by invoking the empirical process theory for so-called $U$-statistics (\cite{arcones1993limit}). On the other hand, using other moment functions (implied by other test statistics) may not preserve this behavior. This connection to function complexity, which informs the usefulness of sampling-based procedures when integrating with output data, is the first of such kind in the model calibration literature as far as we know.


We have focused on a continuous output distribution in Theorem \ref{main guarantee}. The assumption $\mathbb P(\text{for any\ }P_X\in\hat{\mathcal P}_X, \text{supp}(\gamma(P_X))\cap\{y_j\}_{j=1,\ldots,n}\neq\emptyset)=0$ is a technical condition that ensures the distribution of $h(\mathbf X)$ under $P_X\in\hat{\mathcal P}_X$ does not have overlapping support points as $y_j$'s, which allows us to reduce the KS-implied constraint into the $n$ interval constraints depicted in the theorem. This assumption holds in almost every discrete-event performance measure provided that the considered $P_X$ and $P_Y$ are continuous. On the other hand, if $P_Y$ is discrete, then the theorem holds with the first constraints in \eqref{KS RO max discretized} and \eqref{KS RO min discretized} replaced by \eqref{discrete moment} (with $q$ suitably calibrated as discussed there), without needing the assumption $\mathbb P(\text{for any\ }P_X\in\hat{\mathcal P}_X, \text{supp}(\gamma(P_X))\cap\{y_j\}_{j=1,\ldots,n}\neq\emptyset)=0$.
 
We mention that \cite{ghosh2015robust} provides a similar guarantee for robust optimization problems designed for quantifying input uncertainty. In particular, their analysis allows to give confidence bounds on output performance measures. However, they do not consider the asymptotic confidence guarantee in relation to the data size and the randomized support size. As a consequence, they do not need considering the complexity of the constraints. Moreover, since they handle input uncertainty, the uncertainty sets are more elementary, in contrast to ours which serve as a tool to invert the input-output relation. 

We note that, like Proposition \ref{basic guarantee}, all the results in this section can be similarly extended to a simultaneous guarantee when solving $L$ optimization problems, where each problem has a different objective function $\psi_l(P_X)$. For instance, under the same assumptions as Theorem \ref{main guarantee} with $L$ different objectives in \eqref{KS RO max discretized} and \eqref{KS RO min discretized}, and using the same generated set of support points across all optimization problems, we would obtain that
$$\liminf_{n\to\infty,m/n\to\infty}\mathbb P\left(\hat{\underline Z}_l+O_p\left(\frac{1}{\sqrt m}\right)\leq\psi_l(P_X^0)\leq\hat{\overline Z}_l+O_p\left(\frac{1}{\sqrt m}\right),l=1,\ldots,L\right)\geq1-\alpha$$
where $\hat{\underline Z}_l,\hat{\overline Z}_l$ are the minimum and maximum values of the discretized optimization with objective $\psi_l(P_X)$, and each $O_p(1/\sqrt m)$ is the error term corresponding to each optimization program.

Lastly, we point out that all the results in Sections \ref{sec:formulation} and \ref{sec:KS} hold when we consider $h(\mathbf X,\mathbf W)$ and $g(\mathbf X,\mathbf W)$, where $\mathbf W$ consist of other input variate sequences independent from $\mathbf X$ with known probability distributions. This is as long as we treat all the expectations $E_{P_X}[\cdot]$ as taken jointly under both the product measure of $P_X$ and the known distribution of $\mathbf W$. We provide further remarks in the appendix.


\section{Optimization Procedure}\label{sec:procedure}

This section presents our optimization strategy for (locally) solving \eqref{KS RO max discretized} and \eqref{KS RO min discretized}. Without loss of generality, we only focus on the minimization problem \eqref{KS RO min discretized} since maximization can be converted to minimization by negating the objective. Section \ref{sec:transformation} first discusses the transformation of the stochastic constrained program into a sequence of programs with deterministic convex constraints, using the quadratic penalty method in nonlinear programming. Section \ref{sec:MDSA} then investigates how this transformation can be utilized effectively in a fully iterative stochastic algorithm using MDSA. Section \ref{convergence guarantee} provides a convergence theorem. In the appendix, we also provide an alternate approach that has a similar convergence guarantee but differs in the implementation details.



\subsection{A Stochastic Quadratic Penalty Method}\label{sec:transformation}
When restricted to distributions with support points $\{z_i\}_{i=1,\ldots,m}$, the candidate input distribution $P_X$ can be identified by an $m$-dimensional vector $\mathbf p=(p_1,\ldots,p_m)$ on the probability simplex $\mathcal P:=\{\mathbf p:\sum_{i=1}^mp_i=1,p_i\geq 0\text{ for each }i\}$, where the subscript $X$ is suppressed with no ambiguity. By the vector $\mathbf p$, we mean the distribution that assigns probability $p_i$ to the point $z_i$. The optimization program \eqref{KS RO min discretized} can thus be rewritten as
\begin{equation}
\begin{array}{ll}
\min&\psi(\mathbf p)\\
\text{subject to}&\hat F_Y(y_j+)-\frac{q_{1-\alpha}}{\sqrt n}\leq E_{\mathbf p}[I(h(\mathbf X)\leq y_j)]\leq\hat F_Y(y_j-)+\frac{q_{1-\alpha}}{\sqrt n},j=1,\ldots,n\\
&\mathbf p\in\mathcal P.
\end{array}\label{KS RO min discretized p}
\end{equation}
Note that the constraints in \eqref{KS RO min discretized p} are in general non-convex because the i.i.d.~input sequence means that the expectation $E_{\mathbf p}[I(h(\mathbf X)\leq y_j)]$ is a high-dimensional polynomial in $\mathbf p$. Moreover, this polynomial can involve a huge number of terms and hence its evaluation requires simulation approximation. As far as we know, the literature on dealing with stochastic non-convex constraints is very limited. To overcome this difficulty, we first introduce the quadratic penalty method (\cite{bertsekas1999nonlinear}) to transform program \eqref{KS RO min discretized p} into a sequence of penalized programs with deterministic convex constraints
\begin{equation}
\begin{array}{ll}
\min&\lambda\psi(\mathbf p)+\sum_{j=1}^n (E_{\mathbf p}[I(h(\mathbf X)\leq y_j)]-s_j)^2\\
\text{subject to}&\hat F_Y(y_j+)-\frac{q_{1-\alpha}}{\sqrt n}\leq s_j\leq\hat F_Y(y_j-)+\frac{q_{1-\alpha}}{\sqrt n},j=1,\ldots,n\\
&\mathbf p\in\mathcal P
\end{array}\label{KS RO min discretized penalty}
\end{equation}
where $\mathbf s=(s_1,\ldots,s_n)$ are slack variables and $\lambda>0$ is an inverse measure of the cost/penalty of infeasibility. A related scheme is also used by \cite{wang2008stochastic} in the context of nonconvex stochastic objectives (with deterministic constraints). As $\lambda\to 0$, there is an increasing cost of violating the stochastic constraints, therefore the optimal solution of \eqref{KS RO min discretized penalty} converges to that of \eqref{KS RO min discretized p}, as stated in the following proposition.
\begin{proposition}\label{prop:quadratic penalty}
Suppose \eqref{KS RO min discretized p} has at least one feasible solution. Let $(\mathbf p^*(\lambda),\mathbf s^*(\lambda))$ be an optimal solution of \eqref{KS RO min discretized penalty} indexed at $\lambda$. As $\lambda$ decreases to $0$, every limit point of the sequence $\{\mathbf p^*(\lambda)\}$ is an optimal solution of \eqref{KS RO min discretized p}.
\end{proposition}

As suggested in the proof of Proposition \ref{prop:quadratic penalty}, a mathematically equivalent reformulation of \eqref{KS RO min discretized penalty} with the slack variables optimized is
\begin{equation}
\begin{array}{ll}
\min&\lambda\psi(\mathbf p)+\sum_{j=1}^n (E_{\mathbf p}[I(h(\mathbf X)\leq y_j)]-\Pi_j(E_{\mathbf p}[I(h(\mathbf X)\leq y_j)]))^2\\
\text{subject to}&\mathbf p\in\mathcal P
\end{array}\label{KS RO min discretized penalty2}
\end{equation}
where each $\Pi_j$ is the projection onto the interval $[F_Y(y_j+)-\frac{q_{1-\alpha}}{\sqrt n},F_Y(y_j-)+\frac{q_{1-\alpha}}{\sqrt n}]$ defined as
\begin{equation}\label{pie_j}
\Pi_j(x)=
\begin{cases}
\hat F_Y(y_j+)-\frac{q_{1-\alpha}}{\sqrt n}&\text{if }x<\hat F_Y(y_j+)-\frac{q_{1-\alpha}}{\sqrt n}\\
\hat F_Y(y_j-)+\frac{q_{1-\alpha}}{\sqrt n}&\text{if }x>\hat F_Y(y_j-)+\frac{q_{1-\alpha}}{\sqrt n}\\
x&\text{otherwise}.
\end{cases}
\end{equation}

\subsection{Constrained Stochastic Approximation}\label{sec:MDSA}
Although the formulations \eqref{KS RO min discretized penalty}, \eqref{KS RO min discretized penalty2} are still non-convex, their constraints are convex and deterministic, which can be handled more easily using SA than in the original formulation \eqref{KS RO min discretized p}. This section investigates the design and analysis of an MDSA algorithm for finding local optima of \eqref{KS RO min discretized p} by solving \eqref{KS RO min discretized penalty} with decreasing values of $\lambda$. The appendix would illustrate another algorithm that uses formulation \eqref{KS RO min discretized penalty2} instead of \eqref{KS RO min discretized penalty}.


To describe the algorithm, MD finds the next iterate via optimizing the objective function linearized at the current iterate, together with a penalty on the distance of movement of the iterate. When the objective function is only accessible via simulation, the linearized objective function, or the gradient, at each iteration can only be estimated with noise, in which case the procedure becomes MDSA (\cite{nemirovski2009robust}). More precisely, when applied to the formulation \eqref{KS RO min discretized penalty} with slack variables, MDSA solves the following optimization given a current iterate $(\mathbf p^k,\mathbf s^k)$
\begin{equation}
\begin{array}{ll}
\min&\gamma^k(\lambda\hat{\bm\Psi}^k+\hat{\bm\phi}_{\mathbf p}^k)'(\mathbf p-\mathbf p^k)+\beta^k\hat{\bm\phi}_{\mathbf s}^{k\prime}(\mathbf s-\mathbf s^k)+V(\mathbf p^k,\mathbf p)+\frac{1}{2}\Vert\mathbf s-\mathbf s^k\Vert_2^2\\
\text{subject to}&\hat F_Y(y_j+)-\frac{q_{1-\alpha}}{\sqrt n}\leq s_j\leq\hat F_Y(y_j-)+\frac{q_{1-\alpha}}{\sqrt n},j=1,\ldots,n\\
&\mathbf p\in\mathcal P
\end{array}\label{step optimization2}
\end{equation}
where $\hat{\bm\Psi}^k$ carries the gradient information of the target performance measure $\psi$ at $\mathbf p^k$, while $\hat{\bm\phi}_{\mathbf p}^k$ and $\hat{\bm\phi}_{\mathbf s}^k$ contain the gradient information of the penalty function in \eqref{KS RO min discretized penalty} with respect to $\mathbf p,\mathbf s$ respectively. The sum $V(\mathbf p^k,\mathbf p)+\frac{1}{2}\Vert\mathbf s-\mathbf s^k\Vert_2^2$ serves as the penalty on the movement of the iterate, where $\Vert\cdot\Vert_2$ denotes the standard Euclidean distance, and $V(\cdot,\cdot)$ defined as
\begin{equation}
V(\mathbf p,\mathbf q)=\sum_{i=1}^nq_i\log\frac{q_i}{p_i}\label{KL}
\end{equation}
is the KL divergence between two probability measures. This particular choice of $V$ has been shown (\cite{nemirovski2009robust}) to have superior performance to other choices like the Euclidean distance when the decision space is the probability simplex. Different from traditional SA, the step sizes $\gamma^k$ and $\beta^k$, used for updating $\mathbf p$ and $\mathbf s$ in \eqref{step optimization2}, are different, the rationale for which shall be discussed in Section \ref{convergence guarantee}.

However, iterations in the form of \eqref{step optimization2} can only find optima of \eqref{KS RO min discretized penalty} for a particular penalty coefficient $\lambda$ while retrieving the optimal solution of the original problem \eqref{KS RO min discretized p} through \eqref{KS RO min discretized penalty} hinges on sending $\lambda$ to $0$. Literature on deterministic optimization suggests solving the penalized optimization repeatedly for a set of decreasing values of $\lambda$, but it could be difficult to tell when to stop decreasing the $\lambda$ in our stochastic case. In order to output the optimal solution in one single run, we decrease $\lambda$ together with the step size from one iteration to the next, hence arrive at the following sequential joint solution-and-penalty-updating routine
\begin{equation}
\begin{array}{ll}
\min&\gamma^k(\lambda^k\hat{\bm\Psi}^k+\hat{\bm\phi}_{\mathbf p}^k)'(\mathbf p-\mathbf p^k)+\beta^k\hat{\bm\phi}_{\mathbf s}^{k\prime}(\mathbf s-\mathbf s^k)+V(\mathbf p^k,\mathbf p)+\frac{1}{2}\Vert\mathbf s-\mathbf s^k\Vert_2^2\\
\text{subject to}&\hat F_Y(y_j+)-\frac{q_{1-\alpha}}{\sqrt n}\leq s_j\leq\hat F_Y(y_j-)+\frac{q_{1-\alpha}}{\sqrt n},j=1,\ldots,n\\
&\mathbf p\in\mathcal P
\end{array}\label{step optimization2 lambda}
\end{equation}
where $\lambda^k$ is appropriately chosen in conjunction with $\gamma^k,\beta^k$, and decreases to $0$. To implement the fully sequential scheme, we need to investigate: 1) how to obtain $\hat{\bm\Psi}^k,\hat{\bm\phi}_{\mathbf p}^{k}$ and $\hat{\bm\phi}_{\mathbf s}^{k}$, 2) efficient solution method for program \eqref{step optimization2 lambda}, and 3) how to select the parameters $\gamma^k,\beta^k$ and $\lambda^k$. The next two subsections present the first two investigations respectively, while Section \ref{convergence guarantee} will analyze the convergence of the algorithm in relation to the parameter choices.


%

\subsubsection{Gradient Estimation and Restricted Programs. }\label{sec:gradient}
Denote by $W(\mathbf p)$ the penalty function in \eqref{KS RO min discretized penalty2}, and by $W_s(\mathbf p,\mathbf s)$ the quadratic penalty in \eqref{KS RO min discretized penalty} where the subscript $s$ refers to ``slack variable''. These are functions of variables on the probability simplex, for which naive differentiation may not lead to simulable object since an arbitrary perturbation may shoot out of the simplex. \cite{gl15_1} and \cite{ghosh2015mirror} have used the idea of Gateaux derivative (in the sense described in Chapter 6 of \cite{serfling2009approximation}) to obtain simulable representations of gradients of expectation-type performance measures. We generalize their result to sums of functions of expectations:


\begin{proposition}
We have:
\begin{enumerate}
\item Suppose $\psi,W,W_s(\cdot,\mathbf s)$ are differentiable in the probability simplex $\mathcal P$, then
\begin{align}
\nabla \psi(\mathbf p)'(\mathbf q-\mathbf p)&=\bm\Psi(\mathbf p)'(\mathbf q-\mathbf p)\label{derivative1}\\
\nabla W(\mathbf p)'(\mathbf q-\mathbf p)&=\bm\phi(\mathbf p)'(\mathbf q-\mathbf p)\label{derivative2}\\
\nabla_{\mathbf p} W_s(\mathbf p,\mathbf s)'(\mathbf q-\mathbf p)&=\bm\phi_{\mathbf p}(\mathbf p,\mathbf s)'(\mathbf q-\mathbf p)\label{derivative3}
\end{align}
for any $\mathbf p,\mathbf q\in\mathcal P$, where the Gateaux derivatives $\bm\Psi(\mathbf p)=(\Psi_1(\mathbf p),\ldots,\Psi_m(\mathbf p))'$, $\bm\phi(\mathbf p)=(\phi_1(\mathbf p),\ldots,\phi_m(\mathbf p))'$, $\bm\phi_{\mathbf p}(\mathbf p,\mathbf s)=(\phi_{\mathbf p,1}(\mathbf p,\mathbf s),\ldots,\phi_{\mathbf p,m}(\mathbf p,\mathbf s))'$, and
\begin{align}
\Psi_i(\mathbf p)&=\frac{d}{d\epsilon}\psi((1-\epsilon)\mathbf p+\epsilon\mathbf 1_i)\Big|_{\epsilon=0^+}\label{Gateaux preservation1}\\
\phi_i(\mathbf p)&=\frac{d}{d\epsilon}W((1-\epsilon)\mathbf p+\epsilon\mathbf 1_i)\Big|_{\epsilon=0^+}\label{Gateaux preservation2}\\
\phi_{\mathbf p,i}(\mathbf p,\mathbf s)&=\frac{d}{d\epsilon}W_s((1-\epsilon)\mathbf p+\epsilon\mathbf 1_i,\mathbf s)\Big|_{\epsilon=0^+}\label{Gateaux preservation3}
\end{align}
\item Assume $\mathbf p=(p_1,\ldots,p_m)$ where each $p_i>0$. Then the Gateaux derivatives \eqref{Gateaux preservation1}\eqref{Gateaux preservation2}\eqref{Gateaux preservation3} are finite and can be expressed as
\begin{align}
\Psi_i(\mathbf p)&=E_{\mathbf p}[g(\mathbf X)S_i(\mathbf X;\mathbf p)]\label{score function1}\\
\phi_i(\mathbf p)&=2\sum_{j=1}^n(E_{\mathbf p}[I(h(\mathbf X)\leq y_j)]-\Pi_j(E_{\mathbf p}[I(h(\mathbf X)\leq y_j)]))E_{\mathbf p}[I(h(\mathbf X)\leq y_j)S_i(\mathbf X;\mathbf p)]\label{score function2}\\
\phi_{\mathbf p,i}(\mathbf p,\mathbf s)&=2\sum_{j=1}^n(E_{\mathbf p}[I(h(\mathbf X)\leq y_j)]-s_j)E_{\mathbf p}[I(h(\mathbf X)\leq y_j)S_i(\mathbf X;\mathbf p)]\label{score function3}
\end{align}
where
$$S_i(\mathbf x;\mathbf p)=\sum_{t=1}^{S}\frac{I_i(x_t)}{p_i}-S\text{ for \eqref{score function1}},\ \text{and }\sum_{t=1}^{T}\frac{I_i(x_t)}{p_i}-T\text{ for \eqref{score function2}\eqref{score function3}}.$$
Here $I_i(x)=1$ if $x=z_i$ and 0 otherwise, and $\mathbf X$ is the i.i.d. input process generated under $\mathbf p$.
\label{derivative form}
\end{enumerate}
\label{lemma:derivative1}
\end{proposition}
The representations \eqref{score function1} and \eqref{score function3} suggest the following unbiased estimators for the gradient of $\psi$,  $\bm\Psi(\mathbf p)=(\Psi_i(\mathbf p))_{i=1}^m$, and the gradient of the penalty function, $\bm\phi_{\mathbf p}(\mathbf p,\mathbf s)=(\phi_{\mathbf p,i}(\mathbf p,\mathbf s))_{i=1}^m$
\begin{align}
&\hat{\Psi}_i(\mathbf p)=\frac{1}{M_3}\sum_{r=1}^{M_3}g(\mathbf X^{(r)})S_i(\mathbf X^{(r)};\mathbf p)\label{gradient estimator 1}\\
&\hat{\phi}_{\mathbf p,i}(\mathbf p,\mathbf s)=2\sum_{j=1}^n\frac{1}{M_1}\sum_{r=1}^{M_1}(I(h(\mathbf X^{(r)})\leq y_j)-s_j)\frac{1}{M_2}\sum_{r=1}^{M_2}I(h(\tilde{\mathbf X}^{(r)})\leq y_j)S_i(\tilde{\mathbf X}^{(r)};\mathbf p)\label{gradient estimator 3}
\end{align}
where $\mathbf X^{(r)}$ and $\tilde{\mathbf X}^{(r)}$ are independent copies of the i.i.d. input process generated under $\mathbf p$ and are used simultaneously for all $i,j$. By direct differentiation, a straightforward unbiased estimator for $\bm\phi_{\mathbf s}(\mathbf p,\mathbf s)=(\phi_{\mathbf s,j}(\mathbf p,\mathbf s))_{j=1}^n$, the gradient of the penalty function with respect to $\mathbf s$, is
\begin{equation}
\hat{\phi}_{\mathbf s, j}(\mathbf p,\mathbf s)=\frac{-2}{M_1}\sum_{r=1}^{M_1}(I(h(\mathbf X^{(r)})\leq y_j)-s_j).\label{gradient estimator 4}
\end{equation}
Our main procedure (shown in Algorithm \ref{algo2} momentarily) uses the above gradient estimators, while an alternate MDSA depicted in Algorithm \ref{algo1} in the appendix solves \eqref{KS RO min discretized penalty2} using a biased estimator of $\bm\phi(\mathbf p)$ conferred by \eqref{score function2}.

Note that the above gradient estimators are available thanks to the KS-implied constraints we introduced. By the reformulation in Theorem \ref{KS guarantee}, the constraints in \eqref{KS RO max} and \eqref{KS RO min} become ($T$-fold) expectation-type constraints. Thus, when differentiating the squared expectation in the quadratic penalty, the gradient becomes the product of two $T$-fold expectations, one with the extra factor $S_i(\cdot;\mathbf p)$ which can be interpreted as a score function. This then allows unbiased estimation of the gradient by generating two independent batches of simulation runs each for one of the expectations. Using other statistics to induce the constraints may not lead to such a convenient form.

Note that the $S_i(\cdot;\mathbf p)$ in the gradient estimators \eqref{gradient estimator 1} and \eqref{gradient estimator 3} contains $p_i$ at the denominator, so a small $p_i$ can blow up the variances of the estimators and in turn adversely affect the convergence of MDSA. To ensure convergence, we make an adjustment to our procedure and solve the following restricted version of \eqref{KS RO min discretized p}
\begin{equation}
\begin{array}{ll}
\min&\psi(\mathbf p)\\
\text{subject to}&\hat F_Y(y_j+)-\frac{q_{1-\alpha}}{\sqrt n}\leq E_{\mathbf p}[I(h(\mathbf X)\leq y_j)]\leq\hat F_Y(y_j-)+\frac{q_{1-\alpha}}{\sqrt n},j=1,\ldots,n\\
&\mathbf p\in\mathcal P(\epsilon)
\end{array}\label{KS RO min discretized p epsilon}
\end{equation}
where the restricted probability simplex $\mathcal P(\epsilon):=\{\mathbf p\in \mathcal P:p_i\geq \epsilon\text{ for each }i\}$. Accordingly, the full simplex $\mathcal P$ in the penalized program \eqref{KS RO min discretized penalty} and stepwise subproblem \eqref{step optimization2 lambda} has to be replaced by $\mathcal P(\epsilon)$.

To maintain the statistical guarantee provided by Theorem \ref{main guarantee} when solving the restricted programs, the shrinking size $\epsilon$ has to be appropriately chosen. Theorem \ref{main guarantee epsilon} below indicates that it suffices to choose $\epsilon$ smaller than $1/(m\sqrt n)$ in case of bounded $g(\mathbf X)$.
\begin{theorem}\label{main guarantee epsilon}
Denote by $\hat{\overline Z}_{\epsilon}$ and $\hat{\underline Z}_{\epsilon}$ the maximum and minimum of $\psi(\mathbf p)$ in the feasible set of \eqref{KS RO min discretized p epsilon}. In addition to the conditions of Theorem \ref{main guarantee}, further assume that $g(\mathbf X)$ is bounded. If $\epsilon$ is chosen such that $\epsilon=o\big(\frac{1}{m\sqrt n}\big)$ then we have
$$\liminf_{n\to\infty,m/n\to\infty}\mathbb P\left(\hat{\underline Z}_{\epsilon}+O_p\left(m\epsilon+\frac{1}{\sqrt m}\right)\leq\psi(P_X^0)\leq\hat{\overline Z}_{\epsilon}+O_p\left(m\epsilon+\frac{1}{\sqrt m}\right)\right)\geq 1-\alpha.$$
\end{theorem}
In particular, the first type of target quantities we consider has a bounded $g(\mathbf X)$. Note that the original optimization itself already poses an error of size $O_p(1/\sqrt m)$ in the confidence bounds (Theorem \ref{main guarantee}), so to keep the error at the same level one can use an $\epsilon=O(1/m^{\frac{3}{2}})$ (recall that $m/n\to \infty$). Since the variances of our gradient estimators \eqref{gradient estimator 1}\eqref{gradient estimator 3} can be shown inversely proportional to the components $p_i$ (\cite{ghosh2015robust}), such an $\epsilon$ gives rise to variances of order $O(m^{\frac{3}{2}})$. We point out that this is only slightly worse than the best attainable order $O(m)$, which results from the fact that the average size of $\mathbf p$ in the $m$-dimensional probability simplex is $1/m$.

\subsubsection{Solving Stepwise Subproblem in MDSA. }\label{sec:stepwise}
Since we are now solving the restricted version of subproblem \eqref{step optimization2 lambda}, consider the following generic form
\begin{equation}\label{stepwise slack epsilon}
\begin{array}{ll}
\min&\bm\xi'(\mathbf q-\mathbf p)+\bm\eta'(\mathbf t-\mathbf s)+V(\mathbf p,\mathbf q)+\frac{1}{2}\Vert\mathbf t-\mathbf s\Vert_2^2\\
\text{subject to}&\hat F_Y(y_j+)-\frac{q_{1-\alpha}}{\sqrt n}\leq t_j\leq\hat F_Y(y_j-)+\frac{q_{1-\alpha}}{\sqrt n},j=1,\ldots,n\\
&\mathbf q\in\mathcal P(\epsilon).
\end{array}
\end{equation}
Because the objective and the feasible set are both separable in $\mathbf q$ and $\mathbf t$, the above program can be decomposed into two independent programs. One is
\begin{equation}
\begin{array}{ll}
\min&\bm\xi'(\mathbf q-\mathbf p)+V(\mathbf p,\mathbf q)\\
\text{subject to}&\mathbf q\in\mathcal P(\epsilon)
\end{array}\label{generic_epsilon}
\end{equation}
and the other is
\begin{equation}\label{generic_slack}
\begin{array}{ll}
\min&\bm\eta'(\mathbf t-\mathbf s)+\frac{1}{2}\Vert\mathbf t-\mathbf s\Vert_2^2\\
\text{subject to}&\hat F_Y(y_j+)-\frac{q_{1-\alpha}}{\sqrt n}\leq t_j\leq\hat F_Y(y_j-)+\frac{q_{1-\alpha}}{\sqrt n},j=1,\ldots,n.
\end{array}
\end{equation}

Program \eqref{generic_slack} is exactly the step-wise routine that appears in the standard gradient descent whose solution takes the form
\begin{equation*}
t^*_j=\Pi_j(s_j-\eta_j)
\end{equation*}
where $\Pi_j$ is the projection defined in \eqref{pie_j}.

The solution of program \eqref{generic_epsilon} has a semi-explicit expression as shown in the following proposition.
\begin{proposition}\label{sol:constrained stepwise}
The optimal solution of the stepwise subproblem \eqref{generic_epsilon} with $0\leq \epsilon<1/m$ is
\begin{equation}\label{sol2}
q^*_i=\frac{\max\{\eta^*,p_ie^{-\xi_i}\}}{\sum_{i=1}^m\max\{\eta^*,p_ie^{-\xi_i}\}}
\end{equation}
where $\eta^*\in [0,\max_ip_ie^{-\xi_i})$ solves the equation
\begin{equation}\label{threshold:eta}
\epsilon=\mu(\eta^*):=\frac{\eta^*}{\sum_{i=1}^m\max\{\eta^*,p_ie^{-\xi_i}\}}.
\end{equation}
\end{proposition}



Proposition \ref{sol:constrained stepwise} suggests a procedure for solving \eqref{generic_epsilon} that involves a root-finding problem \eqref{threshold:eta}. To design an efficient root-finding routine, note that the function $\mu(\eta)$ is strictly increasing in $\eta$. More importantly, it consists of at most $m$ smooth pieces, and on the $i$-th piece it takes the form
\begin{equation*}
\mu(\eta)=\frac{\eta}{i\eta+\sum_{i'=i+1}^mp_{(i')}e^{-\xi_{(i')}}},\text{ if }p_{(i)}e^{-\xi_{(i)}}\leq\eta\leq p_{(i+1)}e^{-\xi_{(i+1)}}
\end{equation*}
where $(p_{(1)}e^{-\xi_{(1)}},\ldots,p_{(m)}e^{-\xi_{(m)}})$ is obtained by sorting $(p_{1}e^{-\xi_{1}},\ldots,p_{m}e^{-\xi_{m}})$ in ascending order. Thus one can first locate which piece the root $\eta^*$ lies on by comparing the values of $\mu$ with $\epsilon$ at the points $p_{(i)}e^{-\xi_{(i)}}$ and then compute $\eta^*$ in closed form from the above expression on that piece. This efficient sort-and-search procedure is described in Algorithm \ref{sort_search} whose proof follows from straightforward algebraic verification and hence is omitted.
\begin{algorithm}[h]
\caption{Sort-and-search for solving \eqref{generic_epsilon} with $0\leq \epsilon<1/m$}
\begin{algorithmic}

\State \textbf{1.} Sort $(p_1e^{-\xi_1},\ldots,p_me^{-\xi_m})$ into ascending order $(p_{(1)}e^{-\xi_{(1)}},\ldots,p_{(m)}e^{-\xi_{(m)}})$, and let $p_{(0)}e^{-\xi_{(0)}}=0$

\State \textbf{2.} Search for the $i^*$ from $0$ to $m-1$ such that
\begin{equation*}
\frac{p_{(i^*)}e^{-\xi_{(i^*)}}}{i^*p_{(i^*)}e^{-\xi_{(i^*)}}+\sum_{i=i^*+1}^mp_{(i)}e^{-\xi_{(i)}}}\leq \epsilon< \frac{p_{(i^*+1)}e^{-\xi_{(i^*+1)}}}{(i^*+1)p_{(i^*+1)}e^{-\xi_{(i^*+1)}}+\sum_{i=i^*+2}^mp_{(i)}e^{-\xi_{(i)}}}
\end{equation*}
\State \textbf{3.} Output $q^*_i$ according to \eqref{sol2} with
\begin{equation*}
\eta^*=\frac{\epsilon\sum_{i=i^*+1}^mp_{(i)}e^{-\xi_{(i)}}}{1-\epsilon i^*}
\end{equation*}
\end{algorithmic}\label{sort_search}
\end{algorithm}

\subsection{Convergence Analysis}\label{convergence guarantee}
We depict our MDSA procedure in Algorithm \ref{algo2}. Steps 1, 2 and 3 of the procedure estimate the gradients using the estimators proposed in Section \ref{sec:gradient}, and Step 4 updates the decision variable with step size $\gamma^k$ and the slack variables with step size $\beta^k$. Steps 1-4 combined are in effect solving the stepwise subproblem \eqref{step optimization2 lambda} with $\mathcal P$ replaced by $\mathcal P(\epsilon)$. Therefore by iterating with decreasing penalty coefficient $\lambda^k$, Algorithm \ref{algo2} searches for the optimum of the restricted formulation \eqref{KS RO min discretized p epsilon}.




\begin{algorithm}[h]
  \caption{MDSA for solving \eqref{KS RO min discretized penalty}}
  \textbf{Input: }A small parameter $\epsilon>0$, initial solution $\mathbf p^1\in\mathcal P(\epsilon)=\{\mathbf p:\sum_{i=1}^mp_i=1,p_i\geq\epsilon\text{\ for\ }i=1,\ldots,m\}$ and $\mathbf s^1\in [\hat F_Y(y_1+)-\frac{q_{1-\alpha}}{\sqrt n},\hat F_Y(y_1-)+\frac{q_{1-\alpha}}{\sqrt n}]\times \cdots\times [\hat F_Y(y_n+)-\frac{q_{1-\alpha}}{\sqrt n},\hat F_Y(y_n+)-\frac{q_{1-\alpha}}{\sqrt n}]$, a step size sequence $\gamma^k$ for $\mathbf p$, a penalty sequence $\lambda^k$, a step size sequence $\beta^k$ for $\mathbf s$, and sample sizes $M_1,M_2,M_3$.

  \textbf{Iteration: }For $k=1,2,\ldots$, do the following: Given $\mathbf p^k,\mathbf s^k$,
\begin{algorithmic}


\State \textbf{1.} Estimate $\hat{\bm\phi}_{\mathbf p}^k=(\hat\phi_{\mathbf p,1}^k,\ldots,\hat\phi_{\mathbf p,m}^k)$, the gradient of the penalty term with respect to $\mathbf p$, with
$$\hat\phi_{\mathbf p,i}^k=2\sum_{j=1}^n\frac{1}{M_1}\sum_{r=1}^{M_1}(I(h(\mathbf X^{(r)})\leq y_j)-s^k_j)\frac{1}{M_2}\sum_{r=1}^{M_2}I(h(\tilde{\mathbf X}^{(r)})\leq y_j)S_i(\tilde{\mathbf X}^{(r)};\mathbf p^k)$$
where $\mathbf X^{(r)},\tilde{\mathbf X}^{(r)}$ are $M_1$ and $M_2$ independent copies of the input process generated under $\mathbf p^k$.

\State \textbf{2.} Estimate $\hat{\bm\Psi}^k=(\hat\Psi_1^k,\ldots,\hat\Psi_m^k)$, the gradient of $E_{\mathbf p}[g(\mathbf X)]$, with
$$\hat\Psi_i^k=\frac{1}{M_3}\sum_{r=1}^{M_3}g(\tilde{\tilde{\mathbf X}}^{(r)})S_i(\tilde{\tilde{\mathbf X}}^{(r)};\mathbf p^k)$$
where $\tilde{\tilde{\mathbf X}}^{(r)}$ are another $M_3$ independent copies of the input process generated under $\mathbf p^k$.

\State \textbf{3.} Estimate $\hat{\bm\phi}_{\mathbf s}^k=(\hat\phi_{\mathbf s,1}^k,\ldots,\hat\phi_{\mathbf s,n}^k)$, the gradient of the penalty term with respect to $\mathbf s$, with
$$\hat\phi_{\mathbf s,j}^k=-\frac{2}{M_1+M_2}\big(\sum_{r=1}^{M_1}(I(h(\mathbf X^{(r)})\leq y_j)-s^k_j)+\sum_{r=1}^{M_2}(I(h(\tilde{\mathbf X}^{(r)})\leq y_j)-s^k_j)\big)$$
where ${\mathbf X}^{(r)},\tilde{\mathbf X}^{(r)}$ are the same replications used in Step 1.

\State \textbf{4.} Compute
$\mathbf p^{k+1}=(p_1^{k+1},\ldots,p_m^{k+1})$ by running Algorithm \ref{sort_search} with $p_i=p_i^{k}$ and $\xi_i=\gamma^k(\lambda^k\hat\Psi_i^k+\hat\phi_{\mathbf p,i}^k)$,
and compute $\mathbf s^{k+1}=(s_1^{k+1},\ldots,s_n^{k+1})$ by
\begin{equation*}
s_j^{k+1}=\Pi_j(s_j^k-\beta^k\hat\phi_{\mathbf s,j}^k).
\end{equation*}



\end{algorithmic}\label{algo2}
\end{algorithm}

To provide convergence guarantee for Algorithm \ref{algo2}, we assume the following:
\begin{assumption}\label{cond:constrained}
The restricted program \eqref{KS RO min discretized p epsilon} has a unique optimal solution $\mathbf p^*_{\epsilon}\in \mathcal P(\epsilon)$ such that for any feasible $\mathbf p\in \mathcal P(\epsilon)$ and $\mathbf p\neq \mathbf p^*_{\epsilon}$ it holds $\bm\Psi(\mathbf p)'(\mathbf p-\mathbf p^*_{\epsilon})>0$, and for any infeasible $\mathbf p\in \mathcal P(\epsilon)$ it holds $\bm\phi(\mathbf p)'(\mathbf p-\mathbf p^*_{\epsilon})>0$, where $\bm\Psi,\bm\phi$ are respectively the Gateaux derivatives of the target quantity $\psi$ and the quadratic penalty function $\sum_{j=1}^n (E_{\mathbf p}[I(h(\mathbf X)\leq y_j)]-\Pi_j(E_{\mathbf p}[I(h(\mathbf X)\leq y_j)]))^2$ in \eqref{KS RO min discretized penalty2}.
\end{assumption}
\begin{assumption}\label{cond:penalized}
There is some threshold $\lambda_{\epsilon}>0$ such that
\begin{description}
\item[1.]for any $\lambda\in (0,\lambda_{\epsilon}]$ the optimization problem \eqref{KS RO min discretized penalty2} with $\mathcal P$ replaced by $\mathcal P(\epsilon)$ has a unique optimal solution $\mathbf p^*_{\epsilon}(\lambda)\in \mathcal P(\epsilon)$ such that for any $\mathbf p\in \mathcal P(\epsilon)$ it holds $(\lambda\bm\Psi(\mathbf p)+\bm\phi(\mathbf p))'(\mathbf p-\mathbf p^*_{\epsilon}(\lambda))\geq 0$
\item[2.]$\mathbf p^*_{\epsilon}(\lambda)$ as a function of $\lambda\in (0,\lambda_{\epsilon}]$ has finite total variation, meaning that there exists a constant $M_{\epsilon}>0$ such that $\sum_{i=0}^{K-1}\norm{\mathbf p^*_{\epsilon}(\lambda_i)-\mathbf p^*_{\epsilon}(\lambda_{i+1})}\leq M_{\epsilon}$
for any $0<\lambda_K<\cdots<\lambda_1<\lambda_0\leq  \lambda_{\epsilon}$ and $K$.
\end{description}
\end{assumption}
\begin{assumption}\label{cond:lambda}
$\norm{\mathbf p^*_{\epsilon}(\lambda)-\mathbf p^*_{\epsilon}}=O(\lambda)$ as $\lambda\to 0$.
\end{assumption}
The condition $(\lambda\bm\Psi(\mathbf p)+\bm\phi(\mathbf p))'(\mathbf p-\mathbf p^*_{\epsilon}(\lambda))\geq 0$ in Assumption \ref{cond:penalized} is a weakened version of the general convexity criterion that has appeared in online learning (e.g., \cite{bottou1998online}) and SA (e.g., \cite{benveniste2012adaptive,broadie2011general}) literature. For a minimization problem with objective $f(x)$ and minimizer $x^*$, this criterion refers to the condition that $\nabla f(x)'(x-x^*)>0$ for any $x\neq x^*$. A geometric interpretation of it is that the opposite of the gradient direction always points to the optimum. Part 1 of Assumption \ref{cond:penalized} stipulates that the criterion holds weakly for the penalized program \eqref{KS RO min discretized penalty2} when the penalty coefficient $\lambda$ lies in a small neighborhood of zero. Assumption \ref{cond:constrained} can be viewed as the same criterion for the limit case $\lambda=0$. To explain, at a feasible solution $\mathbf p$ of \eqref{KS RO min discretized p epsilon} the derivative $\bm\phi(\mathbf p)$ vanishes hence the criterion in Assumption \ref{cond:penalized} reduces to $\bm\Psi(\mathbf p)'(\mathbf p-\mathbf p^*_{\epsilon}(\lambda))\geq 0$ when $\lambda>0$, which in the limit $\lambda\to 0$ forces $\bm\Psi(\mathbf p)'(\mathbf p-\mathbf p^*_{\epsilon})\geq 0$ since $\mathbf p^*_{\epsilon}(\lambda)\to \mathbf p^*_{\epsilon}$. Whereas for an infeasible solution $\mathbf p$ the derivative $\bm\phi(\mathbf p)$ is non-zero, thus the criterion becomes $\bm\phi(\mathbf p)'(\mathbf p-\mathbf p^*_{\epsilon})\geq 0$ as $\lambda\to 0$ because $\lambda\bm\Psi(\mathbf p)\to \mathbf 0$ and $\mathbf p^*_{\epsilon}(\lambda)\to \mathbf p^*_{\epsilon}$. Note that Assumption \ref{cond:constrained} further requires the two inequalities to hold strictly.

Part 2 of Assumption \ref{cond:penalized} and Assumption \ref{cond:lambda} impose mild regularity conditions on the solution path of \eqref{KS RO min discretized penalty2} parametrized by $\lambda$. In fact, the solution path is expected to be continuously differentiable in $\lambda$, a stronger property than the assumptions. The reason is that the optimal solution $\mathbf p^*_{\epsilon}(\lambda)$ has to satisfy the set of KKT conditions which is smooth in the decision variable $\mathbf p$ and the penalty coefficient $\lambda$, hence an application of the implicit function theorem reveals the continuous differentiability of $\mathbf p^*_{\epsilon}(\lambda)$ in $\lambda$.



When the target quantity $\psi(\mathbf p)=\mathbf c'\mathbf p$ for some $\mathbf c\in \R^m$, which includes the first type of target quantities we consider in Section \ref{sec:formulation}, the condition $\bm\Psi(\mathbf p)'(\mathbf p-\mathbf p^*_{\epsilon})>0$ in Assumption \ref{cond:constrained} is guaranteed to hold. To explain, note that the feasible set of program \eqref{KS RO min discretized p epsilon} is supported by the hyperplane $\{\mathbf p:\mathbf c'\mathbf p=\mathbf c'\mathbf p^*_{\epsilon}\}$ at the optimum $\mathbf p^*_{\epsilon}$ even if the feasible set is non-convex, and any non-optimal solution $\mathbf p$ will lie in the strict half-space $\{\mathbf p:\mathbf c'\mathbf p> \mathbf c'\mathbf p^*_{\epsilon}\}$ which is exactly the condition in Assumption \ref{cond:constrained}. However, the second condition $\bm\phi(\mathbf p)'(\mathbf p-\mathbf p^*_{\epsilon})>0$ could still be hard to verify because of the nonlinearity of the constraint functions $E_{\mathbf p}[I(h(\mathbf X)\leq y_j)]$. In our numerical experiments, we investigate the use of multi-start and show that our procedure appears to perform well empirically.


Our convergence guarantee of Algorithm \ref{algo2} is stated in Theorem \ref{thm:algo2}, whose proof follows the framework in \cite{blum1954multidimensional} that considers SA on unconstrained problems.

\begin{theorem}
Under Assumptions \ref{cond:constrained}, \ref{cond:penalized} and \ref{cond:lambda}, if the step size sequences $\{\gamma^k\},\{\beta^k\}$ and the penalty sequence $\{\lambda^k\}$ of Algorithm \ref{algo2} are chosen as
\begin{equation}\label{stepsize}
\begin{aligned}
&\gamma^k=\frac{a}{k^{\alpha_1}},\ \frac{3}{4}<\alpha_1\leq 1\\
&\beta^k=\frac{b}{k^{\alpha_2}},\ 2-2\alpha_1<\alpha_2<2\alpha_1-1\\
&\lambda^k=
\begin{cases}
\frac{c}{k^{\alpha_3}},\ 0<\alpha_3\leq 1-\alpha_1&\text{if }\frac{3}{4}<\alpha_1<1\\
\frac{c}{\log k}&\text{if }\alpha_1=1
\end{cases}
\end{aligned}
\end{equation}
then $\mathbf p^k$ generated in Algorithm \ref{algo2} converges to $\mathbf p^*_{\epsilon}$ a.s..\label{thm:algo2}
\end{theorem}
Here $\gamma^k$ and $\beta^k$ are chosen in such a way that the slack variables $\mathbf s^k$ is guaranteed to stay close to the projections $\Pi_j(E_{\mathbf p^k}[I(h(\mathbf X)\leq y_j)])$ and hence the MDSA is effectively solving \eqref{KS RO min discretized penalty2}. Note that the choice of penalty coefficient $\lambda^k$ only depends on the step size $\gamma^k$. The rule of thumb is that $\gamma^k\lambda^k$ should sum up to $\infty$, as indicated by the relation between $\alpha_1$ and $\alpha_3$ in \eqref{stepsize}. This ensures sufficient exploration of the feasible region of \eqref{KS RO min discretized p epsilon}, the rationale of which will be further elaborated in Appendix \ref{sec:proofs algorithm}.

Finally, we mention that in the presence of a collection of auxiliary input sequences $\mathbf W$ with known distribution that is independent of $\mathbf X$, namely that we now have $h(\mathbf X,\mathbf W)$ instead of $h(\mathbf X)$ and $g(\mathbf X,\mathbf W)$ instead of $g(\mathbf X)$, all the results in this section hold by viewing $E_{\mathbf p}[\cdot]$ as taken jointly with respect to the product measure of $\mathbf p$ and the true distribution of $\mathbf W$. In Algorithm \ref{algo2} (and also the other algorithms in the appendix), one only needs to simulate the independent $\mathbf W$ in conjunction with $\mathbf X$ in each replication, e.g., $h(\mathbf X^{(r)},\mathbf W^{(r)})$ instead of $h(\mathbf X^{(r)})$. Appendix \ref{sec:procedure proofs} provides further discussion.

\section{Numerical Results}\label{sec:numerics}
This section provides numerical illustration of our methodology. We focus on a stylized M/G/1 queue, where we assume known i.i.d. unit rate exponential interarrival times. Our goal is to calibrate the unknown i.i.d. service time distribution $P_X$ given the output data. Here, we assume the collection of data for the averaged wait time of the first 20 customers, starting from the empty state. Say these observations are i.i.d. (e.g., among different days or work cycles), denoted $y_1,\ldots,y_n$. The data size $n$ varies from $30$ to $100$ in our experiments.

We consider two target quantities of interest $\psi(P_X)$: 1) the expected averaged queue length seen by the first 20 customers. This performance measure, though related to the waiting time data, is not directly observable and depends on the known service time distribution; 2) the distribution function of the service time. We also consider two different ``true" service time distributions, first one is exponential with rate $1.2$, and second one is a mixture of beta distributions that has a bimodal shape. We set the confidence level to be $95\%$, i.e., $\alpha=5\%$.


Since the input distribution of interest and the output distribution are both continuous, we use optimization programs \eqref{KS RO max discretized} and \eqref{KS RO min discretized} to infer the confidence bounds on $\psi(P_X^0)$. From Theorem \ref{main guarantee}, we first randomly sample $m$ support points from some ``safe" input distribution (i.e., distribution believed to have heavier tail than the truth), where $m$ varies from $100$ to $500$ in our experiments. Then we implement Algorithm \ref{algo2}. In our implementation we choose $M_1 = M_2 = M_3 = 100$ ,  $\gamma^k =  a/k^{0.8}$, $\beta^k = b/k^{0.5}$, $\lambda^k = c/k^{0.2}$, in which the constants $a, b, c$ will be determined slightly different in different cases. The iteration stops when $\|\mathbf p^{k+1} - \mathbf p^{k}\|_\infty \leq 0.0005$.

\subsection{Inferring the Average Queue Length}



\begin{figure}[htbp]
	\centering
	\includegraphics{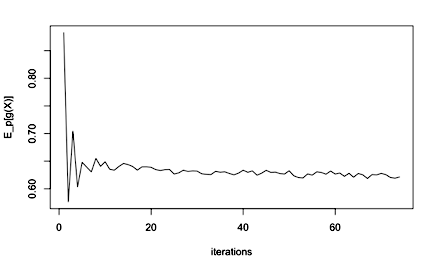}
	\includegraphics{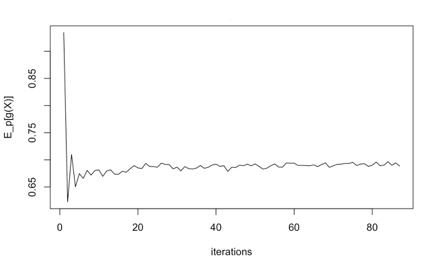}
	\caption{Objective value of the minimization (left) and maximization (right) for the expected queue length using Algorithm \ref{algo2} against the iteration number; $n=30,m=100$; true service time distribution is exponential}\label{trend_queue_100_30}
\end{figure}

We first consider inferring the average queue length $E_{P_X}[g(\mathbf X)]$, and consider a small output data size $n = 30$ for the average waiting time. In this setting, the true service time distribution is set as exponential with rate $1.2$. We generate the input support points with a lognormal distribution with parameter $\mu = 0$ and standard deviation $\sigma = 1$. In light of Theorem \ref{main guarantee}, we choose $m = 100$ to make $m$ bigger than $n$. Figure \ref{trend_queue_100_30} shows the trend of the objective value $E_{\mathbf p}[g(\textbf{X})]$ when we apply Algorithm \ref{algo2} to the max and the min problems. The algorithm appears to converge fairly quickly (within about 10 iterations). The jitter of the trend is due to the evaluation of the objective value, for each of whom we use $100,000$ simulation runs. The minimization stops at $0.622$ and the maximization stops at $0.688$ according to our stopping criterion described above. This gives us an interval $[0.622,0.688]$. The true value in this case is $E_{\mathbf p}[g(\textbf{X})] = 0.636$ (from running 1 million simulation using the true service time distribution), thus demonstrating that the confidence interval we obtained covers the truth. Moreover, the interval we obtained is encouragingly tight.

\begin{figure}[htbp]
	\centering
	\includegraphics{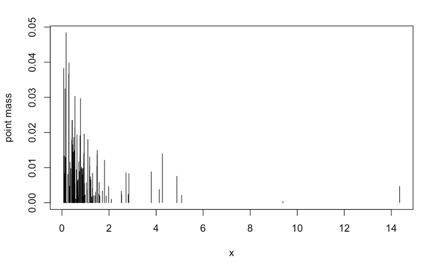}
	\includegraphics{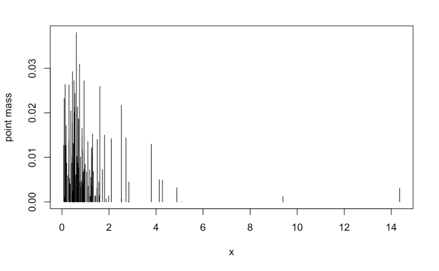}
	\caption{Minimal (left) and maximal (right) distribution of the service time for bounding the expected queue length; $n=30,m=100$; true service time distribution is exponential}\label{dist_queue_100_30}
\end{figure}

We also investigate the shape of the input distribution when the algorithm stops. This is shown in Figure \ref{dist_queue_100_30}. We observe that both the obtained maximal and minimal distributions place more masses on the lower value than the upper, roughly following the true exponential distribution. We should mention, however, that the shapes of the obtained optimal distributions are not indicative of the performance of our method, as the latter intends to compute valid bounds for a target quantity, namely the average queue length in this example, instead of direct recovery of the input distribution. The shapes in Figure \ref{dist_queue_100_30} should be interpreted as the worst-case distributions that give rise to the lower and upper bounds for the queue length. The resemblance of these distributions to the true one leads us to conjecture that the service time distribution could be close to identifiable with the waiting time data. 

\begin{figure}[htbp]
	\centering
	\includegraphics{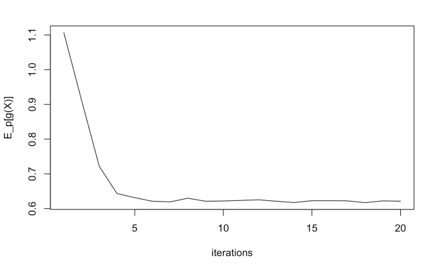}
	\includegraphics{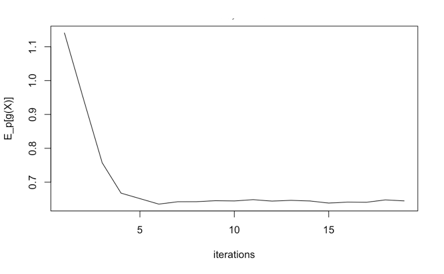}
	\caption{Objective value of the minimization (left) and maximization (right) for the expected queue length using Algorithm \ref{algo2} against the iteration number; $n=30,m=200$; true service time distribution is exponential}\label{trend_queue_200_30}
\end{figure}

\begin{figure}[htbp]
	\centering
	\includegraphics{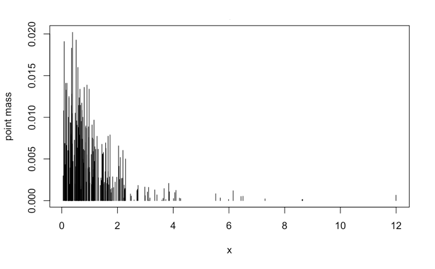}
	\includegraphics{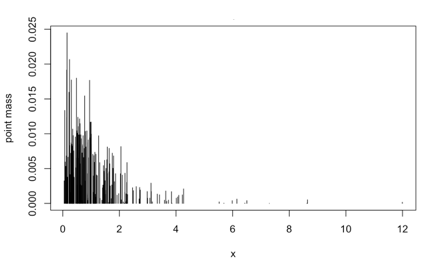}
	\caption{Minimal (left) and maximal (right) distribution of the service time for bounding the expected queue length; $n=30,m=200$; true service time distribution is exponential}\label{dist_queue_200_30}
\end{figure}

Next we increase our support size $m$ to $200$, keeping the output data size $n$ fixed at $30$. Like the previous case, we show the trend of the objective value as the algorithm progresses, in Figure \ref{trend_queue_200_30}. Compared to the case $m=100$, the algorithm appears to stabilize faster, at around 5 iteration, and exhibit a more monotonic trend (which could be due to our initialization). The minimization stops at $0.622$ and the maximization stops at $0.647$. This gives us an interval $[0.622,0.647]$ which again covers the true value $0.636$, and is shorter than the one obtained when $m=100$. Finally, The obtained maximal and minimal distributions, shown in Figure \ref{dist_queue_200_30}, show a pattern even closer to the exponential distribution.


We increase the support size $m$ further to $300$ or the data size $n$ to $100$. Table \ref{table1} shows the obtained optimal values. These runs provide valid lower and upper bounds for the true value $0.636$, except when $m=300$ and $n=30$ that misses marginally. The interval lengths do not seem to vary much; all are around $0.03-0.06$.


\begin{table}[htbp]
	\footnotesize
	\centering
	
	\begin{tabular}{cccccccccc}
		\toprule
		$m$ & $n$ & min value & max value\\
		\midrule
	   100 &	30 &	0.622 &	0.688 \\
       200 &	30 &	0.622 &	0.647 \\
       300 &    30 &	0.593 &	0.629 \\
       100 &	100 &	0.627 &	0.652 \\
		\bottomrule
	\end{tabular}%
    \caption{Optimal values for bounding the expected queue length under different combinations of $n$ and $m$; true service time distribution is exponential}
	\label{table1}%
\end{table}%


The selection of $a, b, c$ in $\gamma^k$, $\beta^k$, $\lambda^k$ depends on $m$ and $n$. We have selected $a = 0.2$  when $m = 100$ and $n = 30$, $a = 0.1$ and  $0.075$ when $m = 200$ and $300$ while $n=30$, and $a = 0.1$ when $m =  100$ and $n = 100$. We always choose $b = 0.2$ and $c = 1$. These choices appear to work well.
Regarding running times, when $m = 100$ and $n = 30$, each iteration takes about 40 seconds. The running time seems to increase linearly as $m$ and $n$ increase.

\begin{figure}[htbp]
	\centering
	\includegraphics[width = .6\textwidth]{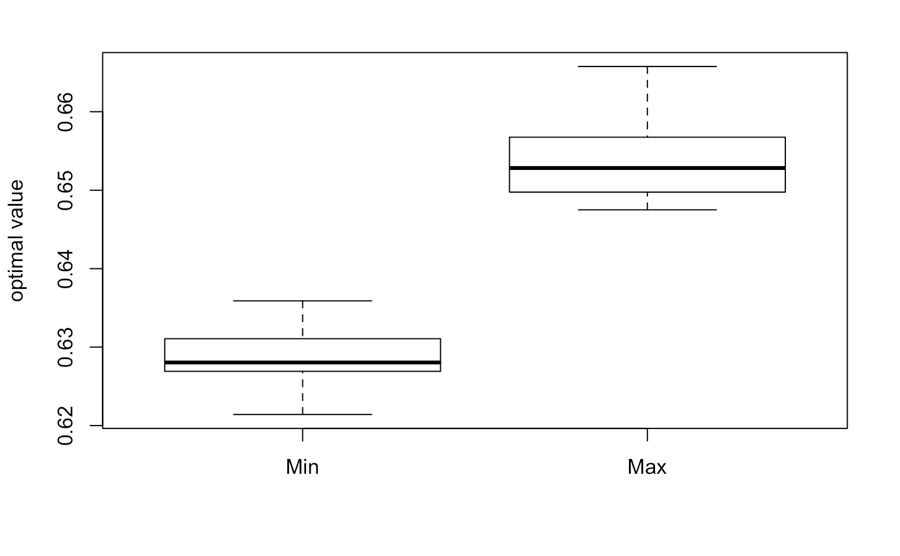}
	\caption{Minimum and maximum values for the expected queue length under different initializations; $n=30,m=100$; true service time distribution is exponential}\label{different initial boxplot}
\end{figure}

Next we check how the initialization of the probability weights in the algorithm affects the obtained optimal values. This is especially important since our algorithm is only guaranteed local convergence. We randomly generate 34 initial distributions of $\mathbf p$ from a Dirichlet distribution to run the algorithm. Figure \ref{different initial boxplot} shows the boxplot of the obtained optimal values under different initial distributions. The minimum value varies from $0.621$ to $0.635$, whereas the maximum value varies from $0.648$ to $0.665$. The differences among the initial distributions seem to be quite small compared to the gap between the minimum and maximum values, and the true value $0.636$ is always covered. This shows that the algorithm tends to converge to the same optimal solution or solutions that have similar objective values.

We then test the coverage of our obtained bounds. For this, we repeatedly sample new output data set of size $n=30$ for $100$ times. For each data set, we generate new support points of size $m=100$. Then we run Algorithm \ref{algo2}. Out of $100$ intervals we obtained, five of them cover the true expected queue length. This gives us a $95\%$ confidence interval for the coverage probability $[0.91,0.99]$, which is consistent with the theoretical guarantee provided by Theorem \ref{main guarantee}.


We have also tested the use of randomized stochastic projected gradient (RSPG), proposed by \cite{ghadimi2016mini}, that has been shown to perform well theoretically and empirically for problems with non-convex stochastic objectives. Specifically, we adapt the algorithm in Section 4.1 and 4.2 of \cite{ghadimi2016mini} heuristically for the current problem we face that has stochastic non-convex constraints. Algorithm \ref{RSPG} in the appendix shows the adaptation of a single run procedure, and Algorithm \ref{2-RSPG-V} shows the adaptation of a post-optimization step to boost the final performance. In our algorithmic specification, we choose $N=30$, $S=5$, $M=500$, $M'=500$, $\bar\gamma=0.03$, and we fix $\lambda$ at $0.03$. We run Algorithm \ref{2-RSPG-V} for two realizations of data and support generation when the true service time distribution is exponential, with $n=30$ and $m=100$. For each realization, we also run Algorithm \ref{algo2} for comparison. For the first realization, we obtained $[0.622, 0.640]$ using RSPG, compared with $[0.626,0.658]$ using Algorithm \ref{algo2}. For the second realization, we obtained  $[0.616,0.644]$ using RSPG, compared with $[0.621,0.660]$ using Algorithm \ref{algo2}. The RSPG thus appears to perform very similarly as our procedure, at least for this particular setup (which shows that RSPG could be an alternative for future investigation).

We test the sensitivity of the algorithm with respect to the bounds in the constraints provided by the KS statistic. More concretely, in Algorithm \ref{algo2}, we increase the number $q_{1-\alpha}/\sqrt n$ in the constraint interval by a small $\delta$. Table \ref{change_of_delta} shows that the obtained bounds are quite stable and do not show significant changes.

\begin{table}[htbp]
	\footnotesize
	\centering
	\begin{tabular}{cccccccccc}
		\toprule
		perturbation size & min value & max value\\
		\midrule
		0.01	& 0.625	& 0.649 \\
		0.02	& 0.628	& 0.649 \\
		0.03	& 0.624	& 0.643 \\
		0.05	& 0.621	& 0.646 \\
		\bottomrule
	\end{tabular}%
	\caption{Effect on optimal values for bounding the expected queue length when perturbing the interval in the optimization constraint; $n=30,m=100$; true service time distribution is exponential}
	\label{change_of_delta}%
\end{table}%

\begin{figure}[htbp]
	\centering
	\includegraphics[width = .6\textwidth]{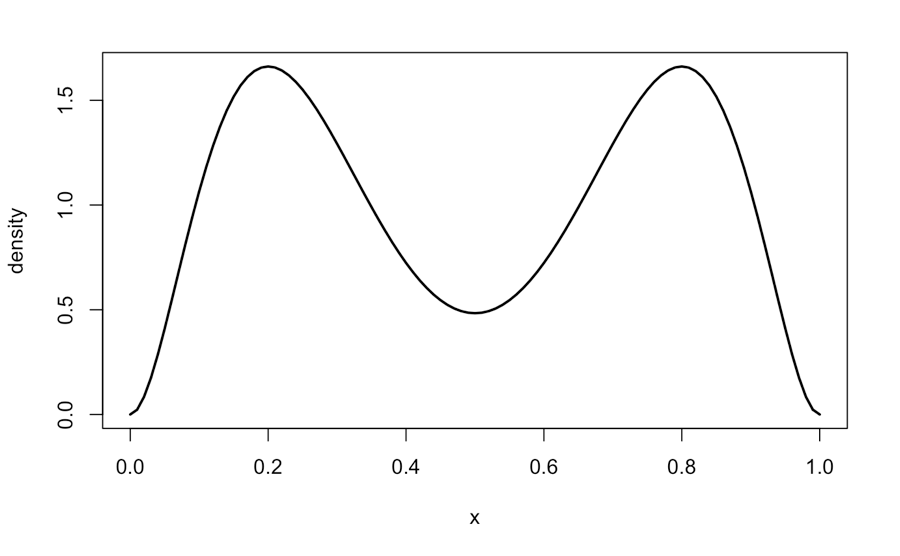}
	\caption{Density of a mixture of two beta distributions }\label{density of bimodal}
\end{figure}

Finally, we test with a more ``challenging" service time distribution that is an equally weighted mixture of two beta distributions with parameters $\alpha = 9,\beta = 3$ and $\alpha = 3,\beta = 9$. This bimodal distribution has highest masses around $0.2$ and $0.8$, with a shape shown in Figure \ref{density of bimodal}.


We consider the setting with $n=50$ output observations. We randomly select $m=100$ input support points from uniform distribution in $[0,1]$, and run Algorithm \ref{algo2}, using the same specifications as in the previous setup. The minimization stops at the value $0.242$ and the maximization stops at $0.284$. These cover the true value $0.274$ (from running 1 million simulation using the true service time distribution). Thus our method appears to continue working in this case. 

\begin{figure}[htbp]
	\centering
	\includegraphics{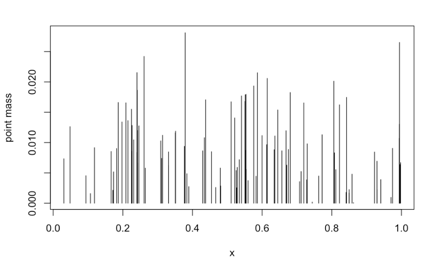}
	\includegraphics{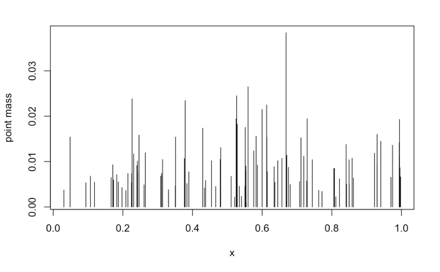}
	\caption{Minimal (left) and maximal (right) distribution of the service time for bounding the expected queue length; $n=50,m=100$; true service time distribution is mixture of betas}\label{dist_bimodal}
\end{figure}

Figure \ref{dist_bimodal} shows the minimal and maximal distributions from Algorithm \ref{algo2}. The distributions are quite spread out throughout the support, though the minimal distribution appears to have a noisy bimodal pattern. As we have discussed before, the shapes of these distributions should be interpreted as the worst-case distributions giving rise to the bounds, but are not indicative of the performance of our approach.


\subsection{Inferring the Input Distribution Function}
We now consider inferring the distribution function of the service time, i.e.,  $P_X(X\leq a)$ for a range of values $a$. We first use a true service time distribution that is exponential with rate $1.2$. We consider a collection of $n=50$ observations from the average waiting time. We randomly generate $m =100$ support points from a lognormal distribution with $\mu=0$ and $\sigma^2=1$. We use Algorithm \ref{algo2} with parameters $\gamma^k =  0.1/k^{0.8}$, $\beta^k = 0.1/k^{0.5}$, $\lambda^k = 1/k^{0.2}$. 

\begin{figure}[htbp]
	\centering
	\includegraphics[width = .6\textwidth]{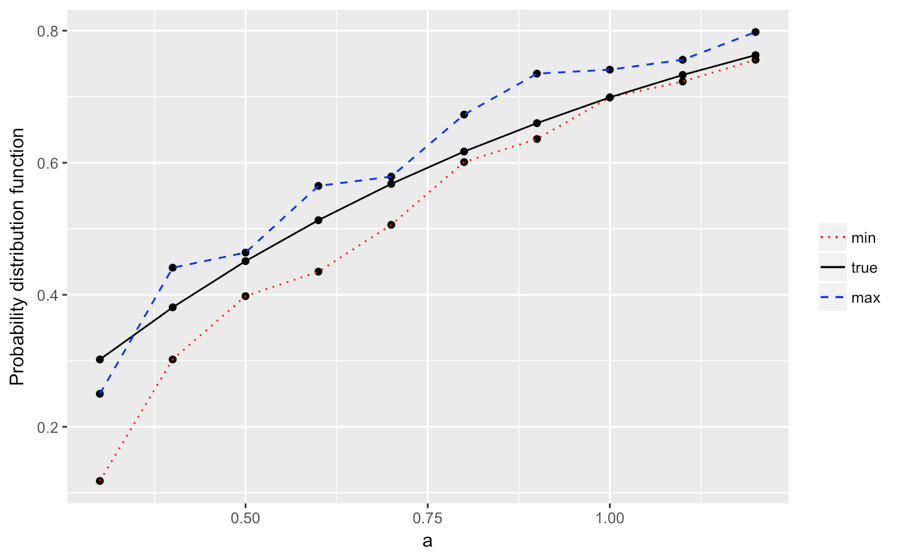}
	\caption{Bounds and true distribution function values for the service time, when the true service time distribution is exponential; $n=50,m=100$}\label{expo_ggplot}
\end{figure}

Table \ref{expo_P(X<=a)} shows the obtained maximum and minimum values compared with the true distribution function evaluated at values $a$ ranging from $0.3$ to $1.2$. Figure \ref{expo_ggplot} further plots the trends of these values.  The dashed lines represent the maximum and minimum values, and the solid line represents the true values. Note that Proposition \ref{extended guarantee}, and the analogous extension of Theorem \ref{main guarantee} to multiple objective functions discussed at the end of Section \ref{sec:discretization}, allow us to compute the bounds for different $a$ values simultaneously with little sacrifice of statistical accuracy. In Table \ref{expo_P(X<=a)} and Figure \ref{expo_ggplot}, the obtained optimal values cover the truth at all points except the leftmost $a=0.3$. This could be due to the challenge in inferring the tail (either left or right), stemming from perhaps the observed output we use (i.e., the waiting time) or the statistic we use to form our uncertainty set (i.e., the KS-statistic, which is known to not capture well the tail region of a distribution).

\begin{table}[htbp]
	\footnotesize
	\centering
	
	\begin{tabular}{cccccccccc}
		\toprule
		$a$ & min value & max value & true value\\
		\midrule
		0.3	& 0.118 	& 0.250 	& 0.302 \\
0.4	& 0.302 	& 0.441 	& 0.381 \\
0.5	& 0.398 	& 0.464 	& 0.451 \\
0.6	& 0.435 	& 0.565 	& 0.513 \\
0.7	& 0.506 	& 0.579 	& 0.568 \\
0.8	& 0.601 	& 0.673 	& 0.617 \\
0.9	& 0.636 	& 0.735 	& 0.660 \\
1	& 0.699 	& 0.741 	& 0.699 \\
1.1	& 0.723 	& 0.756 	& 0.733 \\
1.2	& 0.756 	& 0.798 	& 0.763 \\
		\bottomrule
	\end{tabular}%
	\caption{Minimum, maximum and true values of the distribution function $P_X(X\leq a)$ of the service time across $a$; $n=50,m=100$; true service time distribution that is exponential}
	\label{expo_P(X<=a)}%
\end{table}%

\begin{figure}[htbp]
	\centering
	\includegraphics[width = .4\textwidth]{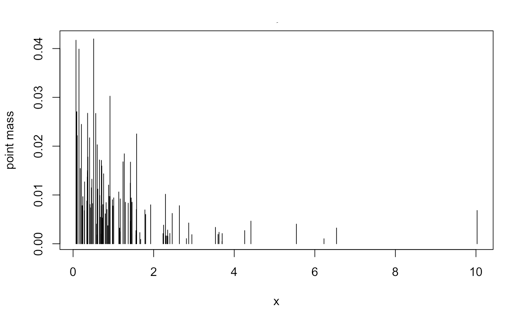}\includegraphics[width = .4\textwidth]{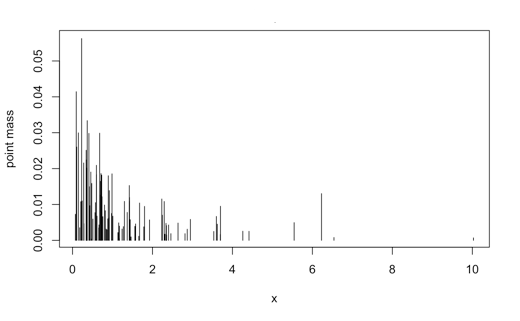}
	\caption{Minimal (left) and maximal (right) distribution of the service time for bounding $P_X(X\leq0.5)$, when the true service time distribution is exponential; $n=50,m=100$}\label{dist_exp_a = 0.5}
\end{figure}

Figure \ref{dist_exp_a = 0.5} shows the minimal and maximal distributions for bounding $P_X(X\leq0.5)$ when the algorithm terminates. We see that the shapes of both distributions resemble exponential, hinting that the service time distribution is close to identifiable in this case.

Next, we investigate the case when the true service time distribution is a mixture of two beta distributions with parameters $\alpha = 9,\beta = 3$ and $\alpha = 3,\beta = 9$. We consider a collection of $n=50$ observations from the average waiting time. We randomly generate $m =100$ support points from a uniform distribution on $[0,1]$. 

\begin{figure}[htbp]
	\centering
	\includegraphics[width = .6\textwidth]{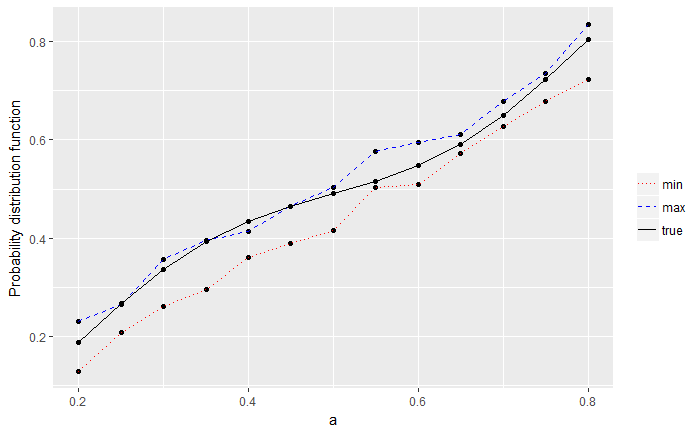}
	\caption{Bounds and true distribution function values for the service time, when the true service time distribution is mixture of betas; $n=50,m=100$}\label{ggplot}
\end{figure}

Like in the previous case, Table \ref{P(x<=a)} shows the maximum and minimum values from Algorithm \ref{algo2}, against the true values of $P_X(X\leq a)$ at different $a$ values. Figure \ref{ggplot} further plots the trends of these values. Here, the obtained optimal values all cover the truth except at $a = 0.35$. The latter could be attributed to the statistical noise when running the many optimization procedures. The point $a=0.35$ is also one that could be ``difficult" to infer intuitively, as it is in between the two modes. Nonetheless, our procedure appears to be reliable in general in bounding the distribution function across the domain of the service time.
\begin{table}[htbp]
	\footnotesize
	\centering
	
	\begin{tabular}{cccccccccc}
		\toprule
		$a$ & min value & max value & true value\\
		\midrule
		0.2	& 0.129 	& 0.231 	& 0.188 \\
0.25	& 0.208 	& 0.266 	& 0.267 \\
0.3	& 0.262 	& 0.358 	& 0.337 \\
0.35	& 0.296 	& 0.395 	& 0.393 \\
0.4	& 0.362 	& 0.413 	& 0.435 \\
0.45	& 0.389 	& 0.464 	& 0.466 \\
0.5	& 0.416 	& 0.503 	& 0.491 \\
0.55	& 0.504 	& 0.577 	& 0.516 \\
0.6	& 0.509 	& 0.594 	& 0.548 \\
0.65	& 0.573 	& 0.611 	& 0.591 \\
0.7	& 0.628 	& 0.679 	& 0.649 \\
0.75	& 0.678 	& 0.736 	& 0.722 \\
0.8	& 0.724 	& 0.834 	& 0.805 \\
		\bottomrule
	\end{tabular}%
	\caption{Minimum, maximum and true values of the distribution function $P_X(X\leq a)$ of the service time across $a$, under a true service time distribution that is mixture of betas; $n=50,m=100$}
	\label{P(x<=a)}%
\end{table}%

\begin{figure}[htbp]
	\centering
	\includegraphics[width = .4\textwidth]{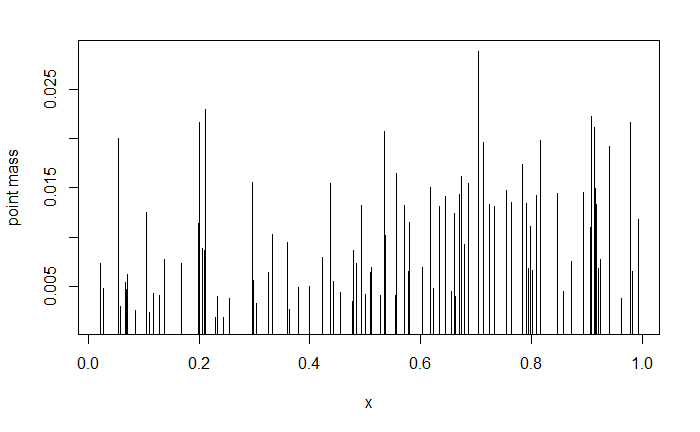}\includegraphics[width = .4\textwidth]{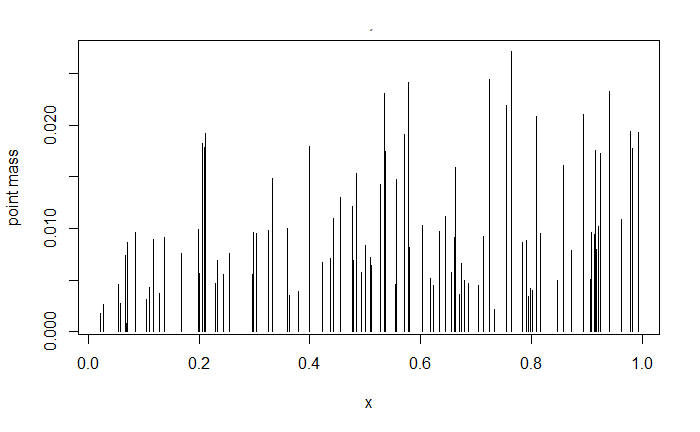}
	\caption{Minimal (left) and maximal (right) distribution of the service time for bounding $P_X(X\leq0.5)$, when the true service time distribution is mixture of betas; $n=50,m=100$}\label{dist_bimodal_a = 0.5}
\end{figure}

Figure \ref{dist_bimodal_a = 0.5} shows the minimal and maximal distributions for bounding $P_X(X\leq0.5)$ when the algorithm terminates. The shapes of these distributions are now considerably noisier than the exponential case in Figure \ref{dist_exp_a = 0.5}. Nonetheless, there is a rough bimodal pattern (around $0.2$ and $0.7$). 


\section{Conclusion}\label{sec:discussion}
We have studied an optimization-based framework to calibrate input quantities in stochastic simulation with only the availability of output data. Our approach uses an output-level uncertainty set, inspired by the DRO literature, to represent the statistical noise of the output data. By expressing the output distribution in terms of a simulable map of the input distribution, we can set up optimization programs cast over the input distribution that infers valid confidence bounds on the input quantities of interest. 

We propose in particular an output-level uncertainty set based on the KS statistic, which exhibits advantages in computation (thanks to reformulation) and statistical accuracy (thanks to a controllable discretization scale needed to retain the confidence guarantee). We have shown these advantages via looking at the complexity of the resulting constraints and invoking the empirical process theory for $U$-statistics. We also study a stochastic quadratic penalty method to solve the resulting optimization problems, including a convergence analysis that informs the suitable tuning of the parameters. Our numerical results demonstrate how our method could provide valid bounds for  input quantities such as the input distribution function and other performance measures that rely on the input.

\ACKNOWLEDGMENT{A preliminary conference version of this work has appeared in \cite{goeva2014reconstructing}. We gratefully acknowledge support from the National Science Foundation under grants CMMI-1542020, CMMI-1523453 and CAREER CMMI-1653339.  We also thank Peter Haas for suggesting the use of quantile-based moments, and Russell Barton, Shane Henderson and Barry Nelson for other helpful suggestions.}


\bibliographystyle{informs2014} 
\bibliography{ms2015} 

\ECSwitch


\ECHead{Supplementary Materials}

\section{Proofs and Additional Results for Section \ref{sec:formulation}}
\proof{Proof of Proposition \ref{basic guarantee}. }
Note that if $P_Y^0=\gamma(P_X^0)\in\mathcal U$, then $P_X^0$ must be a feasible solution for programs \eqref{RO max} and \eqref{RO min}, and consequently $\underline Z\leq\psi(P_X^0)\leq\overline Z$. This implies that
$$\mathbb P_D(\underline Z\leq\psi(P_X^0)\leq\overline Z)\geq\mathbb P_D(P_Y^0\in\mathcal U)$$
concluding the proposition.\Halmos
\endproof

\begin{proposition}
Let $P_X^0$ and $P_Y^0$ be the true input and output distributions. Consider a collection of quantities $\psi_l(P_X),l=1,\ldots,L$ and the collection of optimization programs
\begin{equation}
\begin{array}{ll}
\max&\psi_l(P_X)\\
\text{subject to}&P_Y\in\mathcal U
\end{array}\label{extended RO max}
\end{equation}
and
\begin{equation}
\begin{array}{ll}
\min&\psi_l(P_X)\\
\text{subject to}&P_Y\in\mathcal U
\end{array}\label{extended RO min}
\end{equation}
for $l=1,\ldots,L$. Suppose $\mathcal U$ is a confidence region for $P_Y^0$, i.e., 
$$\mathbb P_D(P_Y^0\in\mathcal U)=1-\alpha$$
where $\mathbb P_D(\cdot)$ denotes the probability with respect to the data $D$. Let $\overline Z_l,\underline Z_l,l=1,\ldots,L$ be the set of optimal values of \eqref{extended RO max} and \eqref{extended RO min} respectively. Then we have
$$\mathbb P_D(\underline Z_l\leq\psi_l(P_X^0)\leq\overline Z_l,l=1,\ldots,L)\geq1-\alpha$$
Similar statements hold if the confidence is approximate, i.e., if
$$\liminf_{n\to\infty}\mathbb P_D(P_Y^0\in\mathcal U)\geq1-\alpha$$
then
$$\liminf_{n\to\infty}\mathbb P_D(\underline Z_l\leq\psi_l(P_X^0)\leq\overline Z_l,l=1,\ldots,L)\geq1-\alpha$$
\label{extended guarantee}
\end{proposition}

\proof{Proof of Proposition \ref{extended guarantee}.}
The proof follows similarly from that of Proposition \ref{basic guarantee}. If $P_Y^0=\gamma(P_X^0)\in\mathcal U$, then $P_X^0$ must be a feasible solution for programs \eqref{extended RO max} and \eqref{extended RO min}, and consequently $\underline Z_l\leq\psi_l(P_X^0)\leq\overline Z_l$, simultaneously for $l=1,\ldots,L$. Therefore
$$\mathbb P_D(\underline Z_l\leq\psi_l(P_X^0)\leq\overline Z_l,l=1,\ldots,L)\geq\mathbb P_D(P_Y^0\in\mathcal U)$$
This concludes the proposition.
\Halmos
\endproof

\section{Proofs for Section \ref{sec:KS}}\label{sec:KS proofs}

\proof{Proof of Theorem \ref{KS guarantee}. }
Note that the first constraints in \eqref{KS RO max direct} and \eqref{KS RO min direct} can be readily replaced by $P_Y\in\mathcal U$ for $\mathcal U$ defined in \eqref{uncertainty set}. We have $\lim_{n\to\infty}\mathbb P_D(P_Y^0\in\mathcal U)=1-\alpha$, where $P_Y^0$  is the true output distribution, as a consequence of the KS statistic asymptotic. By using Proposition \ref{basic guarantee}, we arrive at the guarantee \eqref{KS confidence guarantee}.

The second conclusion comes from a reformulation of \eqref{uncertainty set}. Note that
$$\|F_Y-\hat F_Y\|_\infty\leq\frac{q_{1-\alpha}}{\sqrt n}$$
is equivalent to
$$\sup_{y\in\mathbb R}|E_{P_X}[I(h(\mathbf X)\leq y)]-\hat F_Y(y)|\leq\frac{q_{1-\alpha}}{\sqrt n}$$
By the monotonicity of distribution functions, this is further equivalent to the set of constraints
\begin{equation}
\hat F_Y(y_j+)-\frac{q_{1-\alpha}}{\sqrt n}\leq E_{P_X}[I(h(\mathbf X)\leq y_j)]\leq\hat F_Y(y_j-)+\frac{q_{1-\alpha}}{\sqrt n},j=1,\ldots,n\label{KS}
\end{equation}
which gives \eqref{KS RO max} and \eqref{KS RO min}. 
\Halmos
\endproof

\proof{Proof of Theorem \ref{main guarantee}.}
We will show the conclusion when \eqref{KS RO max discretized} and \eqref{KS RO min discretized} are replaced by 
\begin{equation}
\begin{array}{ll}
\max&\psi(P_X)\\
\text{subject to}&\|E_{P_X}[I(h(\mathbf X)\leq \cdot)]-\hat F_Y(\cdot)\|_\infty\leq\frac{q_{1-\alpha}}{\sqrt n}\\
&P_X\in\hat{\mathcal P}_X
\end{array}\label{KS RO max discretized direct}
\end{equation}
and
\begin{equation}
\begin{array}{ll}
\min&\psi(P_X)\\
\text{subject to}&\|E_{P_X}[I(h(\mathbf X)\leq \cdot)]-\hat F_Y(\cdot)\|_\infty\leq\frac{q_{1-\alpha}}{\sqrt n}\\
&P_X\in\hat{\mathcal P}_X
\end{array}\label{KS RO min discretized direct}
\end{equation}
Then by the assumption that the true output distribution is continuous and that $\mathbb P(\text{for any\ }P_X\in\hat{\mathcal P}_X, \text{supp}(\gamma(P_X))\cap\{y_j\}_{j=1,\ldots,n}\neq\emptyset)=0$, we can use the same argument as in Theorem \ref{KS guarantee} to deduce that the constraints in \eqref{KS RO max discretized direct} and \eqref{KS RO min discretized direct} are equivalent to those in \eqref{KS RO max discretized} and \eqref{KS RO min discretized} with probability 1, from which we conclude the theorem. 

Denote $L=dP_X^0/dQ$. Denote $\hat P_X(\cdot)$ as the empirical distribution on $\{z_j\}$ given by
$$\hat P_X(\cdot)=\frac{1}{m}\sum_{i=1}^m\delta_{z_i}(\cdot)$$
where $\delta_{z_j}(\cdot)$ is the delta mass on $z_j$. Consider 
$$\tilde P_X(\cdot)=\sum_{i=1}^m\frac{L(z_i)}{\sum_{j=1}^mL(z_j)}\delta_{z_i}(\cdot)$$
i.e., $\tilde P_X$ is a discrete probability distribution with mass $L(z_i)/\sum_{j=1}^mL(z_j)$ on each generated support point $z_i$ of $X$. Consider, for any $y\in\mathbb R$,
\begin{eqnarray}
&&E_{\tilde P_X}[I(h(\mathbf X)\leq y)]-E_{P_X^0}[I(h(\mathbf X)\leq y)]\notag\\
&=&\left(E_{\tilde P_X}[I(h(\mathbf X)\leq y)]-E_{\bar P_X}[I(h(\mathbf X)\leq y)]\right)+\left(E_{\bar P_X}[I(h(\mathbf X)\leq y)]-E_{P_X^0}[I(h(\mathbf X)\leq y)]\right)\label{interim update1}
\end{eqnarray}
where $\bar P_X(\cdot)$ is a measure (not necessarily a probability) given by
$$\bar P_X(\cdot)=\frac{1}{m}\sum_{i=1}^mL(z_i)\delta_{z_i}(\cdot)$$
and the expectation $E_{\bar P_X}[I(h(\mathbf X)\leq y)]$ is defined in a general sense as the $T$-fold integral of $I(h(\mathbf X)\leq y)$ with respect to $\bar P_X$. We consider both terms in \eqref{interim update1}. Writing $\mathbf x=(x_1,\ldots,x_T)$, we can write the first term as
\begin{eqnarray}
&&\int\cdots\int I(h(\mathbf x)\leq y)\prod_{t=1}^Td\tilde P_X(x_t)-\int\cdots\int I(h(\mathbf x)\leq y)\prod_{t=1}^Td\bar P_X(x_t)\notag\\
&=&\int\cdots\int I(h(\mathbf x)\leq y)\frac{\prod_{t=1}^TL(x_t)d\hat P_X(x_t)}{\left(\frac{1}{m}\sum_{j=1}^mL(z_j)\right)^T}-\int\cdots\int I(h(\mathbf x)\leq y)\prod_{t=1}^TL(x_t)d\hat P_X(x_t)\notag\\
&=&\int\cdots\int I(h(\mathbf x)\leq y)\prod_{t=1}^TL(x_t)d\hat P_X(x_t)\left(\frac{1}{\left(\frac{1}{m}\sum_{j=1}^mL(z_j)\right)^T}-1\right)\label{interim update2}
\end{eqnarray}
Since $Var_Q(L)<\infty$, and $E_Q[L]=1$ by the definition of likelihood ratio, we have $\sqrt m((1/m)\sum_{j=1}^mL(z_j)-1)\Rightarrow N(0,Var_Q(L))$ by the central limit theorem. By using the delta method (Chapter 3 in \cite{serfling2009approximation}), we also have $\sqrt m(1/((1/m)\sum_{j=1}^mL(z_j))^T-1)\Rightarrow N(0,T^2Var_Q(L))$.

Moreover, $\int\cdots\int I(h(\mathbf x)\leq y)\prod_{t=1}^TL(x_t)d\hat P_X(x_t)$ is bounded by $C^T$ since $\|L\|_\infty\leq C$. Hence \eqref{interim update2} satisfies
\begin{align}
\sup_{y\in\mathbb R}\left|\int\cdots\int I(h(\mathbf x)\leq y)\prod_{t=1}^TL(x_t)d\hat P_X(x_t)\left(\frac{1}{\left(\frac{1}{m}\sum_{j=1}^mL(z_j)\right)^T}-1\right)\right|&\leq C^T\left|\frac{1}{\left(\frac{1}{m}\sum_{j=1}^mL(z_j)\right)^T}-1\right|\notag\\
&=O_p\left(\frac{1}{\sqrt m}\right)\label{interim updated3}
\end{align}

Now consider the second term in \eqref{interim update1}. We have
\begin{eqnarray*}
&&E_{\bar P_X}[I(h(\mathbf X)\leq y)]-E_{P_X^0}[I(h(\mathbf X)\leq y)]\\
&=&E_{\hat P_X}\left[I(h(\mathbf X)\leq y)\prod_{t=1}^TL(X_t)\right]-E_{Q}\left[I(h(\mathbf X)\leq y)\prod_{t=1}^TL(X_t)\right]
\end{eqnarray*}
by the definition of $\bar P_X$, $\hat P_X$ and $Q$. Note that $E_{\hat P_X}\left[I(h(\mathbf X)\leq y)\prod_{t=1}^TL(X_t)\right]$ is the average, over all possible selections with replacement of $x_1,\ldots,x_T$ drawn from $\{z_i\}_{i=1,\ldots,m}$, of the multilinear form $I(h(\mathbf x)\leq y)\prod_{t=1}^TL(x_t)$. This is equivalent to the $V$-statistic (\cite{serfling2009approximation} Chapter 5) with kernel $I(h(\mathbf x)\leq y)\prod_{t=1}^TL(x_t)$. 

Define $\mathcal F$ as the class of functions from $\mathcal X^T$ to $\mathbb R$ given by $\mathcal F=\{I(h(\mathbf x)\leq y)\prod_{t=1}^TL(x_t):y\in\mathbb R\}$. Since $I(h(\mathbf x)\leq y)\prod_{t=1}^TL(x_t)$ is non-decreasing fixing each $\mathbf x$, and the envelope of $\mathcal F$, namely $\sup_{y\in\mathbb R}I(h(\mathbf x)\leq y)\prod_{t=1}^TL(x_t)$, is bounded by $C^T$ a.s., Problem 3 in Chapter 2.7 of \cite{van1996weak} (Theorem \ref{bracketing} in the appendix) implies that $\mathcal F$ has a polynomial bracketing number. Therefore, Theorem 4.10 in \cite{arcones1993limit} (Theorem \ref{CLT U} in the appendix; see also the discussion after therein) concludes the convergence
$$\left\{\sqrt m\left(U_m^T\left[I(h(\mathbf X)\leq y)\prod_{t=1}^TL(X_t)\right]-E_{Q}\left[I(h(\mathbf X)\leq y)\prod_{t=1}^TL(X_t)\right]\right)\right\}_{y\in\mathbb R}\Rightarrow \{\mathbb G(y)\}_{y\in\mathbb R}\text{\ in\ }\ell^\infty(\mathcal F)$$
where $U_m^T$ is the $U$-operator defined in \eqref{u operator} generated from $P_X$, and $\mathbb G$ is a Gaussian process defined as in \eqref{CLT process}.


Following the argument of the lemma in Section 5.7.3 in \cite{serfling2009approximation}, we can write the difference between the $U$-statistic, denoted for simplicity $U_m=U_m^T\left[I(h(\mathbf X)\leq y)\prod_{t=1}^TL(X_t)\right]$, and the $V$-statistic, denoted $V_m=E_{\hat P_X}\left[I(h(\mathbf X)\leq y)\prod_{t=1}^TL(X_t)\right]$, as
$$m^T(U_m-V_m)=(m^T-m_{(T)})(U_m-W_m)$$
where $m_{(T)}=m(m-1)\cdots(m-T+1)$, and $W_m$ is the average of all $I(h(\mathbf x)\leq y)\prod_{t=1}^TL(x_t)$ where $\mathbf x$ are drawn from $\{z_i\}_{i=1,\ldots,m}$ with replacement and at least one overlapping selection. Following \cite{serfling2009approximation}, we can verify $m^T-m_{(T)}=O(m^{T-1})$, and since $\|L\|_\infty\leq C$, we have $U_m-W_m$ bounded a.s. Hence $E\sup_{t\in\mathbb R}|U_m-V_m|^2=O(1/m^2)$, and so $\sup_{t\in\mathbb R}|U_m-V_m|=O_p(1/m)$.

Therefore, we write
\begin{eqnarray*}
&&\sqrt m\left(E_{\hat P_X}\left[I(h(\mathbf X)\leq y)\prod_{t=1}^TL(X_t)\right]-E_{Q}\left[I(h(\mathbf X)\leq y)\prod_{t=1}^TL(X_t)\right]\right)\\
&=&\sqrt m\left(E_{\hat P_X}\left[I(h(\mathbf X)\leq y)\prod_{t=1}^TL(X_t)\right]-U_m^T\left[I(h(\mathbf X)\leq y)\prod_{t=1}^TL(X_t)\right]\right){}\\
&&{}+\sqrt m\left(U_m^T\left[I(h(\mathbf X)\leq y)\prod_{t=1}^TL(X_t)\right]-E_{Q}\left[I(h(\mathbf X)\leq y)\prod_{t=1}^TL(X_t)\right]\right)
\end{eqnarray*}
where $\sqrt m\left(E_{\hat P_X}\left[I(h(\mathbf X)\leq y)\prod_{t=1}^TL(X_t)\right]-U_m^T\left[I(h(\mathbf X)\leq y)\prod_{t=1}^TL(X_t)\right]\right)=o_p(1)$ and $\sqrt m\left(U_m^T\left[I(h(\mathbf X)\leq y)\prod_{t=1}^TL(X_t)\right]-E_{Q}\left[I(h(\mathbf X)\leq y)\prod_{t=1}^TL(X_t)\right]\right)$ converges to a Gaussian process. 
This entails that
\begin{equation}
\sup_{y\in\mathbb R}\left|E_{\hat P_X}\left[I(h(\mathbf X)\leq y)\prod_{t=1}^TL(X_t)\right]-E_{Q}\left[I(h(\mathbf X)\leq y)\prod_{t=1}^TL(X_t)\right]\right|=O_p\left(\frac{1}{\sqrt m}\right)\label{interim updated4}
\end{equation}
From \eqref{interim update1}, and using \eqref{interim updated3} and \eqref{interim updated4}, we get
\begin{eqnarray}
&&\sup_{y\in\mathbb R}\left|E_{\tilde P_X}[I(h(\mathbf X)\leq y)]-E_{P_X^0}[I(h(\mathbf X)\leq y)]\right|\notag\\
&\leq&\sup_{y\in\mathbb R}\left|E_{\tilde P_X}[I(h(\mathbf X)\leq y)]-E_{\bar P_X}[I(h(\mathbf X)\leq y)]\right|+\sup_{y\in\mathbb R}\left|E_{\bar P_X}[I(h(\mathbf X)\leq y)]-E_{P_X^0}[I(h(\mathbf X)\leq y)]\right|\notag\\
&=&O_p\left(\frac{1}{\sqrt m}\right)\label{interim update6}
\end{eqnarray}

For the above chosen $\tilde P_X$, we now have, for any small enough $\delta>0$,
\begin{eqnarray}
&&P\left(\|E_{\tilde P_X}[I(h(\mathbf X)\leq \cdot)]-\hat F_Y(\cdot)\|_\infty\leq\frac{q_{1-\alpha}}{\sqrt n}\right)\notag\\
&\geq&P\left(\|E_{\tilde P_X}[I(h(\mathbf X)\leq \cdot)]-E_{P_X^0}[I(h(\mathbf X)\leq \cdot)]\|_\infty+\|E_{P_X^0}[I(h(\mathbf X)\leq \cdot)]-\hat F_Y(\cdot)\|_\infty\leq\frac{q_{1-\alpha}}{\sqrt n}\right)\notag\\
&\geq&P\left(\|E_{P_X^0}[I(h(\mathbf X)\leq \cdot)]-\hat F_Y(\cdot)\|_\infty\leq\frac{q_{1-\alpha}-\delta}{\sqrt n};\ \|E_{\tilde P_X}[I(h(\mathbf X)\leq \cdot)]-E_{P_X^0}[I(h(\mathbf X)\leq \cdot)]\|_\infty\leq\frac{\delta}{\sqrt n}\right)\notag\\
&\geq&P\left(\|E_{P_X^0}[I(h(\mathbf X)\leq \cdot)]-\hat F_Y(\cdot)\|_\infty\leq\frac{q_{1-\alpha}-\delta}{\sqrt n}\right)-P\left(\|E_{\tilde P_X}[I(h(\mathbf X)\leq \cdot)]-E_{P_X^0}[I(h(\mathbf X)\leq \cdot)]\|_\infty>\frac{\delta}{\sqrt n}\right)\notag\\
&\to&1-\alpha+\zeta(-\delta)\label{interim update5}
\end{eqnarray}
as $n\to\infty$ and $m/n\to\infty$, where $\zeta(\cdot)$ is a function with $\lim_{x\to0}\zeta(x)=0$ that satisfies $P(\sup_{u\in[0,1]}|BB(u)|\leq q_{1-\alpha}+\rho)=1-\alpha+\zeta(\rho)$, which exists by the continuity of the distribution of $\sup_{u\in[0,1]}BB(u)$. The convergence \eqref{interim update5} follows from the definition that $P_X^0$ is the true input distribution and hence $E_{P_X^0}[I(h(\mathbf X)\leq \cdot)]$ is the true output distribution, thus leading to $\sqrt n\|E_{P_X^0}[I(h(\mathbf X)\leq \cdot)]-\hat F_Y(\cdot)\|_\infty\Rightarrow \sup_{u\in[0,1]}|BB(u)|$. It also follows from \eqref{interim update6} so that $P\left(\|E_{\tilde P_X}[I(h(\mathbf X)\leq \cdot)]-E_{P_X^0}[I(h(\mathbf X)\leq \cdot)]\|_\infty>\frac{\delta}{\sqrt n}\right)\to0$ as $m/n\to\infty$.

Similarly, for any small enough $\delta>0$, we have
\begin{eqnarray}
&&P\left(\|E_{\tilde P_X}[I(h(\mathbf X)\leq \cdot)]-\hat F_Y(\cdot)\|_\infty\leq\frac{q_{1-\alpha}}{\sqrt n}\right)\notag\\
&\leq&P\left(\|E_{P_X^0}[I(h(\mathbf X)\leq \cdot)]-\hat F_Y(\cdot)\|_\infty-\|E_{\tilde P_X}[I(h(\mathbf X)\leq \cdot)]-E_{P_X^0}[I(h(\mathbf X)\leq \cdot)]\|_\infty\leq\frac{q_{1-\alpha}}{\sqrt n}\right)\notag\\
&\leq&P\left(\|E_{P_X^0}[I(h(\mathbf X)\leq \cdot)]-\hat F_Y(\cdot)\|_\infty\leq\frac{q_{1-\alpha}+\delta}{\sqrt n}\right)+P\left(\|E_{\tilde P_X}[I(h(\mathbf X)\leq \cdot)]-E_{P_X^0}[I(h(\mathbf X)\leq \cdot)]\|_\infty>\frac{\delta}{\sqrt n}\right)\notag\\
&\to&1-\alpha+\zeta(\delta)\label{interim update7}
\end{eqnarray}
as $n\to\infty$ and $m/n\to\infty$. Since $\delta$ is arbitrary, by combining \eqref{interim update5} and \eqref{interim update7}, we have $$P\left(\|E_{\tilde P_X}[I(h(\mathbf X)\leq \cdot)]-\hat F_Y(\cdot)\|_\infty\leq\frac{q_{1-\alpha}}{\sqrt n}\right)\to1-\alpha$$
as $n\to\infty $ and $m/n\to\infty$.

Lastly, we argue that the objective function satisfies $E_{\tilde P_X}[g(\mathbf X)]-E_{P_X^0}[g(\mathbf X)]=O_p(1/\sqrt m)$. This follows mostly as a special case of the arguments above in showing $\sup_{y\in\mathbb R}|E_{\tilde P_X}[I(h(\mathbf X)\leq y)]-E_{P_X^0}[I(h(\mathbf X)\leq y)]|=O_p(1/\sqrt m)$, by simply replacing $I(h(\mathbf X\leq y)$ with $g(\mathbf X)$ and without considering the uniformity over $y\in\mathbb R$. More precisely, we have
\begin{eqnarray}
&&E_{\tilde P_X}[g(\mathbf X)]-E_{P_X^0}[g(\mathbf X)]\notag\\
&=&\left(E_{\tilde P_X}[g(\mathbf X)]-E_{\bar P_X}[g(\mathbf X)]\right)+\left(E_{\bar P_X}[g(\mathbf X)]-E_{P_X^0}[g(\mathbf X)]\right)
\end{eqnarray}
similar to \eqref{interim update1}, where $E_{\tilde P_X}[g(\mathbf X)]-E_{\bar P_X}[g(\mathbf X)]=O_p(1/\sqrt m)$ similar to \eqref{interim updated3}, and $E_{\bar P_X}[g(\mathbf X)]-E_{P_X^0}[g(\mathbf X)]=O_p(1/\sqrt m)$ by using the standard central limit theorem for $U$-statistic (Theorem A in Section 5.5 in \cite{serfling2009approximation}) and, with the assumption $E_{P_X^0}[g(X_{i_1},\ldots,X_{i_T})^2]<\infty$ for any $1\leq i_1,\ldots,i_T\leq T$, translating it to $V$-statistic (the lemma in Section 5.7.3 in \cite{serfling2009approximation})). Therefore, we have $E_{\tilde P_X}[g(\mathbf X)]-E_{P_X^0}[g(\mathbf X)]=O_p(1/\sqrt m)$. 


In conclusion, we have found a solution $\tilde P_X$ that is feasible for \eqref{KS RO max discretized direct} and \eqref{KS RO min discretized direct} with probability asymptotically $1-\alpha$ as $n\to\infty$ and $m/n\to\infty$. Moreover, $\psi(\tilde P_X)-\psi(P_X^0)=O_p(1/\sqrt m)$. Therefore, we have
\begin{eqnarray*}
1-\alpha&\leq&\lim_{n\to\infty,m/n\to\infty}\mathbb P\left(\|E_{\tilde P_X}[I(h(\mathbf X)\leq \cdot)]-\hat F_Y(\cdot)\|_\infty\leq\frac{q_{1-\alpha}}{\sqrt n}\right)\\
&\leq&\liminf_{n\to\infty,m/n\to\infty}\mathbb P\left(\hat{\underline Z}\leq\psi(\tilde P_X)\leq\hat{\overline Z}\right)\\
&=&\liminf_{n\to\infty,m/n\to\infty}\mathbb P\left(\hat{\underline Z}+O_p\left(\frac{1}{\sqrt m}\right)\leq\psi(P_X^0)\leq\hat{\overline Z}+O_p\left(\frac{1}{\sqrt m}\right)\right)
\end{eqnarray*}
which concludes the theorem.\Halmos
\endproof

We provide some remark on the case where we consider $h(\mathbf X,\mathbf W)$ and $g(\mathbf X,\mathbf W)$ for some collection of auxiliary input variate sequences $\mathbf W$ that is independent of $\mathbf X$ and has a known distribution. In this case, the results in Sections \ref{sec:formulation} and \ref{sec:KS} all hold with the $E_{P_X}[\cdot]$ interpreted as the joint expectation taken with respect to both the product measure of $P_X$ and $P_W^0$, the known distribution of $\mathbf W$. In the proofs above, we keep the expectation $E_{P_X}[\cdot]$ as taken under the product measure of $P_X$ only, but we use a conditioning argument, namely we change $I(h(\mathbf X)\leq y)$ to $P_{P_W^0}(h(\mathbf X,\mathbf W)\leq y|\mathbf X)=E_{P_W^0}[I(h(\mathbf X)\leq y)|\mathbf X]$ and $g(\mathbf X)$ to $E_{P_W^0}[g(\mathbf X,\mathbf W)|\mathbf X]$, where $P_{P_W^0}(\cdot|\mathbf X)$ and $E_{P_W^0}[\cdot|\mathbf X]$ denote the conditional probability and expectation under the true distribution of $\mathbf W$ given $\mathbf X$. In particular, in the proof of Theorem \ref{main guarantee}, we have that $P_{P_W^0}(h(\mathbf X,\mathbf W)\leq y|\mathbf X=\mathbf x)$ is non-decreasing given any $\mathbf x$, and $P_{P_W^0}(h(\mathbf X,\mathbf W)\leq y|\mathbf X=\mathbf x)\leq1$, which, via Problem 3 in Chapter 2.7 of \cite{van1996weak} again, gives a polynomial bracketing number for the class of functions $\{P_{P_W^0}(h(\mathbf X,\mathbf W)\leq y|\mathbf X=\mathbf x)\prod_{t=1}^TL(x_t):y\in\mathbb R\}$. We also have $E_{P_X^0}[E_{P_W^0}[g(X_{i_1},\ldots,X_{i_T},\mathbf W)|X_{i_1},\ldots,X_{i_T}]^2]\leq E_{P_X^0,P_W^0}[g(X_{i_1},\ldots,X_{i_T},\mathbf W)^2]<\infty$ for any $1\leq i_1,\ldots,i_T\leq T$, where $E_{P_X^0,P_W^0}[\cdot]$ denotes the joint expectation under the product measure of $P_X^0$ and $P_W^0$, so that the central limit theorem for ensuring the approximation of the objective value holds in the proof. Other proofs follow quite trivially.


\section{Proofs for Section \ref{sec:procedure}}\label{sec:procedure proofs}
\proof{Proof of Proposition \ref{prop:quadratic penalty}}
Consider the equivalent reformulation of the program \eqref{KS RO min discretized p}
\begin{equation}\label{equality constraints}
\begin{array}{ll}
\min&\psi(\mathbf p)\\
\text{subject to}&E_{\mathbf p}[I(h(\mathbf X)\leq y_j)]-s_j=0,j=1,\ldots,n\\
&\hat F_Y(y_j+)-\frac{q_{1-\alpha}}{\sqrt n}\leq s_j\leq\hat F_Y(y_j-)+\frac{q_{1-\alpha}}{\sqrt n},j=1,\ldots,n\\
&\mathbf p\in\mathcal P
\end{array}
\end{equation}
where both $\mathbf p$ and $\mathbf s$ are viewed as decision variables. An application of the conventional quadratic penalty method (\cite{bertsekas1999nonlinear}) for equality constraints yields the following optimization sequence
\begin{equation}\label{quad penalty min}
\begin{array}{ll}
\min&\psi(\mathbf p)+c\sum_{j=1}^n(E_{\mathbf p}[I(h(\mathbf X)\leq y_j)]-s_j)^2\\
\text{subject to}&\hat F_Y(y_j+)-\frac{q_{1-\alpha}}{\sqrt n}\leq s_j\leq\hat F_Y(y_j-)+\frac{q_{1-\alpha}}{\sqrt n},j=1,\ldots,n\\
&\mathbf p\in\mathcal P
\end{array}
\end{equation}
for $c>0$, which is equivalent to \eqref{KS RO min discretized penalty} with $\lambda=1/c$. Proposition 4.2.1 in \cite{bertsekas1999nonlinear} entails that as $c\to\infty$ ($\lambda\to 0$), every limit point $(\mathbf p^*,\mathbf s^*)$ of the sequence of optimal solutions $\{(\mathbf p^*(\lambda),\mathbf s^*(\lambda))\}$ to \eqref{quad penalty min} is an optimal solution to \eqref{equality constraints}, given that \eqref{equality constraints} is feasible. Note that due to optimality, the optimal slack variables $\mathbf s^*(\lambda)=(s^*_{1}(\lambda),\ldots,s^*_{n}(\lambda))$ must take the following form
\begin{equation*}
s^*_{j}(\lambda)=\Pi_j(E_{\mathbf p^*(\lambda)}[I(h(\mathbf X)\leq y_j)]),\ j=1,\ldots,n
\end{equation*}
where each $\Pi_j$ is the projection defined in \eqref{pie_j}. Since projections are continuous maps, the operations of taking limit points and coordinate projection are interchangeable, i.e.
\begin{align*}
&\{\mathbf p^*:\mathbf p^*\text{ is a limit point of }\{\mathbf p^*(\lambda)\}\}\\
=&\{\mathbf p^*:\text{there exists an }\mathbf s^*\text{ s.t. }(\mathbf p^*,\mathbf s^*)\text{ is a limit point of }\{(\mathbf p^*(\lambda),\mathbf s^*(\lambda))\}\}.
\end{align*}
This allows translation of optimality of the limit point of $\{(\mathbf p^*(\lambda),\mathbf s^*(\lambda))\}$ to optimality of the limit point of $\{\mathbf p^*(\lambda)\}$. The desired conclusion follows.\Halmos
\endproof

\proof{Proof of Proposition \ref{lemma:derivative1}.}
Part 1 and the expression for $\Psi_i$ in part 2 come from a direct application of \cite{gl15_1} and \cite{ghosh2015mirror}. We will prove \eqref{score function2} and \eqref{score function3} in part 2 only, but in the more general setting of differentiable functions of expectations. Let $f(\mathbf X)$ with $\mathbf X=(X_1,\ldots,X_{T_f})$ be a performance function, where $T_f$ is a finite and deterministic time horizon, and $\Phi(y):\R \to \R$ be any differentiable function. By the chain rule
\begin{equation*}
\Phi_i(\mathbf p):=\frac{d}{d\epsilon}\Phi(\mathbf E_{(1-\epsilon)\mathbf p+\epsilon \mathbf 1_i}[f(\mathbf X)])\Big\rvert_{\epsilon=0^+}=\frac{d}{dy}\Phi(\mathbf E_{\mathbf p}[f(\mathbf X)])\frac{d}{d\epsilon}\mathbf E_{(1-\epsilon)\mathbf p+\epsilon \mathbf 1_i}[f(\mathbf X)]\rvert_{\epsilon=0^+}.
\end{equation*}
Similar to \eqref{score function1} we have
\begin{equation*}
\frac{d}{d\epsilon}\mathbf E_{(1-\epsilon)\mathbf p+\epsilon \mathbf 1_i}[f(\mathbf X)]\rvert_{\epsilon=0^+}=E_{\mathbf p}[f(\mathbf X)S_i(\mathbf X;\mathbf p)]
\end{equation*}
where
\begin{equation*}
S_i(\mathbf x;\mathbf p)=\sum_{t=1}^{T_f}\frac{I_i(x_t)}{p_i}-T_f.
\end{equation*}
Therefore the following expression holds for the derivative
\begin{equation*}
\Phi_i(\mathbf p)=\frac{d}{dy}\Phi(\mathbf E_{\mathbf p}[f(\mathbf X)])E_{\mathbf p}[f(\mathbf X)S_i(\mathbf X;\mathbf p)].
\end{equation*}
\eqref{score function2} and \eqref{score function3} follow from applying the above result to $f(\mathbf X)=h(\mathbf X)$, $\Phi(y)=(y-\Pi_j(y))^2$ and $\Phi(y)=(y-s_j)^2$ respectively, together with the linearity of differentiation. Note that $\frac{d}{dy}(y-\Pi_j(y))^2=2(y-\Pi_j(y))$.\Halmos
\endproof

\proof{Proof of Proposition \ref{sol:constrained stepwise}.}First note that the function $\mu(\eta)$ is continuous and strictly increasing in the interval $[0,\max_ip_ie^{-\xi_i}]$, and satisfies $\mu(0)=0,\mu(\max_ip_ie^{-\xi_i})=1/m$ at the endpoints. So indeed there exists a unique $\eta^*$ that solves \eqref{threshold:eta}. Then we show \eqref{sol2} is indeed the optimal solution. Consider the Lagrangian
\begin{equation*}
L(\mathbf q,\lambda,\bm\beta)=\bm\xi'(\mathbf q-\mathbf p)+V(\mathbf p,\mathbf q)+\lambda(\sum_{i=1}^mq_i-1)-\sum_{i=1}^m\beta_i(q_i-\epsilon)
\end{equation*}
defined for $\beta_i\geq 0$ and $\lambda\in \R$. Since \eqref{generic_epsilon} is a convex program with linear constraints and obviously Slater's condition holds, by Proposition 6.2.5 and Proposition 6.4.4 of \cite{bertsekas2003convex} it suffices to find dual variables $\lambda^*$ and $\beta_i^*$ such that the solution given by \eqref{sol2} satisfies the set of KKT conditions
\begin{align}
&\frac{\partial L}{\partial q_i}=\xi_i+\log \frac{q^*_i}{p_i}+1+\lambda^*-\beta^*_i=0\text{ for }i=1,\ldots,m\label{KKT1}\\
&\sum_{i=1}^mq^*_i=1,\ q^*_i\geq \epsilon,\text{ for }i=1,\ldots,m\label{KKT2}\\
&\beta^*_i\geq 0,\ \beta^*_i(q^*_i-\epsilon)=0\text{ for }i=1,\ldots,m.\label{KKT3}
\end{align}
Equations \eqref{KKT2} obviously hold because of equation \eqref{threshold:eta}. Equations \eqref{KKT1} can be rewritten as
\begin{align*}
q^*_i=p_ie^{-\xi_i-1-\lambda^*+\beta^*_i}\text{ for }i=1,\ldots,m
\end{align*}
which hold if $\lambda^*,\beta_i^*$ are chosen such that
\begin{equation*}
e^{1+\lambda^*}=\sum_{i=1}^m\max\{\eta^*,p_ie^{-\xi_i}\},\ e^{\beta_i^*}=\frac{\max\{\eta^*,p_ie^{-\xi_i}\}}{p_ie^{-\xi_i}}.
\end{equation*}
It is obvious that such chosen $\beta_i^*\geq 0$. To show complementary slackness \eqref{KKT3}, note that if $q^*_i>\epsilon$ then \eqref{threshold:eta} forces $p_ie^{-\xi_i}>\eta^*$ which results in $\beta^*_i=0$.\Halmos
\endproof
\proof{Proof of Theorem \ref{main guarantee epsilon}.}Consider the auxiliary programs obtained from replacing $\alpha$ by some $\alpha'>\alpha$ in \eqref{KS RO min discretized p}
\begin{equation*}
\begin{array}{ll}
\max&\psi(\mathbf p)\\
\text{subject to}&\hat F_Y(y_j+)-\frac{q_{1-\alpha'}}{\sqrt n}\leq E_{\mathbf p}[I(h(\mathbf X)\leq y_j)]\leq\hat F_Y(y_j-)+\frac{q_{1-\alpha'}}{\sqrt n},j=1,\ldots,n\\
&\mathbf p\in\mathcal P
\end{array}
\end{equation*}
and
\begin{equation*}
\begin{array}{ll}
\min&\psi(\mathbf p)\\
\text{subject to}&\hat F_Y(y_j+)-\frac{q_{1-\alpha'}}{\sqrt n}\leq E_{\mathbf p}[I(h(\mathbf X)\leq y_j)]\leq\hat F_Y(y_j-)+\frac{q_{1-\alpha'}}{\sqrt n},j=1,\ldots,n\\
&\mathbf p\in\mathcal P.
\end{array}
\end{equation*}
Denote by $\mathbf p^{*'}_{\max}$ and $\mathbf p^{*'}_{\min}$ optimal solutions of the above maximization and minimization programs, which by Theorem \ref{main guarantee} satisfy
$$\liminf_{n\to\infty,m/n\to\infty}\mathbb P\left(\psi(\mathbf p^{*'}_{\min})+O_p\left(\frac{1}{\sqrt m}\right)\leq\psi(P_X^0)\leq\psi(\mathbf p^{*'}_{\max})+O_p\left(\frac{1}{\sqrt m}\right)\right)\geq1-\alpha'.$$

Now, we try to show that $\hat{\underline Z}_{\epsilon}\leq \psi(\mathbf p^{*'}_{\min})+O(m\epsilon)$ and $\hat{\overline Z}_{\epsilon}\geq \psi(\mathbf p^{*'}_{\max})-O(m\epsilon)$, therefore to conclude that
\begin{equation}\label{coverage:alpha'}
\liminf_{n\to\infty,m/n\to\infty}\mathbb P\left(\hat{\underline Z}_{\epsilon}+O_p\left(m\epsilon+\frac{1}{\sqrt m}\right)\leq\psi(P_X^0)\leq\hat{\overline Z}_{\epsilon}+O_p\left(m\epsilon+\frac{1}{\sqrt m}\right)\right)\geq1-\alpha'.
\end{equation}
To avoid repetition, we only prove the minimization case here. To proceed, let $\mathbf p,\mathbf q\in \mathcal P$ be two arbitrary probability distributions in $\mathcal P$, and $\mathbf p^S,\mathbf q^S$ be the corresponding $S$-fold product measure, then we have
\begin{equation*}
\lvert\psi(\mathbf p)-\psi(\mathbf q)\rvert=\lvert E_{\mathbf p}[g(\mathbf X)]-E_{\mathbf q}[g(\mathbf X)]\rvert\leq 2\sup_{\mathbf X}\lvert g(\mathbf X)\rvert\cdot\Vert\mathbf p^S-\mathbf q^S\Vert_{TV}
\end{equation*}
where $\Vert\cdot\Vert_{TV}$ denotes the total variation distance between the product measures. It is well-known that the total variation distance between product measures can be bounded as (see, e.g.~Lemma 3.6.2 of \cite{durrett2010probability})
\begin{equation*}
\Vert\mathbf p^S-\mathbf q^S\Vert_{TV}\leq S\Vert\mathbf p-\mathbf q\Vert_{TV},
\end{equation*}
therefore
\begin{equation*}
\lvert\psi(\mathbf p)-\psi(\mathbf q)\rvert\leq 2S\sup_{\mathbf X}\lvert g(\mathbf X)\rvert\cdot\Vert\mathbf p-\mathbf q\Vert_{TV}=C_1\Vert\mathbf p-\mathbf q\Vert_{TV}.
\end{equation*}
Similarly for the constraint functions we have
\begin{equation*}
\lvert E_{\mathbf p}[I(h(\mathbf X)\leq y_j)]-E_{\mathbf q}[I(h(\mathbf X)\leq y_j)]\rvert\leq 2T\Vert\mathbf p-\mathbf q\Vert_{TV}=C_2\Vert\mathbf p-\mathbf q\Vert_{TV},j=1,\ldots,n.
\end{equation*}
Consider the total variation ball of radius $m\epsilon$ surrounding $\mathbf p^{*'}_{\min}$
\begin{equation*}
B_{TV}(\mathbf p^{*'}_{\min},m\epsilon)=\{\mathbf p\in\mathcal P:\Vert \mathbf p^{*'}_{\min}-\mathbf p\Vert_{TV}\leq m\epsilon\}.
\end{equation*}
It is clear that for all $\mathbf p\in B_{TV}(\mathbf p^{*'}_{\min},m\epsilon)$ it holds
\begin{align}
&\lvert\psi(\mathbf p)-\psi(\mathbf p^{*'}_{\min})\rvert\leq C_1m\epsilon\label{value gap}\\
&\lvert E_{\mathbf p}[I(h(\mathbf X)\leq y_j)]-E_{\mathbf p^{*'}_{\min}}[I(h(\mathbf X)\leq y_j)]\rvert \leq C_2m\epsilon,j=1,\ldots,n.\label{constraint gap}
\end{align}
Note that $\mathbf p^{*'}_{\min}$ is optimal and hence feasible for the program with $\alpha'$, thus the inequality \eqref{constraint gap} ensures for all $\mathbf p\in B_{TV}(\mathbf p^{*'}_{\min},m\epsilon)$
\begin{equation*}
\hat F_Y(y_j+)-\frac{q_{1-\alpha'}}{\sqrt n}-C_2m\epsilon\leq E_{\mathbf p}[I(h(\mathbf X)\leq y_j)]\leq \hat F_Y(y_j+)+\frac{q_{1-\alpha'}}{\sqrt n}+C_2m\epsilon,j=1,\ldots,n.
\end{equation*}
Since $\epsilon=o(1/(m\sqrt n))$, for large enough $m,n$ we have $C_2m\epsilon\leq (q_{1-\alpha}-q_{1-\alpha'})/\sqrt n$ which results in
\begin{equation*}
\hat F_Y(y_j+)-\frac{q_{1-\alpha}}{\sqrt n}\leq E_{\mathbf p}[I(h(\mathbf X)\leq y_j)]\leq \hat F_Y(y_j+)+\frac{q_{1-\alpha}}{\sqrt n},j=1,\ldots,n.
\end{equation*}
That is, all $\mathbf p\in B_{TV}(\mathbf p^{*'}_{\min},m\epsilon)$ satisfy the first constraint in \eqref{KS RO min discretized p epsilon}. In view of inequality \eqref{value gap}, it remains to show that $B_{TV}(\mathbf p^{*'}_{\min},m\epsilon)\cap \mathcal P(\epsilon)\neq \emptyset$ in order to conclude $\hat{\underline Z}_{\epsilon}\leq \psi(\mathbf p^{*'}_{\min})+O(m\epsilon)$. Easily one can verify that for any $\mathbf p\in\mathcal P$ it holds $\inf\{\Vert\mathbf p-\mathbf q\Vert_{TV}:\mathbf q\in \mathcal P(\epsilon)\}\leq (m-1)\epsilon$, and in particular $\inf\{\Vert\mathbf p^{*'}_{\min}-\mathbf q\Vert_{TV}:\mathbf q\in \mathcal P(\epsilon)\}\leq (m-1)\epsilon$ which implies $B_{TV}(\mathbf p^{*'}_{\min},m\epsilon)\cap \mathcal P(\epsilon)\neq \emptyset$.

Lastly note that \eqref{coverage:alpha'} holds true for arbitrary $\alpha'>\alpha$, hence holds for $\alpha$ as well. This concludes the theorem.\Halmos
\endproof

\begin{lemma}
For any $i,j$ and $l=1,2$, the moments of gradient estimators
\begin{equation*}
E_{\mathbf p}\left[( g(\mathbf X)S_i(\mathbf X;\mathbf p))^l\right],\;E_{\mathbf p}\big[\big(I(h(\mathbf X)\leq y_j)S_i(\mathbf X;\mathbf p)\big)^l\big]
\end{equation*}
are continuous in $\mathcal P^o=\{\mathbf p\in \mathcal P:p_i>0\text{ for all }i\}$, the relative interior of $\mathcal P$.
\label{lemma:derivative continuity}
\end{lemma}
\proof{Proof of Lemma \ref{lemma:derivative continuity}.}Restricted to $\mathcal P^o$, each of the moments can be written as the sum of finitely many terms each of which are smooth in $\mathbf p$. A sum of finitely many smooth functions is also smooth, hence continuous.\Halmos
\endproof

\begin{lemma}\label{rate:slack}
Let $\{D^k\}_{k=1}^{\infty}$ be a positive sequence. If for $0<\alpha_2<\alpha_1\leq 1$ and constants $C_1,C_2>0$ it holds $D^{k+1}\leq (1-\frac{C_1}{k^{\alpha_2}})D^k+C_2(\frac{1}{k^{2\alpha_2}}+\frac{1}{k^{2\alpha_1-\alpha_2}})$ for all $k$ large enough, then there exits a constant $C>0$ such that $D^k\leq C(\frac{1}{k^{\alpha_2}}+\frac{1}{k^{2(\alpha_1-\alpha_2)}})$ for all $k$.
\end{lemma}
\proof{Proof of Lemma \ref{rate:slack}.}
Assume $D^k\leq C(\frac{1}{k^{\alpha_2}}+\frac{1}{k^{2(\alpha_1-\alpha_2)}})$, then
\begin{align*}
D^{k+1}&\leq (1-\frac{C_1}{k^{\alpha_2}})D^k+C_2(\frac{1}{k^{2\alpha_2}}+\frac{1}{k^{2\alpha_1-\alpha_2}})\\
&\leq \frac{C}{k^{\alpha_2}}+\frac{C}{k^{2(\alpha_1-\alpha_2)}}-\frac{C_1C-C_2}{k^{2\alpha_2}}-\frac{C_1C-C_2}{k^{2\alpha_1-\alpha_2}}\\
&\leq \frac{C}{(k+1)^{\alpha_2}}+\frac{C\alpha_2}{k^{\alpha_2+1}}+\frac{C}{(k+1)^{2(\alpha_1-\alpha_2)}}+\frac{C\cdot 2(\alpha_1-\alpha_2)}{k^{2(\alpha_1-\alpha_2)+1}}-\frac{C_1C-C_2}{k^{2\alpha_2}}-\frac{C_1C-C_2}{k^{2\alpha_1-\alpha_2}}\\
&\leq \frac{C}{(k+1)^{\alpha_2}}+\frac{C\alpha_2k^{\alpha_2-1}}{k^{2\alpha_2}}+\frac{C}{(k+1)^{2(\alpha_1-\alpha_2)}}+\frac{C\cdot 2(\alpha_1-\alpha_2)k^{\alpha_2-1}}{k^{2\alpha_1-\alpha_2}}-\frac{C_1C-C_2}{k^{2\alpha_2}}-\frac{C_1C-C_2}{k^{2\alpha_1-\alpha_2}}\\
&\leq \frac{C}{(k+1)^{\alpha_2}}+\frac{C}{(k+1)^{2(\alpha_1-\alpha_2)}}-\frac{C(C_1-\alpha_2k^{\alpha_2-1})-C_2}{k^{2\alpha_2}}-\frac{C(C_1-2(\alpha_1-\alpha_2)k^{\alpha_2-1})-C_2}{k^{2\alpha_1-\alpha_2}}\\
&\leq \frac{C}{(k+1)^{\alpha_2}}+\frac{C}{(k+1)^{2(\alpha_1-\alpha_2)}}.
\end{align*}
Note that the above argument goes through when $k$ is large and $C$ is chosen such that $\frac{C_1}{k^{\alpha_2}}<1$, $C(C_1-2(\alpha_1-\alpha_2)k^{\alpha_2-1})-C_2\geq 0$ and $C(C_1-\alpha_2k^{\alpha_2-1})-C_2\geq 0$. By induction $D^k\leq C(\frac{1}{k^{\alpha_2}}+\frac{1}{k^{2(\alpha_1-\alpha_2)}})$ holds for all sufficiently large $k$. By enlarging $C$ one can make it hold for all $k$.\Halmos\endproof

\proof{Proof of Theorem \ref{thm:algo2}.}
We borrow from Lemma 2.1 in \cite{nemirovski2009robust} the inequality
\begin{equation}\label{ineq:p}
V(\mathbf p^{k+1},\mathbf p_{\epsilon}^*(\lambda^k))-V(\mathbf p^k,\mathbf p_{\epsilon}^*(\lambda^k))\leq\gamma^k(\lambda^k\hat{\bm\Psi}^k+\hat{\bm\phi}_{\mathbf p}^k)'(\mathbf p_{\epsilon}^*(\lambda^k)-\mathbf p^k)+\frac{(\gamma^k)^2\|\lambda^k\hat{\bm\Psi}^k+\hat{\bm\phi}_{\mathbf p}^k\|_\infty^2}{2}
\end{equation}
which holds as long as $\mathbf p^{k+1}$ is the prox-mapping of $\mathbf p^k$. The norm $\|\cdot\|_\infty$ is the supremum norm, the dual of the $L_1$-norm that is used in the strong convexity property of $\omega(\mathbf p)=\sum_{i=1}^mp_i\log p_i$, with $\alpha=1$. Note that $V(\mathbf p^{k+1},\mathbf p_{\epsilon}^*(\lambda^k))=\sum_{i=1}^mp^*_i(\lambda^k)(\log p^*_i(\lambda^k)-\log p^{k+1}_i)$ and both $\mathbf p_{\epsilon}^*(\lambda^k),\mathbf p^{k+1}\in\mathcal P(\epsilon)$, by mean value theorem it holds
\begin{align}\label{ktok+1}
V(\mathbf p^{k+1},\mathbf p_{\epsilon}^*(\lambda^{k+1}))-V(\mathbf p^{k+1},\mathbf p_{\epsilon}^*(\lambda^k))\leq C\lvert\log \epsilon\rvert \Vert\mathbf p_{\epsilon}^*(\lambda^{k+1})-\mathbf p_{\epsilon}^*(\lambda^k)\Vert
\end{align}
where $C$ is an absolute constant. This gives
\begin{align}
\nonumber&V(\mathbf p^{k+1},\mathbf p_{\epsilon}^*(\lambda^{k+1}))-V(\mathbf p^k,\mathbf p_{\epsilon}^*(\lambda^k))\\
\leq&\gamma^k(\lambda^k\hat{\bm\Psi}^k+\hat{\bm\phi}_{\mathbf p}^k)'(\mathbf p_{\epsilon}^*(\lambda^k)-\mathbf p^k)+\frac{(\gamma^k)^2\|\lambda^k\hat{\bm\Psi}^k+\hat{\bm\phi}_{\mathbf p}^k\|_\infty^2}{2}+C\lvert\log \epsilon\rvert \Vert\mathbf p_{\epsilon}^*(\lambda^{k+1})-\mathbf p_{\epsilon}^*(\lambda^k)\Vert\label{prox-mapping relation}
\end{align}

Let $\mathcal F^k$ be the filtration generated by $\{\mathbf p^1,\mathbf s^1,\ldots,\mathbf p^k,\mathbf s^k\}$. Taking conditional expectation of \eqref{prox-mapping relation} with respect to $\mathcal F^k$, we have
\begin{align}
\nonumber&E[V(\mathbf p^{k+1},\mathbf p_{\epsilon}^*(\lambda^{k+1}))-V(\mathbf p^k,\mathbf p_{\epsilon}^*(\lambda^k))|\mathcal F^k]\\
\nonumber\leq& \gamma^k(\lambda^k\bm\Psi(\mathbf p^k)+\bm\phi(\mathbf p^k))'(\mathbf p_{\epsilon}^*(\lambda^k)-\mathbf p^k)+\gamma^k(E[\hat{\bm\phi}_{\mathbf p}^k\vert \mathcal F^k]-\bm\phi(\mathbf p^k))'(\mathbf p_{\epsilon}^*(\lambda^k)-\mathbf p^k)\\
&+\frac{1}{2}(\gamma^k)^2E[\|\lambda^k\hat{\bm\Psi}^k+\hat{\bm\phi}_{\mathbf p}^k\|_\infty^2|\mathcal F^k]+C\lvert\log \epsilon\rvert \Vert\mathbf p_{\epsilon}^*(\lambda^{k+1})-\mathbf p_{\epsilon}^*(\lambda^k)\Vert.\label{interim2}
\end{align}
Note that on the right hand side we are still using $\bm\phi(\mathbf p^k)$, the derivative of the quadratic penalty in the formulation \eqref{KS RO min discretized penalty2}, rather than $\bm\phi_{\mathbf p}(\mathbf p^k,\mathbf s^k)$.

In order to use the martingale convergence theorem, we examine the following
\begin{align}
&\sum_{k=1}^\infty E[E[V(\mathbf p^{k+1},\mathbf p_{\epsilon}^*(\lambda^{k+1}))-V(\mathbf p^k,\mathbf p_{\epsilon}^*(\lambda^k))|\mathcal F^k]^+]\label{ineq:convergence}\\
\nonumber\leq& \sum_{k=1}^\infty O(\gamma^k\sqrt{E[\Vert E[\hat{\bm\phi}_{\mathbf p}^k\vert \mathcal F^k]-\bm\phi(\mathbf p^k) \Vert^2]})+\sum_{k=1}^\infty\frac{1}{2}(\gamma^k)^2E[\|\lambda^k\hat{\bm\Psi}^k+\hat{\bm\phi}_{\mathbf p}^k\|_\infty^2|\mathcal F^k]+\sum_{k=1}^\infty C\lvert\log \epsilon\rvert \Vert\mathbf p_{\epsilon}^*(\lambda^{k+1})-\mathbf p_{\epsilon}^*(\lambda^k)\Vert.
\end{align}
We need to bound two quantities, $E[\Vert E[\hat{\bm\phi}_{\mathbf p}^k\vert \mathcal F^k]-\bm\phi(\mathbf p^k)\Vert^2]$ and $E[\|\lambda^k\hat{\bm\Psi}^k+\hat{\bm\phi}_{\mathbf p}^k\|_\infty^2|\mathcal F^k]$. To bound the first one
\begin{align}
\nonumber&E[\Vert E[\hat{\bm\phi}_{\mathbf p}^k\vert \mathcal F^k]-\bm\phi(\mathbf p^k))\Vert^2]\\
\nonumber=&\sum_{i=1}^mE[\lvert E[\hat{\phi}^k_{\mathbf p,i}\vert \mathcal F^k]-\phi_i(\mathbf p^k)\rvert^2]\\
\nonumber=&4\sum_{i=1}^mE\big[\big\lvert\sum_{j=1}^n(\Pi_j(E_{\mathbf p^k}[I(h(\mathbf X)\leq y_j)])-s^k_j)E_{\mathbf p^k}[I(h(\mathbf X)\leq y_i)S_i(\mathbf X;\mathbf p^k)]\big\rvert^2\big]\\
\nonumber\leq &4\sum_{i=1}^mE\big[\sum_{j=1}^n(\Pi_j(E_{\mathbf p^k}[I(h(\mathbf X)\leq y_j)])-s^k_j)^2\sum_{j=1}^n(E_{\mathbf p^k}[I(h(\mathbf X)\leq y_i)S_i(\mathbf X;\mathbf p^k)])^2\big]\\
\nonumber\leq &4E\big[\sum_{j=1}^n(\Pi_j(E_{\mathbf p^k}[I(h(\mathbf X)\leq y_j)])-s^k_j)^2\big]\sum_{i=1}^m\sum_{j=1}^n\sup_{\mathbf p\in \mathcal P(\epsilon)}(E_{\mathbf p}[I(h(\mathbf X)\leq y_i)S_i(\mathbf X;\mathbf p)])^2\\
\leq &Cmn\sum_{j=1}^nE[(\Pi_j(E_{\mathbf p^k}[I(h(\mathbf X)\leq y_j)])-s^k_j)^2]\label{first term bd}
\end{align}
where in the first inequality we use Cauchy Schwartz inequality, and the third inequality holds because each $E_{\mathbf p}[I(h(\mathbf X)\leq y_i)S_i(\mathbf X;\mathbf p)]$ by Lemma \ref{lemma:derivative continuity} is continuous in $\mathbf p$ and hence by a compactness argument is uniformly bounded in $\mathcal P(\epsilon)$. Therefore the key step lies in deriving an upper bound for each $E[(\Pi_j(E_{\mathbf p^k}[I(h(\mathbf X)\leq y_j)])-s^k_j)^2]$, for which we need the counterpart of \eqref{ineq:p} for $s^k_j$, i.e.
\begin{align*}
&\frac{1}{2}( s^{k+1}_j-\Pi_j(E_{\mathbf p^k}[I(h(\mathbf X)\leq y_j)]))^2-\frac{1}{2}(s^{k}_j-\Pi_j(E_{\mathbf p^k}[I(h(\mathbf X)\leq y_j)]))^2\\
\leq &\beta^k\hat\phi_{\mathbf s,j}^k(\Pi_j(E_{\mathbf p^k}[I(h(\mathbf X)\leq y_j)]-s^k_j)+\frac{1}{2}(\beta^k)^2(\hat\phi_{\mathbf s,j}^k)^2.
\end{align*}
Taking expectation with respect to $\mathcal F_k$ gives
\begin{align}
\nonumber&\frac{1}{2}E[( s^{k+1}_j-\Pi_j(E_{\mathbf p^k}[I(h(\mathbf X)\leq y_j)]))^2\vert \mathcal F_k]-\frac{1}{2}(s^{k}_j-\Pi_j(E_{\mathbf p^k}[I(h(\mathbf X)\leq y_j)]))^2\\
\nonumber\leq &-2\beta^k(E_{\mathbf p^k}[I(h(\mathbf X)\leq y_j)]-s^k_j)(\Pi_j(E_{\mathbf p^k}[I(h(\mathbf X)\leq y_j)])-s^k_j)+\frac{1}{2}(\beta^k)^2(2+q_{1-\alpha}/\sqrt n)^2\\
\leq &-2\beta^k(\Pi_j(E_{\mathbf p^k}[I(h(\mathbf X)\leq y_j)])-s^k_j)^2+\frac{C}{2}(\beta^k)^2.\label{ineq2:s}
\end{align}
Note that with step size $\gamma^k$ we have
\begin{align*}
&E[( s^{k+1}_j-\Pi_j(E_{\mathbf p^k}[I(h(\mathbf X)\leq y_j)]))^2\vert \mathcal F_k]\\
=&E[( s^{k+1}_j-\Pi_j(E_{\mathbf p^{k+1}}[I(h(\mathbf X)\leq y_j)]))^2\vert \mathcal F_k]\\
&+2E[( s^{k+1}_j-\Pi_j(E_{\mathbf p^{k+1}}[I(h(\mathbf X)\leq y_j)]))( \Pi_j(E_{\mathbf p^{k+1}}[I(h(\mathbf X)\leq y_j)])-\Pi_j(E_{\mathbf p^{k}}[I(h(\mathbf X)\leq y_j)]))\vert \mathcal F_k]\\
&+E[( \Pi_j(E_{\mathbf p^{k+1}}[I(h(\mathbf X)\leq y_j)])-\Pi_j(E_{\mathbf p^{k}}[I(h(\mathbf X)\leq y_j)]))^2\vert \mathcal F_k]\\
\geq &E[( s^{k+1}_j-\Pi_j(E_{\mathbf p^{k+1}}[I(h(\mathbf X)\leq y_j)]))^2\vert \mathcal F_k]\\
&-2\sqrt{E[( s^{k+1}_j-\Pi_j(E_{\mathbf p^{k+1}}[I(h(\mathbf X)\leq y_j)]))^2\vert \mathcal F_k]}\sqrt{E[(E_{\mathbf p^{k+1}}[I(h(\mathbf X)\leq y_j)]-E_{\mathbf p^{k}}[I(h(\mathbf X)\leq y_j)])^2\vert \mathcal F_k]}\\
\geq &E[( s^{k+1}_j-\Pi_j(E_{\mathbf p^{k+1}}[I(h(\mathbf X)\leq y_j)]))^2\vert \mathcal F_k]-2(\sqrt{E[( s^{k+1}_j-\Pi_j(E_{\mathbf p^{k+1}}[I(h(\mathbf X)\leq y_j)]))^2\vert \mathcal F_k]}C\gamma^k)\\
\geq &E[( s^{k+1}_j-\Pi_j(E_{\mathbf p^{k+1}}[I(h(\mathbf X)\leq y_j)]))^2\vert \mathcal F_k]-2\beta^kE[( s^{k+1}_j-\Pi_j(E_{\mathbf p^{k+1}}[I(h(\mathbf X)\leq y_j)]))^2\vert \mathcal F_k]-\frac{C^2(\gamma^k)^2}{2\beta^k}
\end{align*}
where the second last inequality follows from
\begin{align*}
\lvert E_{\mathbf p^{k+1}}[I(h(\mathbf X)\leq y_j)]-E_{\mathbf p^{k}}[I(h(\mathbf X)\leq y_j)]\rvert&\leq \Vert \mathbf p^{k+1}-\mathbf p^k\Vert\cdot\sup_{\mathbf p\in\mathcal P(\epsilon)}\Vert\nabla E_{\mathbf p}[I(h(\mathbf X)\leq y_j)]\Vert\\
&\leq C \Vert \mathbf p^{k+1}-\mathbf p^k\Vert=O(\gamma^k)
\end{align*}
and in the last inequality we use Young's inequality. Substituting the above into \eqref{ineq2:s} gives
\begin{align*}
&E[( s^{k+1}_j-\Pi_j(E_{\mathbf p^{k+1}}[I(h(\mathbf X)\leq y_j)]))^2\vert \mathcal F_k]\\
\leq &\frac{1-4\beta^k}{1-2\beta^k}( s^{k}_j-\Pi_j(E_{\mathbf p^{k}}[I(h(\mathbf X)\leq y_j)]))^2+C((\beta^k)^2+\frac{(\gamma^k)^2}{\beta^k}).
\end{align*}
Hence taking full expectation we have the following recursion
\begin{align*}
E[( s^{k+1}_j-\Pi_j(E_{\mathbf p^{k+1}}[I(h(\mathbf X)\leq y_j)]))^2]\leq (1-2\beta^k)E[( s^{k}_j-\Pi_j(E_{\mathbf p^{k}}[I(h(\mathbf X)\leq y_j)]))^2]+C((\beta^k)^2+\frac{(\gamma^k)^2}{\beta^k}).
\end{align*}
Denote by $D_j^k=E[( s^{k}_j-\Pi_j(E_{\mathbf p^{k}}[I(h(\mathbf X)\leq y_j)]))^2]$. When the sequences $\gamma^k$ and $\beta^k$ are taken to be \eqref{stepsize}, the recursion reduces to
\begin{align*}
D_j^{k+1}\leq (1-\frac{2b}{k^{\alpha_2}})D_j^k+C(\frac{1}{k^{2\alpha_2}}+\frac{1}{k^{2\alpha_1-\alpha_2}})
\end{align*}
which by Lemma \ref{rate:slack} implies that $D_j^k=O(\frac{1}{k^{\alpha_2}}+\frac{1}{k^{2(\alpha_1-\alpha_2)}})$. Therefore from \eqref{first term bd} we conclude
\begin{equation}\label{first bd}
E[\Vert E[\hat{\bm\phi}_{\mathbf p}^k\vert \mathcal F^k]-\bm\phi(\mathbf p^k))\Vert^2]\leq Cmn\sum_{j=1}^nD_j^k=O(\frac{1}{k^{\alpha_2}}+\frac{1}{k^{2(\alpha_1-\alpha_2)}}).
\end{equation}

To bound the term $E[\|\lambda^k\hat{\bm\Psi}^k+\hat{\bm\phi}_{\mathbf p}^k\|_\infty^2|\mathcal F^k]$, we use Minkowski inequality to get
\begin{eqnarray*}
E[\|\hat{\bm\phi}_{\mathbf p}^k\|_\infty^2|\mathcal F^k]&\leq& E\brac{\sum_{i=1}^m\prth{\hat{\phi}_{\mathbf p,i}^k}^2\bigg\vert\mathcal F^k}\\
&\leq&4n\sum_{i=1}^m\sum_{j=1}^n(2+q_{1-\alpha}/\sqrt n)^2E_{\mathbf p^k}[\big(I(h(\mathbf X)\leq y_j)S_i(\mathbf X;\mathbf p^k)\big)^2]
\end{eqnarray*}
and
\begin{equation*}
E[\|\hat{\bm\Psi}^k\|_\infty^2|\mathcal F^k]\leq E\brac{\sum_{i=1}^m\prth{\hat{\Psi}_i^k}^2\bigg\vert\mathcal F^k}\leq4m\sum_{i=1}^m\sum_{j=1}^nE_{\mathbf p^k}[\big(g(\mathbf X)S_i(\mathbf X;\mathbf p^k)\big)^2].
\end{equation*}
Again by Proposition \ref{lemma:derivative continuity}, each expectation in the sum is continuous in $\mathbf p^k$, hence uniformly bounded in $\mathcal P(\epsilon)$ by compactness. Therefore $E[\|\hat{\bm\phi}_{\mathbf p}^k\|^2_\infty|\mathcal F^k]\leq C$ and $E[\|\hat{\bm\Psi}^k\|^2_\infty|\mathcal F^k]\leq C$ uniformly holds for some $C>0$. This implies
\begin{equation}\label{second bd}
E[\|\lambda^k\hat{\bm\Psi}^k+\hat{\bm\phi}_{\mathbf p}^k\|_\infty^2|\mathcal F^k]\leq 2(E[\|\hat{\bm\phi}_{\mathbf p}^k\|_\infty^2|\mathcal F^k]+(\lambda^k)^2E[\|\hat{\bm\Psi}^k\|_\infty^2|\mathcal F^k])\leq C.
\end{equation}

Assumption \ref{cond:penalized} entails $\gamma^k(\lambda^k\bm\Psi(\mathbf p^k)+\bm\phi(\mathbf p^k))'(\mathbf p_{\epsilon}^*(\lambda^k)-\mathbf p^k)\leq0$. Substituting \eqref{first bd} and \eqref{second bd} into \eqref{ineq:convergence} we arrive at
\begin{align*}
&\sum_{k=1}^\infty E[E[V(\mathbf p^{k+1},\mathbf p_{\epsilon}^*(\lambda^{k+1}))-V(\mathbf p^k,\mathbf p_{\epsilon}^*(\lambda^k))|\mathcal F^k]^+]\\
\leq&\sum_{k=1}^\infty O(\frac{1}{k^{\alpha_1+\frac{1}{2}\alpha_2}}+\frac{1}{k^{2\alpha_1-\alpha_2}})+\sum_{k=1}^\infty\frac{1}{2}(\gamma^k)^2E[\|\lambda^k\hat{\bm\Psi}^k+\hat{\bm\phi}_{\mathbf p}^k\|_\infty^2|\mathcal F^k]+\sum_{k=1}^\infty C\lvert\log \epsilon\rvert \Vert\mathbf p_{\epsilon}^*(\lambda^{k+1})-\mathbf p_{\epsilon}^*(\lambda^k)\Vert\\
\leq&C\sum_{k=1}^\infty\big(\frac{1}{k^{\alpha_1+\frac{1}{2}\alpha_2}}+\frac{1}{k^{2\alpha_1-\alpha_2}}+\frac{1}{k^{2\alpha_1}}+\Vert\mathbf p_{\epsilon}^*(\lambda^{k+1})-\mathbf p_{\epsilon}^*(\lambda^k)\Vert\big)<\infty.
\end{align*}
By martingale convergence theorem (Corollary in Section 3 in \cite{blum1954multidimensional}, restated in Theorem \ref{martingale convergence} in the Appendix), we have $V(\mathbf p^k,\mathbf p_{\epsilon}^*(\lambda^k))$ converges a.s. to some random variable $V_\infty$. Because of $\mathbf p_{\epsilon}^*(\lambda^k)\to \mathbf p^*_{\epsilon}\in \mathcal P(\epsilon)$ and inequality \eqref{ktok+1} which holds uniformly for $\mathbf p^{k+1}\in \mathcal P(\epsilon)$, we conclude that $V(\mathbf p^k,\mathbf p^*_{\epsilon})$ converges a.s. to the same variable $V_\infty$.

Now we would like to argue that the limit $V_{\infty}=0$ a.s.. To this end it suffices to show that a.s. there exists a subsequence of $\mathbf p^{k}$ converging to $\mathbf p^*_{\epsilon}$. Taking expectation and summing up on both sides of \eqref{interim1} and using similar bounding techniques, we have
\begin{align*}
&\sum_{k=1}^\infty E[\gamma^k(\lambda^k\bm\Psi(\mathbf p^k)+\bm\phi(\mathbf p^k))'(\mathbf p^k-\mathbf p_{\epsilon}^*(\lambda^k))]\\
\leq& V(\mathbf p^1,\mathbf p_{\epsilon}^*(\lambda^1))+C\sum_{k=1}^\infty\big(\frac{\gamma^k}{\sqrt{M_1^k}}+(\gamma^k)^2+\Vert\mathbf p_{\epsilon}^*(\lambda^{k+1})-\mathbf p_{\epsilon}^*(\lambda^k)\Vert\big)<\infty.
\end{align*}
Since each $(\lambda^k\bm\Psi(\mathbf p^k)+\bm\phi(\mathbf p^k))'(\mathbf p^k-\mathbf p_{\epsilon}^*(\lambda^k))\geq0$, it follows that
\begin{align*}
\sum_{k=1}^\infty \gamma^k(\lambda^k\bm\Psi(\mathbf p^k)+\bm\phi(\mathbf p^k))'(\mathbf p^k-\mathbf p_{\epsilon}^*(\lambda^k))<\infty\;a.s..
\end{align*}
Define the (random) set of feasible-solution indices
\begin{align*}
\mathcal K_1=\{k\geq 1:\mathbf p^k\text{ is feasible for }\eqref{KS RO min discretized p epsilon}\}.
\end{align*}
Note that when $\mathbf p^k$ is feasible for \eqref{KS RO min discretized p epsilon}, it holds $\bm\phi(\mathbf p^k)=\mathbf 0$, hence
\begin{align}
&\sum_{k\in \mathcal K_1} \gamma^k\lambda^k\bm\Psi(\mathbf p^k)'(\mathbf p^k-\mathbf p_{\epsilon}^*(\lambda^k))<\infty,\;a.s.\label{feasiblesum}\\
&\sum_{k\notin \mathcal K_1} \gamma^k(\lambda^k\bm\Psi(\mathbf p^k)+\bm\phi(\mathbf p^k))'(\mathbf p^k-\mathbf p_{\epsilon}^*(\lambda^k))<\infty,\;a.s.\label{infeasiblesum}
\end{align}

If $\sum_{k\in \mathcal K_1}\gamma^k\lambda^k=\infty$, then due to \eqref{feasiblesum} there must exist a subsequence $k_i\in \mathcal K_1$ such that $\bm\Psi(\mathbf p^{k_i})'(\mathbf p^{k_i}-\mathbf p_{\epsilon}^*(\lambda^{k_i}))\to0$. Since $\mathbf p_{\epsilon}^*(\lambda^k)\to \mathbf p^*_{\epsilon}$, this implies that $\bm\Psi(\mathbf p^{k_i})'(\mathbf p^{k_i}-\mathbf p^*_{\epsilon})\to0$, which by Assumption \ref{cond:constrained} further implies that $\mathbf p^{k_i}\to \mathbf p^*_{\epsilon}$.

Otherwise if $\sum_{k\in \mathcal K_1}\gamma^k\lambda^k<\infty$ then it must hold $\sum_{k\notin \mathcal K_1}\gamma^k\lambda^k=\infty$ because the parameters stated in the theorem satisfy $\sum_{k=1}^{\infty}\gamma^k\lambda^k=\infty$. Due to \eqref{infeasiblesum} there exists a subsequence $k_i\notin \mathcal K_1$ such that
\begin{equation}\label{to be contradicted}
(\bm\Psi(\mathbf p^{k_i})+\frac{1}{\lambda^{k_i}}\bm\phi(\mathbf p^{k_i}))'(\mathbf p^{k_i}-\mathbf p^*_{\epsilon}(\lambda^{k_i}))\to 0.
\end{equation}
By a compactness argument, there exists a subsubsequence $k_i'\notin \mathcal K_1$ such that $\mathbf p^{k_i'}$ converges to some $\mathbf q\in \mathcal P(\epsilon)$. First we argue that $\mathbf q$ must be feasible for \eqref{KS RO min discretized p epsilon}. Since $\lambda^{k}\to0$ and $\bm\Psi(\mathbf p^{k_i'}),\bm\phi(\mathbf p^{k_i'})$ are uniformly bounded, it is clear that $(\lambda^{k_i'}\bm\Psi(\mathbf p^{k_i'})+\bm\phi(\mathbf p^{k_i'}))'(\mathbf p^{k_i'}-\mathbf p^*_{\epsilon}(\lambda^{k_i'}))\to 0$ and $\lambda^{k_i'}\bm\Psi(\mathbf p^{k_i'})'(\mathbf p^{k_i'}-\mathbf p^*_{\epsilon}(\lambda^{k_i'}))\to 0$ hold. Therefore the difference $\bm\phi(\mathbf p^{k_i'})'(\mathbf p^{k_i'}-\mathbf p^*_{\epsilon}(\lambda^{k_i'}))\to 0$. On the other hand $\bm\phi(\mathbf p^{k_i'})'(\mathbf p^{k_i'}-\mathbf p^*_{\epsilon}(\lambda^{k_i'}))\to \bm\phi(\mathbf q)'(\mathbf q-\mathbf p^*_{\epsilon})$ because $\bm\phi(\cdot)$ is continuous. This means $\bm\phi(\mathbf q)'(\mathbf q-\mathbf p^*_{\epsilon})=0$ so $\mathbf q$ must be feasible in view of Assumption \ref{cond:constrained}. Then we argue $\mathbf q=\mathbf p^*_{\epsilon}$ in fact. If $\mathbf q\neq \mathbf p^*_{\epsilon}$ then $\bm\Psi(\mathbf q)'(\mathbf q-\mathbf p^*_{\epsilon})>0$ by Assumption \ref{cond:constrained}, and we derive a contradiction as follows. Recall that each $\mathbf p^{k_i'}$ is infeasible for \eqref{KS RO min discretized p epsilon} and $\bm\phi(\mathbf p^{k_i'})\to \bm\phi(\mathbf q)=0$, where $\bm\phi(\mathbf q)$ vanishes since $\mathbf q$ is feasible. We have
\begin{align*}
&\liminf_i (\bm\Psi(\mathbf p^{k_i'})+\frac{1}{\lambda^{k_i'}}\bm\phi(\mathbf p^{k_i'}))'(\mathbf p^{k_i'}-\mathbf p^*_{\epsilon}(\lambda^{k_i'}))\\
=&\liminf_i \big\{\bm\Psi(\mathbf p^{k_i'})'(\mathbf p^{k_i'}-\mathbf p^*_{\epsilon}(\lambda^{k_i'}))+\frac{1}{\lambda^{k_i'}}\bm\phi(\mathbf p^{k_i'})'(\mathbf p^{k_i'}-\mathbf p^*_{\epsilon})+\frac{1}{\lambda^{k_i'}}\bm\phi(\mathbf p^{k_i'})'(\mathbf p^*_{\epsilon}-\mathbf p^*_{\epsilon}(\lambda^{k_i'}))\big\}\\
\geq& \liminf_i\bm\Psi(\mathbf p^{k_i'})'(\mathbf p^{k_i'}-\mathbf p^*_{\epsilon}(\lambda^{k_i'}))+\liminf_i\frac{1}{\lambda^{k_i'}}\bm\phi(\mathbf p^{k_i'})'(\mathbf p^{k_i'}-\mathbf p^*_{\epsilon})+\liminf_i\frac{1}{\lambda^{k_i'}}\bm\phi(\mathbf p^{k_i'})'(\mathbf p^*_{\epsilon}-\mathbf p^*_{\epsilon}(\lambda^{k_i'}))\\
\geq &\bm\Psi(\mathbf q)'(\mathbf q-\mathbf p^*_{\epsilon})+0+\liminf_i\frac{1}{\lambda^{k_i'}}o(1)O(\lambda^{k_i'})=\bm\Psi(\mathbf q)'(\mathbf q-\mathbf p^*_{\epsilon})>0
\end{align*}
which contradicts \eqref{to be contradicted}.

The above argument shows that a.s.~there exists a subsequence of $\mathbf p^{k}$ converging to $\mathbf p^*_{\epsilon}$, hence the corresponding $V(\mathbf p^{k},\mathbf p^*_{\epsilon})\to0$. Since we have proved above that $V(\mathbf p^k,\mathbf p^*_{\epsilon})$ converges a.s., the limit must be identically 0. Therefore, by Pinsker's inequality, we have $\mathbf p^k\to\mathbf p^*_{\epsilon}$ in total variation a.s.. This concludes the theorem.\Halmos\endproof


As discussed at the end of Section \ref{sec:procedure}, our results and algorithms still hold in the presence of a collection of auxiliary independent input processes $\mathbf W$ distributed according to known distributions. Like in Sections \ref{sec:KS} and \ref{sec:KS proofs}, all proofs in this section still apply by invoking the same conditioning argument. Specifically, in Proposition \ref{lemma:derivative1} the expressions \eqref{score function1},\eqref{score function2},\eqref{score function3} are still valid with $h(\mathbf X),g(\mathbf X)$ replaced by $E_{P_W^0}[h(\mathbf X,\mathbf W)|\mathbf X],E_{P_W^0}[g(\mathbf X,\mathbf W)|\mathbf X]$, so are the estimators \eqref{gradient estimator 1},\eqref{gradient estimator 3}. In Lemma \ref{lemma:derivative continuity}, the continuity of moments of gradient estimators can be similarly established by conditioning. For example, the moment $E_{\mathbf p}\left[( g(\mathbf X,\mathbf W)S_i(\mathbf X;\mathbf p))^2\right]$ is equal to $E_{\mathbf p}\big[E_{P_{W}^0}[g^2(\mathbf X,\mathbf W)\vert \mathbf X]S_i^2(\mathbf X;\mathbf p)\big]$, hence the same proof applies viewing $E_{P_{W}^0}[g^2(\mathbf X,\mathbf W)\vert \mathbf X]$ as the performance measure. Similarly, the boundedness condition in Theorem \ref{main guarantee epsilon} is made on $E_{P_{W}^0}[g(\mathbf X,\mathbf W)\vert \mathbf X]$ instead.

\section{An Alternate MDSA Algorithm and Some Further Discussion}\label{sec:proofs algorithm}
Algorithm \ref{algo1} shows an alternate MDSA algorithm that does not use slack variables, but at the expense of increasing the simulation replication size per iteration.

When applied to the (restricted) penalized minimization problem \eqref{KS RO min discretized penalty2}, MDSA solves the following optimization given a current iterate $\mathbf p^k$
\begin{equation}
\begin{array}{ll}
\min&\gamma^k(\lambda\hat{\bm\Psi}^k+\hat{\bm\phi}^k)'(\mathbf p-\mathbf p^k)+V(\mathbf p^k,\mathbf p)\\
\text{subject to}&\mathbf p\in\mathcal P(\epsilon)
\end{array}\label{step optimization1}
\end{equation}
where $\hat{\bm\Psi}^k$ carries the gradient information of the target performance measure $\psi$ at $\mathbf p^k$, $\hat{\bm\phi}^k$ contains the gradient information of the quadratic penalty function in \eqref{KS RO min discretized penalty2} at $\mathbf p^k$, and $V(\cdot,\cdot)$ is the KL divergence defined in \eqref{KL}. The step-wise subproblem \eqref{step optimization1} without stochastic noise is also called the entropic descent algorithm (\cite{beck2003mirror}). To make it a single-run procedure, we decrease the penalty coefficient $\lambda$ as the iteration goes on, and thereby arrive at the following counterpart of \eqref{step optimization2 lambda}
\begin{equation}
\begin{array}{ll}
\min&\gamma^k(\lambda^k\hat{\bm\Psi}^k+\hat{\bm\phi}^k)'(\mathbf p-\mathbf p^k)+V(\mathbf p^k,\mathbf p)\\
\text{subject to}&\mathbf p\in\mathcal P(\epsilon)
\end{array}\label{step optimization1 lambda}
\end{equation}

Inspired by \eqref{score function2} in Proposition \ref{lemma:derivative1}, we use the following estimator for the gradient of the penalty function $\bm\phi(\mathbf p)=(\phi_i(\mathbf p))_{i=1}^m$
\begin{equation}
\hat{\phi}_i(\mathbf p)=2\sum_{j=1}^n(u_j-\Pi_j(u_j))\frac{1}{M_2}\sum_{r=1}^{M_2}I(h(\tilde{\mathbf X}^{(r)})\leq y_j)S_i(\tilde{\mathbf X}^{(r)};\mathbf p),\ u_j=\frac{1}{M_1}\sum_{r=1}^{M_1}I(h(\mathbf X^{(r)})\leq y_j)\label{gradient estimator 2}\\
\end{equation}
where $\mathbf X^{(r)}$ and $\tilde{\mathbf X}^{(r)}$ are independent copies of the i.i.d. input process generated under $\mathbf p$ and are used simultaneously for all $i,j$. Since we are using the plug-in estimator $\Pi_j(u_j)$ for the projection, in general \eqref{gradient estimator 2} has a bias. In particular, the bias can be shown to vanish as slow as $O(1/\sqrt{M_1})$ if $E_{\mathbf p}[I(h(\mathbf X)\leq y_j)]$ is close to either $\hat F_Y(y_j+)-q_{1-\alpha}/\sqrt n$ or $\hat F_Y(y_j-)+q_{1-\alpha}/\sqrt n$. Due to this biasedness, the batch size $M_1$ has to grow to $\infty$ in the course of iteration in order for the algorithm to converge properly.


\begin{algorithm}[h]
  \caption{Alternate MDSA for solving \eqref{KS RO min discretized penalty2}}
  \textbf{Input: }A small parameter $\epsilon>0$, initial solution $\mathbf p^1\in\mathcal P(\epsilon)=\{\mathbf p:\sum_{i=1}^mp_i=1,p_i\geq\epsilon\text{\ for\ }i=1,\ldots,m\}$, a step size sequence $\gamma^k$, a penalty sequence $\lambda^k$, a sample size sequences $M^k_1$, and sample sizes $M_2,M_3$.

  \textbf{Iteration: }For $k=1,2,\ldots$, do the following: Given $\mathbf p^k$,
\begin{algorithmic}

\State \textbf{1.} Estimate the probabilities $E_{\mathbf p^k}[I(h(\mathbf X)\leq y_j)],j=1,\ldots,n$ with
$$u^k_j=\frac{1}{M_1^k}\sum_{r=1}^{M_1^k}I(h(\mathbf X^{(r)})\leq y_j)$$
where $\mathbf X^{(r)}$ are $M_1^k$ independent copies of the input process generated under $\mathbf p^k$.

\State \textbf{2.} Estimate $\hat{\bm\phi}^k=(\hat\phi_1^k,\ldots,\hat\phi_m^k)$, the gradient of the penalty term, with
$$\hat\phi_i^k=2\sum_{j=1}^n(u^k_j-\Pi_j(u^k_j))\frac{1}{M_2}\sum_{r=1}^{M_2}I(h(\tilde{\mathbf X}^{(r)})\leq y_j)S_i(\tilde{\mathbf X}^{(r)};\mathbf p^k)$$
where $\mathbf X^{(r)}$ are the same set of replications used in Step 1, and $\tilde{\mathbf X}^{(r)}$ are another $M_2$ independent copies of the input process generated under $\mathbf p^k$.

\State \textbf{3.} Estimate $\hat{\bm\Psi}^k=(\hat\Psi_1^k,\ldots,\hat\Psi_m^k)$, the gradient of $E_{\mathbf p}[g(\mathbf X)]$, with
$$\hat\Psi_i^k=\frac{1}{M_3}\sum_{r=1}^{M_3}g(\tilde{\tilde{\mathbf X}}^{(r)})S_i(\tilde{\tilde{\mathbf X}}^{(r)};\mathbf p^k)$$
where $\tilde{\tilde{\mathbf X}}^{(r)}$ are another $M_3$ independent copies of the input process generated under $\mathbf p^k$.

\State \textbf{4.} Compute
$\mathbf p^{k+1}=(p_1^{k+1},\ldots,p_m^{k+1})$ by running Algorithm \ref{sort_search} with $p_i=p_i^k$ and $\xi_i=\gamma^k(\lambda^k\hat\Psi_i^k+\hat\phi_{i}^k)$.



\end{algorithmic}\label{algo1}
\end{algorithm}

Like for Algorithm \ref{algo2}, the following provides the convergence guarantee of Algorithm \ref{algo1}:
\begin{theorem}
Under Assumptions \ref{cond:constrained}, \ref{cond:penalized} and \ref{cond:lambda}, if the step size sequence $\{\gamma^k\}$, the penalty sequence $\{\lambda^k\}$ and the sample size sequence $\{M_1^k\}$ of Algorithm \ref{algo1} are chosen such that
$$\sum_{k=1}^\infty\gamma^k\lambda^k=\infty,\ \ \sum_{k=1}^\infty(\gamma^k)^2<\infty,\ \ \sum_{k=1}^\infty\frac{\gamma^k}{\sqrt{M_1^k}}<\infty,\ \ \lambda^k\to 0\text{ and non-increasing}$$
then $\mathbf p^k$ generated in Algorithm \ref{algo1} converges to $\mathbf p^*_{\epsilon}$ a.s.. In particular, when the sequences are chosen as
\begin{equation*}
\begin{aligned}
&\gamma^k=\frac{a}{k^{\alpha_1}},\ \frac{1}{2}<\alpha_1\leq 1\\
&M_1^k=bk^{\alpha_2},\ \alpha_2>2(1-\alpha_1)\\
&\lambda^k=
\begin{cases}
\frac{c}{k^{\alpha_3}},\ 0<\alpha_3\leq 1-\alpha_1&\text{if }\frac{1}{2}<\alpha_1<1\\
\frac{c}{\log k}&\text{if }\alpha_1=1
\end{cases}
\end{aligned}
\end{equation*}
$\mathbf p^k$ converges to $\mathbf p^*_{\epsilon}$ a.s..\label{thm:algo1}
\end{theorem}
Here are some discussions on the parameter choices of Algorithm \ref{algo1}. $\sum_{k=1}^\infty(\gamma^k)^2<\infty$ is a standard condition in SA which ensures that the effect of stochasticity will vanish eventually, whereas the condition $\sum_{k=1}^\infty\gamma^k/\sqrt{M_1^k}<\infty$ is meant to eliminate the effect of biasedness of the gradient estimator \eqref{gradient estimator 2}. What is special about our MDSA is the condition $\sum_{k=1}^\infty \gamma^k\lambda^k=\infty$. The rationale for this condition is as follows. When $\mathbf p$ is feasible for $\eqref{KS RO min discretized p epsilon}$, the gradient of the penalty function vanishes, i.e.~$\bm\phi(\mathbf p)=\mathbf 0$, hence the effective step size in \eqref{step optimization1 lambda} is $\gamma^k\lambda^k$. Under the condition $\sum_{k=1}^\infty \gamma^k\lambda^k=\infty$, the algorithm is able to fully explore the feasible set of \eqref{KS RO min discretized p epsilon}.

The difference between Algorithm \ref{algo2} and \ref{algo1} lies in how the projection $\Pi_j(E_{\mathbf p^k}[I(h(\mathbf X)\leq y_j)])$ at the current iterate $\mathbf p^k$ is estimated. Algorithm \ref{algo1} computes the projection by directly simulating $E_{\mathbf p^k}[I(h(\mathbf X)\leq y_j)]$ from scratch and substituting into the projection $\Pi_j$ in each iteration, whereas Algorithm \ref{algo2} iteratively updates the slack variables $s_j^k$ together with the decision variable in such a way that eventually each $s_j^k$ consistently estimates the projection $\Pi_j(E_{\mathbf p^k}[I(h(\mathbf X)\leq y_j)])$.

We point out that both Algorithm \ref{algo2} and \ref{algo1} are essentially solving the formulation \eqref{KS RO min discretized penalty2}, despite the fact that the design of Algorithm \ref{algo2} is mostly based on \eqref{KS RO min discretized penalty}. The reason that neither of Algorithm \ref{algo2} and \ref{algo1} solves the formulation \eqref{KS RO min discretized penalty} has to do with the fact that algorithmically the formulation \eqref{KS RO min discretized penalty} with slack variables in general is not as well behaved as the formulation \eqref{KS RO min discretized penalty2} with the projections, despite their mathematical equivalence. To see this, consider a generic inequality constraint $f(x)\leq 0$ where $x$ is some decision variable. It is easy to see that the quadratic penalty $(\max\{f(x),0\})^2$ expressed via projection preserves the convexity of $f(x)$, whereas the one with slack variable $s\leq 0$, $(f(x)-s)^2$, can very likely lose convexity even if $f(x)$ itself is convex. In fact, if $(f(x)-s)^2$ is jointly convex in $x$ and $s$, $(\max\{f(x),0\})^2$ is guaranteed to be convex. This also explains why the general convexity criterion in Assumptions \ref{cond:constrained} and \ref{cond:penalized} is imposed on formulation \eqref{KS RO min discretized penalty2}.

\proof{Proof of Theorem \ref{thm:algo1}.}
The proof resembles that of Theorem \ref{thm:algo2}. Let $\mathcal F^k$ be the filtration generated by $\{\mathbf p^1,\ldots,\mathbf p^k\}$. Following the same line of argument, we have the following counterpart of \eqref{interim2}
\begin{align}
\nonumber&E[V(\mathbf p^{k+1},\mathbf p_{\epsilon}^*(\lambda^{k+1}))-V(\mathbf p^k,\mathbf p_{\epsilon}^*(\lambda^k))|\mathcal F^k]\\
\nonumber\leq& \gamma^k(\lambda^k\bm\Psi(\mathbf p^k)+\bm\phi(\mathbf p^k))'(\mathbf p_{\epsilon}^*(\lambda^k)-\mathbf p^k)+\gamma^k(E[\hat{\bm\phi}^k\vert \mathcal F^k]-\bm\phi(\mathbf p^k))'(\mathbf p_{\epsilon}^*(\lambda^k)-\mathbf p^k)\\
&+\frac{1}{2}(\gamma^k)^2E[\|\lambda^k\hat{\bm\Psi}^k+\hat{\bm\phi}^k\|_\infty^2|\mathcal F^k]+C\lvert\log \epsilon\rvert \Vert\mathbf p_{\epsilon}^*(\lambda^{k+1})-\mathbf p_{\epsilon}^*(\lambda^k)\Vert.\label{interim1}
\end{align}

We need to bound $E[\hat{\bm\phi}^k\vert \mathcal F^k]-\bm\phi(\mathbf p^k)$ and $E[\|\lambda^k\hat{\bm\Psi}^k+\hat{\bm\phi}^k\|_\infty^2|\mathcal F^k]$. By independence of $\mathbf X^{(r)}$ and $\tilde{\mathbf X}^{(r)}$ and conditional Jensen's inequality
\begin{align*}
\lvert E[\hat{\phi}^k_i\big\vert \mathcal F^k]-\phi_i(\mathbf p^k)\rvert&=2\big\lvert\sum_{j=1}^n(\Pi_j(E_{\mathbf p^k}[I(h(\mathbf X)\leq y_j)])-E[\Pi_j(u^k_j)\big\vert \mathcal F^k]) E_{\mathbf p^k}[I(h(\mathbf X)\leq y_j)S_i(\mathbf X;\mathbf p^k)]\big\rvert\\
&\leq C\sum_{j=1}^n\lvert\Pi_j(E_{\mathbf p^k}[I(h(\mathbf X)\leq y_j)])-E[\Pi_j(u^k_j)\big\vert \mathcal F^k]\rvert\\
&\leq C\sum_{j=1}^n E[\lvert\Pi_j(E_{\mathbf p^k}[I(h(\mathbf X)\leq y_j)])-\Pi_j(u^k_j)\rvert\big\vert \mathcal F^k]\\
&\leq C\sum_{j=1}^n E[\lvert E_{\mathbf p^k}[I(h(\mathbf X)\leq y_j)]-u^k_j\rvert\big\vert \mathcal F^k]\\
&\leq C\sum_{j=1}^n \sqrt{E[(E_{\mathbf p^k}[I(h(\mathbf X)\leq y_j)]-u^k_j)^2\big\vert \mathcal F^k]}=O(\frac{1}{\sqrt{M_1^k}})
\end{align*}
where in the second last inequality we use the contraction property of projection, i.e.~$\lvert \Pi_j(a)-\Pi_j(b)\rvert\leq \lvert a-b\rvert$ for any $a,b\in \R$. The first inequality holds because each derivative $E_{\mathbf p^k}[I(h(\mathbf X)\leq y_j)S_i(\mathbf X;\mathbf p^k)]$ by Proposition \ref{lemma:derivative continuity} is continuous in $\mathbf p$ and by a compactness argument is hence uniformly bounded in $\mathcal P(\epsilon)$. Following the proof of Theorem \ref{thm:algo2}, one can show that
\begin{equation*}
E[\|\lambda^k\hat{\bm\Psi}^k+\hat{\bm\phi}^k\|_\infty^2|\mathcal F^k]\leq C
\end{equation*}
as the counterpart of \eqref{second bd}.

Therefore, taking expectation and summing up on both sides of \eqref{interim1}, we have
\begin{align*}
&\sum_{k=1}^\infty E[E[V(\mathbf p^{k+1},\mathbf p_{\epsilon}^*(\lambda^{k+1}))-V(\mathbf p^k,\mathbf p_{\epsilon}^*(\lambda^k))|\mathcal F^k]^+]\\
\leq&\sum_{k=1}^\infty O(\frac{\gamma^k}{\sqrt{M_1^k}})+\sum_{k=1}^\infty\frac{1}{2}(\gamma^k)^2E[\|\lambda^k\hat{\bm\Psi}^k+\hat{\bm\phi}^k\|_\infty^2|\mathcal F^k]+\sum_{k=1}^\infty C\lvert\log \epsilon\rvert \Vert\mathbf p_{\epsilon}^*(\lambda^{k+1})-\mathbf p_{\epsilon}^*(\lambda^k)\Vert\\
\leq&C\sum_{k=1}^\infty\big(\frac{\gamma^k}{\sqrt{M_1^k}}+(\gamma^k)^2+\Vert\mathbf p_{\epsilon}^*(\lambda^{k+1})-\mathbf p_{\epsilon}^*(\lambda^k)\Vert\big)<\infty.
\end{align*}
The rest of the proof is the same as that of Theorem \ref{thm:algo2}.\Halmos\endproof

\section{A Randomized Stochastic Projected Gradient Algorithm for the Comparison in Section \ref{sec:numerics}}\label{sec:mini-batch}

We show a randomized stochastic projected gradient (RSPG) algorithm that we compare with in the numerical section. Algorithm \ref{RSPG} shows the procedure for a single run. Algorithm \ref{2-RSPG-V} includes a post-optimization step to boost its performance. As a rough guidance, we use $\bar\gamma<1/L$ where $L$ is the Lipschitz constant of the gradient function, and $M=O(Nm)$ (the $m$ here could possibly be removed), where $m$ is the dimension of the decision space. $S$ could be a small number like $5,10$, and the post-optimization batch size $M'$ is chosen to be some big number. The penalty $\lambda$ is chosen small and fixed.


\begin{algorithm}[h]
  \caption{Randomized stochastic projected gradient (RSPG) for solving \eqref{KS RO min discretized penalty}}
  \textbf{Input: }A small parameter $\epsilon>0$, initial solution $\mathbf p^1\in\mathcal P(\epsilon)=\{\mathbf p:\sum_{i=1}^mp_i=1,p_i\geq\epsilon\text{\ for\ }i=1,\ldots,m\}$ and $\mathbf s^1\in [\hat F_Y(y_1+)-\frac{q_{1-\alpha}}{\sqrt n},\hat F_Y(y_1-)+\frac{q_{1-\alpha}}{\sqrt n}]\times \cdots\times [\hat F_Y(y_n+)-\frac{q_{1-\alpha}}{\sqrt n},\hat F_Y(y_n+)-\frac{q_{1-\alpha}}{\sqrt n}]$, step size $\bar\gamma$ for both $\mathbf p$ and $\mathbf s$, penalty $\lambda$, batch size $M$, and number of iterations $N$.

\textbf{Generate random stopping time: }Draw $\tau$ uniformly from $\{1,\ldots,N\}$

  \textbf{Iteration: }For $k=1,\ldots,\tau-1$ do the following: Given $\mathbf p^k,\mathbf s^k$,
\begin{algorithmic}


\State \textbf{1.} Estimate $\hat{\bm\phi}_{\mathbf p}^k=(\hat\phi_{\mathbf p,1}^k,\ldots,\hat\phi_{\mathbf p,m}^k)$, the gradient of the penalty term with respect to $\mathbf p$, with
$$\hat\phi_{\mathbf p,i}^k=2\sum_{j=1}^n\frac{1}{M}\sum_{r=1}^{M}(I(h(\mathbf X^{(r)})\leq y_j)-s^k_j)\frac{1}{M}\sum_{r=1}^{M}I(h(\tilde{\mathbf X}^{(r)})\leq y_j)S_i(\tilde{\mathbf X}^{(r)};\mathbf p^k)$$
where each of $\mathbf X^{(r)},\tilde{\mathbf X}^{(r)}$ are $M$ independent copies of the input process generated under $\mathbf p^k$.

\State \textbf{2.} Estimate $\hat{\bm\Psi}^k=(\hat\Psi_1^k,\ldots,\hat\Psi_m^k)$, the gradient of $E_{\mathbf p}[g(\mathbf X)]$, with
$$\hat\Psi_i^k=\frac{1}{M}\sum_{r=1}^{M}g(\tilde{\tilde{\mathbf X}}^{(r)})S_i(\tilde{\tilde{\mathbf X}}^{(r)};\mathbf p^k)$$
where $\tilde{\tilde{\mathbf X}}^{(r)}$ are another $M$ independent copies of the input process generated under $\mathbf p^k$.

\State \textbf{3.} Estimate $\hat{\bm\phi}_{\mathbf s}^k=(\hat\phi_{\mathbf s,1}^k,\ldots,\hat\phi_{\mathbf s,n}^k)$, the gradient of the penalty term with respect to $\mathbf s$, with
$$\hat\phi_{\mathbf s,j}^k=-\frac{1}{M}\big(\sum_{r=1}^{M}(I(h(\mathbf X^{(r)})\leq y_j)-s^k_j)+\sum_{r=1}^{M}(I(h(\tilde{\mathbf X}^{(r)})\leq y_j)-s^k_j)\big)$$
where ${\mathbf X}^{(r)},\tilde{\mathbf X}^{(r)}$ are the same replications used in Step 1.

\State \textbf{4.} Compute
$\mathbf p^{k+1}=(p_1^{k+1},\ldots,p_m^{k+1})$ by running Algorithm \ref{sort_search} with $\xi_i=\bar\gamma(\lambda\hat\Psi_i^k+\hat\phi_{\mathbf p,i}^k)$
and compute $\mathbf s^{k+1}=(s_1^{k+1},\ldots,s_n^{k+1})$ by
\begin{equation*}
s_j^{k+1}=\Pi_j(s_j^k-\bar\gamma\hat\phi_{\mathbf s,j}^k)
\end{equation*}



\end{algorithmic}\label{RSPG}
\textbf{Output: }$\mathbf p^{\tau},\mathbf s^{\tau}$
\end{algorithm}

\begin{algorithm}[h]
  \caption{Two-phase RSPG for solving \eqref{KS RO min discretized penalty}}
  \textbf{Input: }A small parameter $\epsilon>0$, initial solution $\mathbf p^1\in\mathcal P(\epsilon)=\{\mathbf p:\sum_{i=1}^mp_i=1,p_i\geq\epsilon\text{\ for\ }i=1,\ldots,m\}$ and $\mathbf s^1\in [\hat F_Y(y_1+)-\frac{q_{1-\alpha}}{\sqrt n},\hat F_Y(y_1-)+\frac{q_{1-\alpha}}{\sqrt n}]\times \cdots\times [\hat F_Y(y_n+)-\frac{q_{1-\alpha}}{\sqrt n},\hat F_Y(y_n+)-\frac{q_{1-\alpha}}{\sqrt n}]$, step size $\bar\gamma$ for both $\mathbf p$ and $\mathbf s$, penalty $\lambda$, batch size $M$, number of RSPG runs $S$, and number of iterations $N$ per run. Batch size $M'$ in the post-optimization phase.

\begin{algorithmic}


\State   \textbf{1. Optimization phase: }For $s=1,\ldots,S$, run Algorithm \ref{RSPG} with initial point $\mathbf p^1,\mathbf s^1$, step size $\bar\gamma$, penalty $\lambda$, batch size $M$, and number of iterations $N$. Let $\mathbf p_s,\mathbf s_s$ be the output of the $s$-th run of Algorithm \ref{RSPG}.

\State \textbf{2. Post-optimization phase: } For $s=1,\ldots,S$, run one iteration of Step 1,2,3,4 of Algorithm \ref{RSPG} but with batch size $M'$ at $\mathbf p_s,\mathbf s_s$. Let $\mathbf p'_s,\mathbf s'_s$ be the output from Step 4 at $\mathbf p_s,\mathbf s_s$. Then compute
\begin{equation*}
(g_{\mathbf p}(\mathbf p_s,\mathbf s_s),g_{\mathbf s}(\mathbf p_s,\mathbf s_s))=(\frac{1}{\bar\gamma}(\mathbf p'_s-\mathbf p_s),\frac{1}{\bar\gamma}(\mathbf s'_s-\mathbf s_s))
\end{equation*}

\end{algorithmic}\label{2-RSPG-V}
\textbf{Output: }the $\mathbf p_{s^*},\mathbf s_{s^*}$ where $s^*=\argmin_{s}\{\Vert g_{\mathbf p}(\mathbf p_s,\mathbf s_s)\Vert_1^2+\Vert g_{\mathbf s}(\mathbf p_s,\mathbf s_s)\Vert_2^2\}$

\end{algorithm}

\section{Auxiliary Results}\label{sec:theorems}
\subsection{Results on Empirical Processes and $U$-Statistics}\label{sec:complexity}
We first introduce some definitions. Using Definition 2.1.6 in \cite{van1996weak}, given two functions $l$ and $u$, the bracket $[l,u]$ is defined as the set of all functions $f$ with $l\leq f\leq u$. An $\epsilon$-bracket is a bracket $[l,u]$ with $\|l-u\|<\epsilon$ for some norm $\|\cdot\|$. For a class of measurable functions $\mathcal F$ on $\mathcal Y\to\mathbb R$, the bracketing number $N_{[]}(\epsilon,\mathcal F,\|\cdot\|)$ is the minimum number of $\epsilon$-brackets needed to cover $\mathcal F$. Moreover, define the envelope of $\mathcal F$ as $F(\cdot)=\sup_{f\in\mathcal F}|f(\cdot)|$. 

We have the following theorem:

\begin{theorem}[Problem 3 in Chapter 2.7 of \cite{van1996weak}]
Let $\mathcal F$ be a class of measurable functions $f(\cdot,r)$ on $\mathcal Y\to\mathbb R$, indexed by $0\leq r\leq1$, such that $f(x,\cdot)$ is monotone for each $x$. If the envelope function of $\mathcal F$ is square integrable, then the bracketing number of $\mathcal F$ is polynomial.\label{bracketing}
\end{theorem}







To introduce the next theorem, we define several additional notions. For any function $f:\mathcal X^T\to\mathbb R$, and $X_1,\ldots,X_m$ generated i.i.d. from $P$, define the $U$-operator $U_T^m$ by
\begin{equation}
U_T^mf=U_T^m(f,P)=\frac{(m-T)!}{m!}\sum_{(i_1,\ldots,i_T)\in I_T^m}f(X_{i_1},\ldots,X_{i_T})\label{u operator}
\end{equation}
where $I_T^m=\{(i_1,\ldots,i_m):1\leq i_j\leq m,i_j\neq i_k\text{\ if\ }j\neq k\}$.
For convenience we denote $P^Tf=E_{P}[f]$, where $E_{P}[\cdot]$ is the expectation with respect to the $T$-fold product measure of $P$.

We say that a central theorem holds for $\{\sqrt m(U_T^mf-P^Tf)\}_{f\in\mathcal F}$ if
\begin{equation}
\{\sqrt m(U_T^mf-P^Tf)\}_{f\in\mathcal F}\Rightarrow\{\mathbb G(f)\}_{f\in\mathcal F}\text{\ in\ }\ell^\infty(\mathcal F)\label{CLT process}
\end{equation}
where $\ell^\infty(\mathcal F)$ is the space (for functionals on $\mathcal F$) defined by
$$\ell^\infty(\mathcal F)=\left\{y:\mathcal F\to\mathbb R:\sup_{f\in\mathcal F}|y(f)|<\infty\right\}$$
$\mathbb G(f)$ is a Gaussian process indexed by $\mathcal F$ that is centered and has covariance function
$$Cov(\mathbb G(f_1),\mathbb G(f_2)) = Cov(TP^{T-1}S_Tf_1,TP^{T-1}S_Tf_2)$$
where $P^{T-1}$ is defined by $P^{T-1}f(x)=\int\cdots\int f(x_1,\ldots,x_{T-1},x)\prod_{t=1}^{T-1}dP(x_t)$
and
$$S_Tf(x_1,\ldots,x_T)=\frac{1}{T!}\sum f(x_{i_1},\ldots,x_{i_T})$$
where the sum is taken over all permutations $(x_{i_1},\ldots,x_{i_T})$ of $(x_1,\ldots,x_T)$. Moreover, the process $\mathbb G(\cdot)$ is sample continuous with respect to the canonical semi-metric
$$\tau^2_{P,T}(f_1,f_2)=Var(P^{T-1}S_T(f_1-f_2))$$
where $Var(\cdot)$ is taken with respect to the probability $P$. These discussions follow from \cite{arcones1993limit}. We have ignored some measurability issues; see \cite{van1996weak} for more details.


We have the following theorem:
\begin{theorem}[Theorem 4.10 in \cite{arcones1993limit}]
Let $\mathcal F$ be a class of functions on $\mathcal X^T\to\mathbb R$. If
$$\int_0^1\sqrt{\log N_{[]}(\epsilon,\mathcal F,\|\cdot\|_{P^T,2})}d\epsilon<\infty$$
where $\|\cdot\|_{P^T,2}$ is the norm induced in the $L_2$-space under $P^T$, the $T$-fold product measure of $P$. Then the central limit theorem holds for 
$\{\sqrt m(U_T^mf-P^Tf)\}_{f\in\mathcal F}$ in the sense of \eqref{CLT process}.\label{CLT U}
\end{theorem}

From Theorem \ref{CLT U}, it is immediate that for a class of functions on $\mathcal X^T\to\mathbb R$, a bracketing number that is polynomial in $\epsilon$ implies the central limit theorem \eqref{CLT process}.

\subsection{Results Needed in the Convergence Proofs of the MDSA}
\begin{theorem}[Corollary in Section 3 in \cite{blum1954multidimensional}]
Let $Y_k$ be a sequence of integrable random variables that satisfy
$$\sum_{k=1}^\infty E[E[Y_{k+1}-Y_k|Y_1,\ldots,Y_k]^+]<\infty$$
where $x^+=x$ if $x>0$ and 0 otherwise, and are bounded below uniformly in $k$. Then $Y_k$ converges a.s. to a random variable.\label{martingale convergence}
\end{theorem}

\begin{lemma}[Adapted from Lemma 2.1 in \cite{nemirovski2009robust}]
Let $V$ be the KL divergence defined in \eqref{KL}. For every $\mathbf q\in\mathcal P$, $\mathbf p\in\mathcal P^\circ$, and $\bm\xi\in\mathbb R^m$, one has
$$V(\tilde{\mathbf p},\mathbf q)\leq V(\mathbf p,\mathbf q)+\bm\xi'(\mathbf q-\mathbf p)+\frac{\|\bm\xi\|_{\infty}^2}{2}$$
where $\tilde{\mathbf p}=\argmin_{\mathbf u\in\mathcal P}\bm\xi'(\mathbf u-\mathbf p)+V(\mathbf p,\mathbf u)$, and $\|\cdot\|_{\infty}$ is the sup norm.\label{iteration lemma}
\end{lemma}

\end{document}